\documentclass[11pt]{amsart}

\usepackage[top=1.1in,bottom=1.1in,left=1.2in,right=1.2in]{geometry}
\usepackage{amsmath, amscd, amssymb}
\usepackage[frame,cmtip,arrow,matrix,line,graph,curve]{xy}
\usepackage{graphpap, color}
\usepackage{graphicx}
\usepackage[mathscr]{eucal}
\usepackage{verbatim}

\usepackage{float}
\usepackage{booktabs}
\usepackage{multirow}
\usepackage{tabularx}
\usepackage{tikz-cd, tikz}
\usepackage{placeins}

\usepackage{hyperref}
\hypersetup{hidelinks,hypertexnames=false}

\usepackage{color}

\usepackage{marginnote}

\usepackage{lscape} 

\numberwithin{equation}{section}

\newtheorem{theorem}{Theorem}[section]
\newtheorem{corollary}[theorem]{Corollary}
\newtheorem{proposition}[theorem]{Proposition}
\newtheorem{lemma}[theorem]{Lemma}

\newtheorem{conjecture}[theorem]{Conjecture}

\theoremstyle{definition}
\newtheorem{definition}[theorem]{Definition}

\newtheorem{example}[theorem]{Example}
\newtheorem{remark}[theorem]{Remark}

\newcommand\hecke{\mathbf{H}}  

\newcommand\grs{\mathfrak{g}^{\mathrm{rs}}}
\newcommand\tgrs{\widetilde{\mathfrak{g}}^{\mathrm{rs}}}
\newcommand\rs{{\mathrm{rs}}}

\def\and{\quad{\rm and}\quad}
\def\lra{\longrightarrow }

\def\beq{\begin{equation}}
\def\eeq{\end{equation}}
\def\ben{\begin{enumerate}}
\def\een{\end{enumerate}}

\DeclareMathOperator{\Irr}{\mathrm{Irr}}
\DeclareMathOperator{\CS}{\mathrm{CS}}
\DeclareMathOperator{\Rep}{\mathrm{Rep}}
\DeclareMathOperator{\car}{\mathrm{char}}
\DeclareMathOperator{\codim}{\mathrm{codim}}
\DeclareMathOperator{\supp}{\mathrm{supp}}
\DeclareMathOperator{\Spec}{\mathrm{Spec}}

\DeclareMathOperator{\Ad}{\mathrm{Ad}}
\DeclareMathOperator{\ad}{\mathrm{ad}}

\DeclareMathOperator{\Ind}{\mathrm{Ind}}
\DeclareMathOperator{\ind}{\mathrm{ind}}
\DeclareMathOperator{\Perv}{\mathrm{Perv}}
\DeclareMathOperator{\Res}{\mathrm{Res}}

\def\sp{\mathrm{sp}}

\DeclareMathOperator{\coh}{\mathrm{H}}

\DeclareMathOperator{\image}{Im}
\DeclareMathOperator{\id}{id}
\DeclareMathOperator{\rank}{rank}
\DeclareMathOperator{\IC}{\mathrm{IC}}
\DeclareMathOperator{\ICperv}{\mathbb{IC}}
\DeclareMathOperator{\Gr}{\mathsf{Gr}}

\DeclareMathOperator{\W}{\mathsf{W}}

\DeclareMathOperator{\Hom}{Hom}

\DeclareMathOperator{\Sym}{Sym}
\DeclareMathOperator{\Semis}{Semis}
\DeclareMathOperator{\mon}{mon}
\DeclareMathOperator{\Four}{Four}

\def\cA{\mathcal {A}}

\def\cB{\mathcal {B}}
\def\cH{\mathcal {H}}

\def\Tr{\mathrm {Tr}}
\DeclareMathOperator{\Hess}{\mathsf{Hess}}

\def\A{\mathbb{A}}

\def\C{\mathbb{C}}
\def\CC{\mathbb{C}}
\def\EE{\mathbb{E}}
\def\FF{\mathbb{F}}
\def\Gm{\mathbb{G}_m}
\def\PP{\mathbb{P}}
\def\Q{\mathbb{Q}}
\def\QQ{\mathbb{Q}}

\def\ZZ{\mathbb{Z}}

\def\fg{\mathfrak {g}}
\def\fb{\mathfrak {b}}
\def\ft{\mathfrak {t}}
\def\fu{\mathfrak {u}}
\def\fX{\mathfrak {X}}
\def\pr{\mathrm{pr}}

\def\Grs{G^{\mathrm{rs}}}

\DeclareMathOperator{\IH}{\mathrm{IH}}
\DeclareMathOperator{\Ch}{\mathrm{Ch}}
\DeclareMathOperator{\Spr}{\mathrm{Spr}}
\DeclareMathOperator{\Groth}{\mathrm{Groth}}

\newlength\cellsize  

\savebox2{%
\begin{picture}(18,18)\linethickness{0.6pt} %
\put(0,0){\line(1,0){18}}
\put(0,0){\line(0,1){18}}
\put(18,0){\line(0,1){18}}
\put(0,18){\line(1,0){18}}
\end{picture}}

\newcommand\cellify[1]{\def\thearg{#1}\def\nothing{}%
\ifx\thearg\nothing\vrule width0pt height\cellsize depth0pt%
  \else\hbox to 0pt{\usebox2\hss}\fi%
  \vbox to \cellsize{\vss\hbox to \cellsize{\hss$_{#1}$\hss}\vss}}

\newcommand\tableau[1]{\vtop{\let\\=\cr
\setlength\baselineskip{-1000pt}
\setlength\lineskiplimit{1000pt}
\setlength\lineskip{0pt}
\ialign{&\cellify{##}\cr#1\crcr}}}

\title{Hessenberg character}
\date{November 4, 2026}

\author{Patrick Brosnan}
\address{Department of Mathematics\\
  University of Maryland\\
  College Park, MD USA}
\email{pbrosnan@math.umd.edu}

\author{Jaehyun Hong}
\address{Center for Complex Geometry, Institute for Basic Science (IBS), 55 Expo-ro, Yuseong-gu, Daejeon 34126, Korea}
\email{jhhong00@ibs.re.kr}

\author{Donggun Lee}
\address{June E Huh Center for Mathematical Challenges, Korea Institute for Advanced Study, 85 Hoegiro, Dongdaemun-gu, Seoul 02455, Republic of Korea}
\email{dglee@kias.re.kr}

\begin{document}

\begin{abstract}
Let $G$ be a connected reductive group with Borel subgroup $B$, flag variety $\mathcal B=G/B$, and Weyl group $W$. We define the \emph{Hessenberg functor} on the geometric Hecke category by microlocalizing along the diagonal in $\mathcal B\times\mathcal B$ and pushing forward along the Springer map. It induces the \emph{Hessenberg character} of the Iwahori--Hecke algebra, with values in graded $W$-representations. On Kazhdan--Lusztig basis elements indexed by smooth elements of $W$, this character recovers the cohomology of the associated regular semisimple Hessenberg varieties. We identify the Hessenberg functor with the microlocalization of Lusztig's character functor and compute the Hessenberg character using Lusztig's nonabelian Fourier pairing.

More generally, for families of Hessenberg schemes defined by pure-dimensional $B$-stable closed cones in the tangent space $T_{eB}\mathcal B$, we show that monodromy on the intersection cohomology of regular semisimple fibers factors through $W$. When the cone is the tangent cone of a Schubert variety, we give a local condition at $eB$ under which the same character formula computes the graded $W$-character of this intersection cohomology.

Combined with the methods of Precup and Sommers, our character formula can be used to compute the graded $W$-characters of all regular semisimple Hessenberg varieties in all Lie types of rank at most ten, including all exceptional types. These computations also motivate positivity and unimodality conjectures for Hessenberg varieties and Lusztig varieties in types $B$, $C$, and $D$, together with positivity conjectures for their partial flag versions.
\end{abstract}

\maketitle
\tableofcontents

\section{Introduction}
\renewcommand{\thetheorem}{\Alph{theorem}}


Let $G$ be a connected reductive algebraic group over an algebraically closed field $\Bbbk$. 
Let $T$ be a maximal torus in a Borel subgroup $B$ of $G$. 
Write $\fg$, $\ft$, and $\fb$ for the Lie algebras of $G$, $T$, and $B$, respectively.
Write $W$ for the Weyl group of $G$ and $\cB=G/B$ for the flag variety. 

A \emph{Hessenberg space} is a $B$-stable linear subspace $H$ of $\fg$ containing $\fb$. 
Write $X_H:=G\times^B H$ for the quotient of $G\times H$ by the action
$(g,z)\cdot b=(gb,\Ad(b^{-1})z)$. 
Consider the morphism
\beq\label{e:mainmor}
\pi_H:X_H\lra\fg, \quad[g,z]\mapsto\Ad(g)z.
\eeq
The \emph{Hessenberg scheme associated with} $H$ and $x\in\fg$ is the fiber $X_H(x):=\pi_H^{-1}(x)$. 

We may also view $X_H$ as a closed subscheme of $\fg\times\cB$ via $[g,z]\mapsto(\Ad(g)z,gB)$.
Under this identification, $\pi_H$ is the first projection, and $X_H(x)$ is the closed subscheme of $\cB$ defined by the condition $\Ad(g^{-1})x\in H$.
In particular, both $\pi_H$ and $X_H(x)$ are projective.

An element $s\in\fg$ is regular semisimple if the identity component $Z_G(s)^\circ$ of its centralizer $Z_G(s)$ is a maximal torus in $G$. 
Write $\grs$ for the open subset of such elements. We assume throughout that $\grs\neq\emptyset$; this is automatic if $\car \Bbbk\ne2$ or $G$ is adjoint
\cite[Section~13.3]{jantzen-nilpotent}.
For $G=\mathrm{SL}_2$ in characteristic $2$, however, $\grs=\emptyset$.

For $s\in\grs$, we call $X_H(s)$ a \emph{regular semisimple Hessenberg variety}. 
These fibers are smooth and projective and admit affine pavings; we recall the argument of De Mari, Procesi and Shayman \cite{demari-procesi-shayman} in Theorem~\ref{thm:DMPS}.
In particular, $\pi_H$ is smooth and projective over $\grs$.
This smooth projective family over $\grs$ gives rise to a monodromy action on $\coh^*(X_H(s))$ for each $s\in \grs$.
We use Betti cohomology with coefficients in $\EE=\C$ when
$\Bbbk\subseteq\C$, and \'etale cohomology with coefficients in
$\EE=\overline{\mathbb Q}_\ell$ for a prime $\ell$ which is invertible in $\Bbbk$, otherwise. Accordingly,
$\pi_1(\grs,s)$ denotes the topological or \'etale fundamental group respectively.

When $H=\fb$, the map $\pi_\fb:G\times^B\fb\to\fg$ is the Grothendieck simultaneous resolution.
Its fiber over $s\in\grs$ consists of the Borel subgroups containing the maximal torus $Z_G(s)^\circ$.
Its restriction to $\grs$ is a $W$-torsor (as recalled in Lemma~\ref{lem:W-torsor}), and hence induces a quotient
\[
 \pi_1(\grs,s)\twoheadrightarrow W.
\]

For every Hessenberg space $H$, the monodromy action factors through this quotient, giving a $W$-representation on each cohomology group. 
Over $\C$, this was proved in \cite{BrCh18,balibanu-crooks}, and we record that the same argument applies over $\Bbbk$ (Theorem~\ref{thm:Hess.monodromy}).
Henceforth, we define
\[
P_H(u)=\sum_{j\geq0}\coh^{2j}(X_H(s))u^j\in(\Rep W)[u],
\]
where $\Rep W$ denotes the representation ring of $W$ and $u$ is a formal variable.
Note that the odd cohomology of $X_H(s)$ vanishes, and this graded $W$-character is independent of the choice of
$s\in\grs$. Our problem is to compute $P_H(u)$.

\medskip

The $W$-action on $\coh^*(X_H(s))$ can be understood in another way, which is \emph{a priori} more amenable to direct computation.
Tymoczko \cite{tymoczko} defined an action, called the \emph{dot action}, through $Z_G(s)^\circ$-equivariant cohomology. 
This action coincides with the monodromy action above \cite{BrCh18,balibanu-crooks}; see also Remark~\ref{rem:Hess.dot.action} for an arbitrary algebraically closed field $\Bbbk$. 
The advantage of this formulation is that it can be described combinatorially using the moment graph. 
For groups of small rank, the computation can be carried out using a computer.
However, its complexity grows quickly with the rank of $G$ 
and the size of $W$. 

Over $\C$, a second approach, due to Precup and Sommers \cite{PS}, uses Fourier transform to relate $P_H(u)$ to the cohomology of Hessenberg ideal fibers. 
Choose a $G$-invariant nondegenerate symmetric bilinear form on $\fg$, and let $I:=H^\perp$, which is contained in the nilpotent radical of $\fb$. Let $\mathcal N$ be the nilpotent cone of $\fg$, and consider the projective morphism
\[
\pi_I:G\times^B I\lra\mathcal N,\quad[g,z]\mapsto\Ad(g)z.
\]
Precup and Sommers prove that every simple perverse constituent of $\pi_{I*}\underline\C_{G\times^B I}$%
\footnote{Throughout the paper, all functors are derived; for example, we write $f_*$ for the derived pushforward $Rf_*$.} is associated with a pair occurring in the Springer correspondence.
Fourier transform identifies this direct image with $\pi_{H*}\underline\C_{X_H}$, up to shift, and hence the latter records the $W$-isotypic decomposition of $\coh^*(X_H(s))$.
As a consequence, the decomposition theorem and the Springer correspondence reduce the computation of $P_H(u)$ to determining the cohomology of the fibers of $\pi_I$ as representations of the component groups of the stabilizers of the corresponding nilpotent elements.
In practice, however, determining these cohomology groups of the fibers of $\pi_I$ and the component-group actions on them can be quite
involved.

Along the way, they also prove Brosnan's conjecture: every simple perverse constituent of $\pi_{H*}\underline\C_{X_H}$ has full support. 
This conjecture arose from the same program of computing $P_H(u)$ via the fibers of $\pi_{H^\perp}$; the type $G_2$ case had previously been worked out by Xue~\cite{xue}.

\subsection{Hecke category}

One of the goals of this paper is to package this picture into a single functor defined on the Hecke category.
For a $G$-variety $X$, let $\mathscr D_G(X)$ denote the bounded derived category of $G$-equivariant constructible sheaves on $X$ with coefficients in $\EE$, and let $\Semis_G(X)$ be the full subcategory of direct sums of shifts of simple perverse sheaves. 
Set
\[
\mathsf H_G:=\Semis_G(\cB\times\cB).
\]
Equipped with convolution, $\mathsf H_G$ is the geometric Hecke category in the following sense.

Let $\hecke$ be the Iwahori--Hecke algebra of $W$ over $\cA=\ZZ[u^{\pm1/2}]$.
For $w\in W$, write $\ell(w)$ for its length and $C'_w$ for the Kazhdan--Lusztig basis element.
Let $\fX_w=G\cdot(eB,wB)$ be the corresponding diagonal $G$-orbit in $\cB\times\cB$, and let $\overline \fX_w$ be its closure.
Let
\[\mathcal F_w:=\IC_{\overline{\fX}_w}[\dim\cB],\] 
where $[m]$ denotes cohomological shift, and $\IC_X$ is the intersection cohomology complex, normalized so that $\IC_X[\dim X]$ is perverse. 
Convolution makes the split Grothendieck group $\mathscr K_\oplus(\mathsf H_G)$ into an $\cA$-algebra, with $u^{1/2}[F]=[F[-1]]$, and there is a natural isomorphism
\[
\mathscr K_\oplus(\mathsf H_G)\xrightarrow{~\sim~}\hecke,
\quad [\mathcal F_w]\mapsto\widehat C_w:=u^{\ell(w)/2}C'_w.
\]

\subsection{Hessenberg functor and Hessenberg character}
Let $\Delta\subset\cB\times\cB$ be the diagonal, whose normal bundle is the tangent bundle $T_{\cB}$.
Microlocalization along $\Delta$ consists of specialization to $T_{\cB}$ followed by Fourier transform to its dual $T^\vee_{\cB}$;
we denote this functor by $\mu_\Delta$ and recall it in Section~\ref{s:microloc}. 
Consider the moment map associated to the $G$-action on $\cB$
\[
\eta:T^\vee_{\cB}\longrightarrow\mathcal N^\vee,
\]
where $\mathcal N^\vee\subset\fg^\vee$ is the coadjoint nilpotent cone, consisting of the linear forms that vanish on the Lie algebra of some Borel subgroup of $G$.
Over $\C$, a $G$-invariant nondegenerate form identifies $\eta$ with the Springer map.
We define the \emph{Hessenberg functor} by
\[
\Hess=\eta_*\mu_\Delta:\mathsf H_G\longrightarrow\mathscr D_G(\mathcal N^\vee).
\]

Let $\Spr^\vee:=\eta_*\underline{\EE}_{T^\vee_{\cB}}[2\dim\cB]$ be the coadjoint Springer sheaf. 
It has a natural $W$-action obtained by Fourier transform from the $W$-action on the Grothendieck sheaf associated with the $W$-torsor over $\grs$ discussed above; see Subsection~\ref{ss Hecke}.
Write $\Perv_G(\mathcal N^\vee)$ for the category of $G$-equivariant perverse sheaves on $\mathcal N^\vee$.
For $A\in \Perv_G(\mathcal N^\vee)$, the Springer functor
\[
\mathsf S(A):=\Hom_{\Perv_G(\mathcal N^\vee)}(\Spr^\vee,A)
\]
gives a $W$-representation.
By taking graded perverse cohomology of the Hessenberg functor and applying $\mathsf S$, we obtain an $\cA$-linear map
\[
\widehat h:\hecke\longrightarrow\Rep W\otimes\cA.
\]
We call $\widehat h$ the \emph{Hessenberg character}.

\medskip

An element $w\in W$ is called \emph{smooth} if its Schubert variety $\overline X_w=\overline{BwB}/B\subset\cB$ is smooth.
For such an element $w$, let $H_w\subset\fg$ be the Hessenberg space determined by
\[
H_w/\fb=T_{eB}\overline X_w\subset\fg/\fb.
\]

The following theorem gives a geometric interpretation of $\widehat h(\widehat C_w)$ for smooth $w$.

\begin{theorem}[Theorem~\ref{thm:smooth}]\label{thm:smooth.intro}
For every smooth $w\in W$, we have
\[
\widehat h(\widehat C_w)=P_{H_w}(u).
\]
\end{theorem}
The term ``Hessenberg'' reflects this geometric interpretation, while the term ``character'' refers to the trace property proved below.

For smooth $w$, the normal bundle of $\Delta$ in $\overline{\fX}_w$ is $G\times^B(H_w/\fb)$.
Specialization along $\Delta$ and compatibility of Fourier transform with proper direct image identify $\Hess(\mathcal F_w)$, up to shift and Tate twist, with the Fourier transform of
$\pi_{H_w*}\underline\EE_{X_{H_w}}$.
The theorem then follows from the monodromy action and the Springer correspondence. 
It holds over any algebraically closed field $\Bbbk$, and its proof does not use the support theorem of Precup and Sommers.

\subsection{Relation to Lusztig's character}

Consider the action morphism
\[
\sigma_G:G\times\cB\lra\cB\times\cB,\quad(g,x)\mapsto(x,gx).
\]
For $w\in W$, let
$\overline Y_w=\sigma_G^{-1}(\overline{\fX}_w)$, and
let 
\[\overline\pi_w:\overline Y_w\lra G\] be the restriction of
the first projection. For $\mathbf s\in G$, its fiber is
\[
\overline Y_w(\mathbf s)= \{x\in\cB:(x,\mathbf sx)\in\overline{\fX}_w\},
\]
which we call the \emph{Lusztig scheme associated with} $w$ and $\mathbf s$. 
Write $G^{\mathrm{rs}}$ for the set of regular
semisimple elements of $G$, where regularity means that the
centralizer has dimension equal to the rank of $G$.
For $\mathbf s\in \Grs$, we call $\overline Y_w(\mathbf s)$ the \emph{regular semisimple Lusztig variety}.

In our previous work \cite{BHL25}, we proved that,
for every smooth $w\in W$, there is a degree-preserving
$W$-equivariant isomorphism
\begin{equation}\label{eq:intro.smooth.comparison}
\coh^*(\overline Y_w(\mathbf s))
\simeq \coh^*(X_{H_w}(s)),
\end{equation}
where $s\in\grs$ and $\mathbf s\in G^{\mathrm{rs}}$.
This was obtained over $\C$ using equivariant localization, and the same proof applies over $\Bbbk$ with coefficients in $\EE=\overline \Q_\ell$, using equivariant $\ell$-adic localization
\cite[Lemme~3.1]{chaudouard-laumon} in the $\ell$-adic setting;
see also Remark~\ref{rem:Hess.dot.action}. 
We also constructed a flat degeneration from $\overline Y_w(\mathbf s)$ to $X_{H_w}(s)$.
The isomorphism \eqref{eq:intro.smooth.comparison} compares the representations indexed by smooth elements of $W$, and it is natural to ask whether it lifts and extends to a functorial comparison on the entire Hecke category.
Our main result gives such a comparison.

To formulate it, 
consider the categorical correspondence
\[
\Ch=\pr_{1*}\sigma_G^*:
\mathsf H_G\longrightarrow\mathscr D_G(G),
\]
which we call \emph{Lusztig's character functor},
where $\pr_1:G\times\cB\to G$ is the first projection.
Then 
\[
\Ch(\IC_{\overline{\fX}_w})
\simeq \overline\pi_{w*}\IC_{\overline Y_w} =:\overline K_w
\]
for every $w\in W$.
By the decomposition theorem \cite{BBD}, $\overline K_w$ is semisimple.
The simple perverse sheaves occurring in its perverse cohomology, as $w$ varies, are precisely Lusztig's \emph{unipotent character sheaves}.
Lusztig's multiplicity formula expresses their graded multiplicities in $\overline K_w$ in terms of irreducible Hecke characters and his nonabelian Fourier pairing \cite{Lus85ii,Lus86v}.
We recall Lusztig's original construction of his character in \cite{Lus85}, its description in terms of $\overline K_w$, and the multiplicity formula in Section~\ref{s lusztig character}.

Let $\mu_e$ denote microlocalization at the identity of $G$ (see Subsection~\ref{ss:micro}).
We write $(m)$ for Tate twist, which is omitted in the Betti setting.
Our main result is the following.

\begin{theorem}[Theorem~\ref{thm:lus.hess} and
Corollary~\ref{cor:trace}]\label{thm:main}
There is a natural isomorphism of functors
\[
\mu_e\circ\Ch\simeq\Hess[-\dim\fb](-\dim\fb).
\]
In particular, for $A,B\in\mathsf H_G$, there is a natural
isomorphism
\[
\Hess(A\star B)\simeq\Hess(B\star A),
\]
and $\widehat h(ab)=\widehat h(ba)$ for $a,b\in\hecke$.
\end{theorem}

The trace property follows from the corresponding property of $\Ch$ \cite[Lemma~2.7.1]{trinh2021heckecategoryunipotentlocus}.

Together with Theorem~\ref{thm:smooth.intro}, the comparison shows that the assignment $\widehat C_w\mapsto P_{H_w}(u)$ for smooth $w$ is the restriction of a $\Rep W\otimes\cA$-valued trace on the entire Iwahori--Hecke algebra.

\subsection{Character formulas}

Let $\Irr W$ denote the set of irreducible characters of $W$.
Writing $\widehat h$ in this basis, 
each coefficient $\hecke\to\cA$ is a trace.
After extending scalars to $\EE(u^{1/2})$, the irreducible
characters of $\hecke$ form a basis of the space of traces.
Hence we may write
\[
 \widehat h(a)=\sum_{\chi,\psi\in\Irr W}
 c_{\chi,\psi}\widetilde\psi(a)\chi
\]
for $c_{\chi,\psi}\in\EE(u^{1/2})$,
where $\widetilde\psi$ denotes the irreducible character of
$\hecke$ corresponding to $\psi$. Thus computing
$\widehat h$ amounts to determining the coefficients
$c_{\chi,\psi}$.

Lusztig's multiplicity formula, together with Theorem~\ref{thm:main} and the classification of unipotent character sheaves with full support (Lemma~\ref{lem:fullsupp.CS}), identifies these coefficients with his nonabelian Fourier pairing:
$c_{\chi,\psi}=\{\chi,\psi\}$.
We obtain the following in Section~\ref{s lusztig character}.
\begin{corollary}\label{cor:hess.formula}
For every $w\in W$, we have
\[
\widehat h(\widehat C_w)
=\sum_{\chi,\psi\in\Irr W}
	\{\chi,\psi\}\widetilde\psi(\widehat C_w)\chi.
\]
\end{corollary}
Together with Theorem~\ref{thm:smooth.intro}, this computes $P_{H_w}(u)$ for every smooth $w$.
In type $A$, the Fourier pairing restricts to the Kronecker delta on $\Irr W$, and every Hessenberg space is of the form $H_w$ for a suitable smooth element $w$ \cite[Proposition~3.1]{haiman}.
Consequently, the formula computes $P_H(u)$ for every Hessenberg space $H$ in type $A$.

The same formula gives the graded $W$-characters of the intersection cohomology of regular semisimple Lusztig varieties.
In Proposition~\ref{prop:fiber.IC}, we prove that $\mathcal H^i(\overline K_w)|_{\Grs}$ is a local system whose stalk at $\mathbf s\in\Grs$ is canonically 
isomorphic to the intersection cohomology of the fiber:
\begin{equation}\label{eq:stalk}
 \mathcal H^i(\overline K_w)_{\mathbf s}
 \simeq \IH^i(\overline Y_w(\mathbf s)),
\end{equation}
where $\IH^*$ denotes intersection cohomology with middle perversity and coefficients in $\EE$.
These local systems on $\Grs$ define the monodromy action in the following result.

\begin{corollary}\label{cor:lus.formula}
For every $w\in W$ and $\mathbf s\in\Grs$, the odd intersection cohomology of $\overline Y_w(\mathbf s)$ vanishes.
The monodromy action on $\IH^*(\overline Y_w(\mathbf s))$ factors through the quotient
$\pi_1(\Grs,\mathbf s)\twoheadrightarrow W$
induced by the Grothendieck--Springer $W$-torsor $\overline Y_e|_{\Grs}\to\Grs$.
With this $W$-action, we have
\[
\sum_i \IH^{2i}(\overline Y_w(\mathbf s))u^i =\sum_{\chi,\psi\in\Irr W}
	\{\chi,\psi\}\widetilde\psi(\widehat C_w)\chi.
\]
\end{corollary}

Over $\C$, Abreu and Nigro established local constancy, the stalk identification \eqref{eq:stalk}, and the factorization of monodromy through $W$, using topological fibrations and localization in torus-equivariant cohomology \cite{AN}.
Their affine paving of the corresponding Bott--Samelson varieties implies the vanishing of odd intersection cohomology by the decomposition theorem \cite[Theorem~1.6(1)]{AN}.
They also computed the graded dimensions of the invariants of these intersection cohomology groups under parabolic subgroups of $W$, obtaining the full graded $W$-character in type $A$ \cite[Theorem~1.5]{AN}.

We give algebraic proofs of these assertions valid over any algebraically closed field $\Bbbk$.
The key observation for the stalk identification is that the map
\beq\label{eq:smooth.map.Lusztig}
G\times\overline Y_w(\mathbf s)\lra
\overline{\fX}_w,\quad (g,x)\mapsto(gx,g\mathbf s x),
\eeq
is smooth.
Smooth pullback then identifies the restriction of $\IC_{\overline Y_w}$ to the fiber with $\IC_{\overline Y_w(\mathbf s)}$.
The Bott--Samelson resolution and the decomposition theorem give local constancy.
We deduce the factorization of monodromy from the classification of unipotent character sheaves with full support (Proposition~\ref{prop:mono.fact}).

Taking stalks in Lusztig's multiplicity formula then gives Corollary~\ref{cor:lus.formula}.
Independently of the present work, Trinh \cite{trinh2026haimansconjecturespringersrepresentations}
obtained the same formula over $\C$ by combining Abreu and Nigro's results with Lusztig's multiplicity formula.

\subsection{Hessenberg schemes defined by cones}
\label{ss:nonlinear.intro}
The construction of Hessenberg varieties extends naturally by replacing the Hessenberg space $H$ with the inverse image in $\fg$ of a $B$-stable closed cone in $\fg/\fb$ \cite[Section~5.1]{BHL25}.
For such schemes, there is a Lie algebra analogue of the decomposition of $\overline K_w$ into shifts of unipotent character sheaves.

Let $\overline C\subset\fg/\fb$ be a nonempty pure-dimensional $B$-stable closed cone, and let $C\subset\fg$ be its inverse image.
Set $X_C=G\times^BC$, and consider the proper morphism
\[
\pi_C:X_C\lra\fg,\quad
[g,z]\mapsto\Ad(g)z.
\]
Its fibers $X_C(s):=\pi_C^{-1}(s)$ are Hessenberg schemes in the sense of \cite{BHL25}.

\begin{theorem}[Theorem~\ref{thm:cone.character}]
\label{thm:cone.character.intro}
For every $i$, the sheaf $\mathcal H^i(\pi_{C*}\IC_{X_C})|_{\grs}$ is a local system whose stalk at $s\in\grs$ is canonically isomorphic to $\IH^i(X_C(s))$.
The resulting monodromy action factors through $\pi_1(\grs,s)\twoheadrightarrow W$.
\end{theorem}

To identify the stalks with the intersection cohomology of the fibers, we use the Lie algebra analogue of the smooth map in \eqref{eq:smooth.map.Lusztig}.
The monodromy statement follows from a decomposition of $\pi_{C*}\IC_{X_C}$ into shifts of \emph{nilpotent character sheaves} on $\fg$, the Lie algebra analogues of unipotent character sheaves on $G$.
Their nonzero restrictions to $\grs$ are, up to shift, the local systems associated with $\chi\in\Irr W$; see Subsection~\ref{ss:nonlinear.support}.

\medskip

We apply this construction to the tangent cones of Schubert varieties at $eB$.
Let $D_w\subset\fg$ be the inverse image of $C_{eB}\overline X_w\subset\fg/\fb$.
When $w$ is smooth, we have $D_w=H_w$ and $X_{D_w}(s)=X_{H_w}(s)$.
Using Theorem~\ref{thm:main}, we extend \eqref{eq:intro.smooth.comparison} to $X_{D_w}(s)$ whenever specialization at $eB$ takes $\IC_{\overline X_w}$ to $\IC_{C_{eB}\overline X_w}$.
We then prove that this holds under the following local hypothesis.

Two pointed schemes $(X,x)$ and $(Y,y)$ are \emph{smoothly equivalent} if there is a pointed scheme $(Z,z)$ admitting morphisms to both that are smooth at $z$.
An \emph{affine cone} is a closed subscheme of a vector space defined by a homogeneous ideal with all coordinates of degree one.

\begin{theorem}[Theorem~\ref{thm:local.conical} and
Proposition~\ref{prop:nonlinear.specialization}]
\label{thm:local.conical.intro}
Suppose that $(\overline X_w,eB)$ is smoothly equivalent to an affine cone pointed by the origin.
Then, for every $s\in\grs$ and $\mathbf s\in\Grs$, there is an isomorphism of graded $W$-representations
\[
\IH^*(X_{D_w}(s))
\simeq \IH^*(\overline Y_w(\mathbf s)).
\]
In particular, these groups vanish in odd degrees, and their graded $W$-character is $\widehat h(\widehat C_w)$.
\end{theorem}

The hypothesis holds for smooth Schubert varieties, inverse images of Schubert varieties in minuscule or cominuscule flag varieties, and covexillary ($3412$-avoiding) permutations in type $A$ (Corollary~\ref{cor:local.conical.classes}).
Thus the theorem extends \eqref{eq:intro.smooth.comparison} to these classes of singular $w\in W$.

\medskip

We also record that for $G=\mathrm{GL}_n$, every simple perverse constituent of $\pi_{C*}\IC_{X_C}$ has full support (Proposition~\ref{prop:cone.typeA.support}), extending the full-support theorem of Precup and Sommers to this setting. 
In type $G_2$, this fails even under the hypothesis of Theorem~\ref{thm:local.conical.intro} (Corollary~\ref{cor:G2.full.support}).

\subsection{Explicit computations}
For the computations described below, we work over $\C$.

For codimension-one Hessenberg spaces in arbitrary type, we obtain a uniform formula for $P_H(u)$ (Proposition~\ref{prop:codim1}).
We also determine when every graded piece is a permutation representation.
For an irreducible root system, this holds exactly in types $A_n$, $B_n$, $C_n$, $F_4$, and $G_2$ (Corollary~\ref{cor:codim1:permutation}).
The proof of Proposition~\ref{prop:codim1} also determines the stalk at the origin of the intersection cohomology complex of the minimal nilpotent orbit closure and relates it to the dimension of the corresponding Springer representation (Corollary~\ref{cor:min.orbit.stalk}).

Outside type $A$, complete computations of $P_H(u)$ for all Hessenberg spaces were previously known only in types $G_2$ and $B_3$, due to Xue \cite[Table~5]{xue} and Precup and Sommers \cite[Section~8]{PS}, respectively.
We present a complete computation in type $C_3$ in detail.
To this end, we use the character formula to compute $P_{H_w}(u)$ for smooth $w\in W$ and then apply the modular law \cite{PS} to reduce the remaining computations.

We first illustrate this approach in type $G_2$.
The codimension-one formula and the modular law give the remaining two characters, recovering the calculations of Xue \cite[Table~5]{xue}.
We also determine, for every $w\in W$, the \emph{nonprincipal} part of $\Hess(\mathcal F_w)$, namely the part whose simple perverse constituents do not occur in $\Spr^\vee$ (see Subsection~\ref{ss:Hess.blocks} and Proposition~\ref{prop:G2.nonprincipal}).

In type $C_3$, the modular law reduces the remaining cases to the computation of a single character, which we determine using the geometry of nilpotent orbits (Proposition~\ref{prop:C3.codim3}).
For exactly nineteen of the twenty Hessenberg spaces, every graded piece of the cohomology is a permutation representation.
For the remaining space, even the total cohomology is not a permutation representation (Corollary~\ref{cor:C3.permutation}).
The complete tables are given in Section~\ref{s:characters}.

\smallskip

More generally, we combine our character formula for smooth $w$ with Precup--Sommers' modular law and decomposition formulas for Hessenberg ideal fibers \cite[Corollary~3.8, (4.7), and Propositions~5.2 and~6.5]{PS} to compute $P_H(u)$.
We explain how these methods are combined in Subsection~\ref{ss:higher.rank} and illustrate their use in type $C_3$ in Example~\ref{ex:C3.ideal}.

We have verified computationally that this approach determines $P_H(u)$ for every Hessenberg space $H$ in types $B_n$, $C_n$, and $D_n$ for $n\leq10$, and in all exceptional types.
We expect this to hold in every Lie type and in arbitrary rank.

In forthcoming joint work with Eunjeong Lee \cite{BHLL}, we use this approach to give an explicit algorithm for computing $P_H(u)$ for all abelian Hessenberg spaces in type $C_n$, for arbitrary $n$.

\subsection{Positivity and unimodality conjectures}
The positivity phenomena predicted by the Stanley--Stembridge and Shareshian--Wachs conjectures for Hessenberg varieties \cite{stanley-stembridge,SW} and by Haiman's conjecture for Lusztig varieties \cite{haiman}, motivate analogous questions in other Lie types.
In types $B$ and $C$, Lesnevich asks whether the graded cohomology of regular semisimple Hessenberg varieties is double-$h$ positive \cite[Question~A.4]{Lesnevich}.

Our computations support an affirmative answer in ranks up to 10 and suggest that the coefficient polynomials are also unimodal.
In Section~\ref{s:positivity}, we strengthen Lesnevich's question by adding coefficientwise unimodality and formulate a corresponding conjecture in type $D$. 
We also propose analogous conjectures for the intersection cohomology of regular semisimple Lusztig varieties in types $B$, $C$, and $D$, together with partial flag extensions of the positivity assertions. See Conjectures~\ref{conj:classical.positivity} and~\ref{conj:classical.partial.positivity}.

We also note that Trinh proposes a different extension of Haiman's conjecture to other Weyl groups, using characters arising from the cohomology of Springer fibers \cite[Conjecture~B]{trinh2026haimansconjecturespringersrepresentations}.

\subsection{Organization}
The paper is organized as follows.
Section~\ref{s:rs.Hess} reviews the geometric and cohomological properties of regular semisimple Hessenberg varieties.
In Section~\ref{s:microloc}, we review specialization, Fourier transform, and microlocalization.
In Section~\ref{s:hess}, we construct the Hessenberg functor and prove Theorems~\ref{thm:smooth.intro} and~\ref{thm:main}.
Section~\ref{s lusztig character} computes the Hessenberg character using Lusztig's multiplicity formula and interprets its values $\widehat h(\widehat C_w)$ as the graded $W$-characters of the intersection cohomology of regular semisimple Lusztig varieties.
Section~\ref{s:nonlinear} studies nilpotent character sheaves and monodromy for Hessenberg schemes defined by cones.
It gives sufficient conditions for extending \eqref{eq:intro.smooth.comparison} to those defined by tangent cones of Schubert varieties. 
Section~\ref{s:codim1} computes $P_H(u)$ for codimension-one Hessenberg spaces in arbitrary type.
In Section~\ref{s:characters}, we compute $P_H(u)$ for every Hessenberg space in types $G_2$ and $C_3$ and describe the approach used for computations in higher rank.
Section~\ref{s:positivity} formulates positivity and unimodality conjectures in types $B$, $C$, and $D$, supported by computational evidence.

\subsection*{Acknowledgements}
We thank Minseong Kwon for answering questions about Schubert varieties in type $G_2$, Eunjeong Lee for writing and running code for the character formula, Wille Liu for helpful explanations about character sheaves on Lie algebras, and Hyeonjun Park for discussions on the Fourier--Sato transform.
We also thank Minh-T\^am Trinh for sharing an earlier draft of his preprint \cite{trinh2026haimansconjecturespringersrepresentations} and for helpful discussions.

Hong was supported by the Institute for Basic Science IBS-R032-D1. 
Lee was supported by the Institute for Basic Science IBS-R032-D1 and the KIAS Individual Grant (HP109201).

\subsection*{AI disclosure}
We used ChatGPT (OpenAI) for language editing,
proofreading, and literature searches while preparing
the first draft.
During revision, we used ChatGPT with GPT 6 Astra
to discuss and develop additional results and arguments
in Subsection~\ref{ss:Hess.blocks},
Section~\ref{s:nonlinear}
(primarily Subsection~\ref{ss:nonlinear.conical}),
Proposition~\ref{prop:G2.nonprincipal}, Corollary~\ref{cor:G2.full.support},
and Section~\ref{s:positivity}.
We also used it to derive
Corollaries~\ref{cor:codim1:permutation}
and~\ref{cor:C3.permutation}
and construct Table~\ref{tab:C3.permutation}
from our preexisting character computations.
We also used ChatGPT to generate and execute code
for the computations described in
Subsection~\ref{ss:higher.rank}.
The authors are responsible for the final text
and mathematical content.

\renewcommand{\thetheorem}{\thesection.\arabic{theorem}}
\section{Regular semisimple Hessenberg varieties} \label{s:rs.Hess}

We record the basic geometric properties of regular semisimple Hessenberg varieties and the monodromy $W$-representation on their cohomology in the generality of the introduction.
These are standard when the base field is $\C$ \cite{demari-procesi-shayman,BrCh18,balibanu-crooks}, and we adapt the arguments there.

We keep the notation of the introduction and set
\[\ft^{\mathrm{rs}}:=\ft\cap\grs.\]
For $s\in \grs$, we write $T_s:=Z_G(s)^\circ$ for the identity component of its centralizer $Z_G(s)$.

\begin{theorem}\label{thm:DMPS}
	Let $H\subset \fg$ be a Hessenberg space.
	Let $X_H^{\mathrm{rs}}:=\pi_H^{-1}(\grs)$. Then  the restriction 
	\[\pi_H^{\mathrm{rs}}:X_H^{\mathrm{rs}}\lra \grs\]
	of $\pi_H$ over $\grs$ is smooth and projective of relative dimension $\dim H/\fb$.
	In particular, $X_H(s)$ is smooth for every $s\in \grs$. 
	Moreover, its $T_s$-fixed locus is $\cB^{T_s}$, consisting of $\lvert W\rvert$ reduced points, and it admits an affine paving indexed by these points. 
\end{theorem}
\begin{proof}
It is immediate that $\pi_H^\rs$ is projective, since $\pi_H$ is projective. For the remaining assertions, the arguments of \cite[Sections~II--III]{demari-procesi-shayman} apply over $\Bbbk$, using $T_s=Z_G(s)^\circ$ as the acting torus.
By conjugation, we may assume that $s\in\mathfrak t^{\mathrm{rs}}$ and $T_s=T$. 
Since $\mathfrak t\subset H$, we have
\[X_H(s)^T=\cB^T=\{nB:n\in N_G(T)\}.\]

The morphism $\pi_H$ is smooth at every such fixed point, since
$H+[\fg,\Ad(n^{-1})s]=\fg$.
If its non-smooth locus met $X_H(s)$, the intersection would be a nonempty closed $T$-stable subset of the projective variety $X_H(s)$ and hence would contain a $T$-fixed point, a contradiction.
Thus $\pi_H$ is smooth over $\grs$, with relative dimension $\dim H/\fb$.

Choose a one-parameter subgroup of $T$ whose fixed locus on
$X_H(s)$ is $X_H(s)^T$. Since $X_H(s)$ is smooth and projective
with isolated fixed points, the Bia{\l}ynicki--Birula decomposition
gives an affine paving indexed by these points.
\end{proof}

We recall  standard facts about the Grothendieck--Springer resolution. 
Let $\tgrs:=\pi_\fb^{-1}(\grs)$ be the regular semisimple locus of $X_\fb=G\times^B\fb$. Then $\tgrs\simeq G\times^B\fb^{\rs}$, where $\fb^{\rs}=\fb\cap \grs$.

\begin{lemma}[{\cite[Lemmas~13.3--13.4 and Section~13.5]{jantzen-nilpotent}}]\label{lem:W-torsor}
For every $s\in\grs$, we have $X_\fb(s)=\cB^{T_s}$.
The morphism $G\times^T\ft^{\mathrm{rs}}\to\tgrs$,
$[g,s]\mapsto[g,s]$
is an isomorphism.
Moreover, the restriction 
\[\pi_\fb^{\mathrm{rs}}:\tgrs\lra \grs\] 
of the Grothendieck--Springer resolution over $\grs$
is a $W$-torsor.
\end{lemma}

Since $\pi_H^\rs$ is smooth and proper, 
$R^i\pi_{H*}^{\mathrm{rs}}\underline \EE$ is locally constant for every $i\in \ZZ$, and its stalk at $s\in \grs$ is isomorphic to $\coh^i(X_H(s))$ equipped with the monodromy representation of $\pi_1(\grs,s)$.

\smallskip

The $W$-torsor in Lemma~\ref{lem:W-torsor} induces a quotient $\pi_1(\grs,s)\twoheadrightarrow W$, which is well-defined up to an inner automorphism of $W$. 
We prove the following theorem.

\begin{theorem}\label{thm:Hess.monodromy}
For every $i\in \ZZ$, the local system $R^i\pi_{H*}^{\mathrm{rs}}\underline{\EE}$ becomes constant when pulled back along the $W$-torsor $\tgrs\to \grs$, that is, $(\pi_{\mathfrak b}^{\mathrm{rs}})^*R^i\pi_{H*}^{\mathrm{rs}}\underline{\EE}$ is constant.
Equivalently, the monodromy action of
$\pi_1(\grs,s)$ on $\coh^i(X_H(s))$ factors through the quotient 
$\pi_1(\mathfrak g^{\mathrm{rs}},s)\twoheadrightarrow W$.
\end{theorem}

The following corollary is immediate from this.
\begin{corollary}
For every $i$, there are finite-dimensional $\EE$-vector spaces
$M_{H,\chi}^i$ such that
\[R^i\pi_{H*}^{\mathrm{rs}}\underline{\EE} \cong \bigoplus_{\chi\in\Irr W}\mathcal L_\chi\otimes_{\EE} M_{H,\chi}^i,\]
where $\mathcal L_\chi$ is the local system associated with the
Grothendieck--Springer $W$-torsor and $\chi$.
\end{corollary}

The rest of this section is devoted to the proof of Theorem~\ref{thm:Hess.monodromy}.
Over $\C$, the theorem is proved in \cite{BrCh18,balibanu-crooks}, and we follow their approach. 

We recall the following standard fact about the $T_s$-equivariant cohomology $\coh^*_{T_s}(X_H(s))$.

\begin{lemma}\label{lem:equiv.formality}
Let $H\subset \fg$ be a Hessenberg space, and let $s\in \grs$. The restriction map
\beq\label{eq:localization}\coh_{T_s}^*(X_H(s))\lra\coh_{T_s}^*(X_H(s)^{T_s})\eeq
is injective, and the natural map
\beq\label{eq:equivtosing}\coh_{T_s}^*(X_H(s)) \lra \coh^*(X_H(s))\eeq
is surjective. 
More precisely, $\coh^*(X_H(s))\cong \coh_{T_s}^*(X_H(s)) \otimes_{\coh^*(BT_s)}\EE$.
Here $\coh^*(BT_s)=\coh^*_{T_s}(\mathrm{pt})$, and $\EE$ is a
module over this ring through the augmentation.
\end{lemma}

\begin{proof}
The affine paving implies that $\coh^{\mathrm{odd}}(X_H(s))=0$, and hence the $T_s$-action is equivariantly formal.
The assertion about \eqref{eq:equivtosing} follows.
The equivariant localization theorem shows that the kernel of the restriction map is $\coh^*(BT_s)$-torsion, while equivariant formality implies that $\coh_{T_s}^*(X_H(s))$ is torsion-free, therefore the restriction map \eqref{eq:localization} is injective.
\end{proof}

\begin{proof}[Proof of Theorem~\ref{thm:Hess.monodromy}]
	Pulling back $\pi_H^{\mathrm{rs}}$ by the Grothendieck--Springer $W$-torsor, we get
\[\widetilde \pi_H:\widetilde X_H^{\mathrm{rs}}:=X_H^{\mathrm{rs}}\times_{\grs}\tgrs\lra \tgrs.\]
Using the isomorphism $\tgrs\simeq G\times^T\ft^{\mathrm{rs}}$ from Lemma~\ref{lem:W-torsor}, let $t\in T$ act by multiplication by $gtg^{-1}$ on the fiber over $[g,s]$.
This defines a fiberwise $T$-action on $\widetilde X_H^{\mathrm{rs}}$, with fixed locus
\[(\widetilde X_H^{\mathrm{rs}})^T\simeq W\times\tgrs.\]

Consider the induced morphisms from the quotient stacks
\[\overline \pi_H:\big[\widetilde X_H^{\mathrm{rs}}/T\big]\lra \tgrs \and \overline \pi_H^T:\big[(\widetilde X_H^{\mathrm{rs}})^T/T\big]\lra \tgrs.\]
Then their higher direct images of $\EE$ are locally constant and commute with base change. Indeed, in each fixed degree they can be computed by the smooth proper Borel approximations
\[
(\widetilde X_H^{\mathrm{rs}}\times U_N)/T\lra \tgrs \and ((\widetilde X_H^{\mathrm{rs}})^T\times U_N)/T\lra \tgrs,
\]
where $U_N=(\mathbb A^{N+1}\setminus\{0\})^{\rank G}$, equipped with a componentwise scaling action of $T\simeq \Gm^{\rank G}$, 
and $N$ is sufficiently large.
For the definition of direct images from quotient stacks, see \cite{bernstein-lunts} in the Betti setting and \cite[Section~8]{LOII} in the $\ell$-adic setting.

Moreover, there are natural morphisms of these local systems
\[R^i\overline \pi_{H*}\underline \EE \hookrightarrow R^i\overline\pi_{H*}^T\EE \and R^i\overline \pi_{H*}\EE \twoheadrightarrow R^i\widetilde\pi_{H*}\EE,\]
which are injective and surjective respectively by Lemma~\ref{lem:equiv.formality}.
Since $(\widetilde X_H^{\mathrm{rs}})^T\simeq W\times \tgrs$ by Lemma~\ref{lem:W-torsor}, the local system $R^i\overline \pi_{H*}^T\EE$ is constant.
It follows that $R^i\overline \pi_{H*}\EE$ is constant, and so is
its quotient $R^i\widetilde\pi_{H*}\EE$.
By proper base change, this proves the theorem.
\end{proof}

\begin{remark}\label{rem:Hess.dot.action}
For $s\in\ft^\rs$, the $W$-action on $\coh_T^*(X_H(s)^T)$ obtained by descending the constant local systems $R^i\overline\pi_{H*}^T\EE$ along $\tgrs\to\grs$ is Tymoczko's dot action.
Indeed, the $W$-action induced by the $W$-torsor permutes the fixed points in $X_H(s)^T\cong W$ and changes the identification of the acting torus with $T$ by conjugation, inducing the natural Weyl action on $\coh^*(BT)$.
Together, these give Tymoczko's dot action on $\coh_T^*(X_H(s)^T)\cong \coh^*(X_H(s)^T)\otimes \coh^*(BT)$.
Compatibility with \eqref{eq:localization} and \eqref{eq:equivtosing} then identifies the monodromy action on $\coh^*(X_H(s))$ with the dot action.

Moreover, the image of \eqref{eq:localization} has the same moment-graph description over $\Bbbk$ as over $\C$, by \cite[Lemme~3.1]{chaudouard-laumon}.
This description and the dot action depend only on the root datum and the roots defining $H$.
Thus the graded $W$-character $\coh^*(X_H(s))$ is independent of the choice of $s\in\grs$ and the algebraically closed base field $\Bbbk$, provided that $\grs\neq\emptyset$.
\end{remark}

\begin{remark}\label{rem:partial.Hess}
Let $P\supset B$ be a parabolic subgroup such that $H$ is $P$-stable. Consider
\[
\pi_{H,P}:X_{H,P}:=G\times^P H\lra\fg,\quad [g,z]\mapsto\Ad(g)z.
\]
Its fiber at $s\in\grs$ is
\[
X_{H,P}(s)=\{gP\in G/P:\Ad(g^{-1})s\in H\}.
\]
The projection $X_H\to X_{H,P}$ is a
$P/B$-bundle. Hence $\pi_{H,P}$ is smooth and projective
over $\grs$.

Products of first Chern classes of line bundles on $G/B$
give a Leray--Hirsch isomorphism
\[
\coh^*(X_H(s))
\simeq
\coh^*(X_{H,P}(s))\otimes\coh^*(P/B)
\]
compatible with monodromy, which acts trivially on the
second factor.
In particular, pullback is injective, so
Theorem~\ref{thm:Hess.monodromy} implies that monodromy on
$\coh^*(X_{H,P}(s))$ factors through $W$.
Set
\[
P_{H,P}(u):=
\sum_{j\geq0}\coh^{2j}(X_{H,P}(s))u^j.
\]
Then, writing $p_{P/B}(u)= \sum_{j\geq0}\dim\coh^{2j}(P/B)u^j$, we have
\beq\label{eq:flagbundle}
P_H(u)=p_{P/B}(u)P_{H,P}(u).
\eeq
\end{remark}

\section{Specialization and microlocalization}\label{s:microloc}

In this section, we review the specialization functor and the microlocalization functor.
\subsection{Notations} \label{ss notations}
For a scheme $X$ of finite type over $\Bbbk$, let $\mathscr D(X)$ denote the bounded derived category of constructible sheaves on $X$ (with coefficients in $\EE$).
For $K\in\mathscr D(X)$, $\cH^i(K)$ denotes the $i$-th cohomology sheaf of $K$.
We write $\Perv(X) \subset \mathscr D(X)$ for the abelian full subcategory of perverse sheaves, and denote the perverse cohomology functors by
\[{}^p\cH^i : \mathscr D(X) \lra \Perv(X).\]
We say that a complex $K \in \mathscr D(X)$ is \emph{semisimple} if
$K \cong \bigoplus_i {}^p\cH^i(K)[-i]$
and each ${}^p\cH^i(K)$ is semisimple in $\Perv(X)$.

For a locally closed smooth irreducible subvariety $V \subset X$ and an irreducible local system $\mathcal L$ on $V$, we denote by
$\mathrm{IC}(\overline V, \mathcal L)$
the \emph{intersection cohomology complex}, normalized so that
\[\ICperv(\overline V,\mathcal L):=\mathrm{IC}(\overline V, \mathcal L)[\dim V]\in \Perv(X).\]
When $\mathcal L$ is trivial, we simply write them as $\IC_{\overline V}$ and $\ICperv_{\overline V}=\IC_{\overline V}[\dim V]$.
For a scheme $X$, we take $\IC_X$ to be the direct sum of the unshifted intersection cohomology complexes of the irreducible components of its reduced subscheme $X_{\mathrm{red}}$.
Some authors instead call $\ICperv(\overline V,\mathcal L)$ the intersection cohomology complex (see, e.g., \cite{Achar,D04,M19}).

The split Grothendieck group $\mathscr K_\oplus(\Perv(X))$ is generated by the isomorphism classes $[K]$, subject to the relations
\[[K_1 \oplus K_2] = [K_1] + [K_2].\]
We write $\mathscr K(\Perv(X))$ for the (usual) Grothendieck group, in which the relations come from short exact sequences, and denote the class of an object $K$ by $\{K\}$.
The natural homomorphism from the split Grothendieck group to the Grothendieck group sends $[K]$ to $\{K\}$.

\medskip
When an algebraic group $G$ acts on $X$,  let $\mathscr D_G(X)$ denote the derived category of bounded $G$-equivariant constructible complexes on $X$ \cite{bernstein-lunts,Joshua} (see also \cite[Definition 6.4.4]{Achar}), and let $\Perv_G(X)$ denote the abelian full subcategory of $G$-equivariant perverse sheaves. 
We write $\Semis_G(X)\subset \mathscr D_G(X)$ for the full subcategory of semisimple complexes, and $\mathscr K_{\oplus}(\Semis_G(X))$ for its split Grothendieck group.

All standard functors on $\mathscr D(X)$ admit equivariant analogues on $\mathscr D_G(X)$; see \cite[Section~6.5]{Achar}.

\subsection{Specialization}\label{ss:specialization}

Let $X$ be a scheme of finite type over $\Bbbk$ and let $Z\subset X$ be a closed subscheme. 
Let $M_{Z/X}$ denote the total space of the deformation of $X$ to the normal cone $C_{Z/X}$.
It is the complement of the strict transform of $X\times\{0\}$ in the blowup 
of $X\times \A^1$ along $Z\times \{0\}$.
The natural projection $\pi:M_{Z/X}\to\A^1$ fits into the following Cartesian diagram:
\[\begin{tikzcd}
    	C_{Z/X}\arrow[r,hook,"i"]\arrow[d]&M_{Z/X}\arrow[d,"\pi"]& X\times \Gm\arrow[l,hook',"j"']\arrow[d]\\ \{0\}\arrow[r,hook]&\A^1&\Gm\arrow[l,hook']
    \end{tikzcd}\]
where $\Gm=\A^1-\{0\}$. When $Z\hookrightarrow X$ is a regular immersion (for example, when $X$ and $Z$ are smooth), $C_{Z/X}$ is the normal bundle $N_{Z/X}$. 
We write $\pr_1$ for the projection $X\times \Gm \to X$.

\smallskip

Let $C$ be a cone with its scaling $\Gm$-action.
A complex $K\in\mathscr D(C)$ is called \emph{monodromic} if, for every $x\in C(\Bbbk)$, the pullback $a_x^*K$ along the orbit map
\[
a_x:\Gm\longrightarrow C,\quad t\mapsto t\cdot x,
\]
has locally constant cohomology sheaves.
In the \'etale setting, these local systems are required to be tame.
We denote the full subcategory of monodromic complexes by
$\mathscr D(C)^{\mon}$; see
\cite[Definition~3.1 and Sections~5 and~9]{verdier}.

\begin{definition}[Specialization {\cite[(8.2)]{verdier}}]
	The \emph{specialization} functor is defined by 
    \[\sp_{Z/X}:\mathscr D(X)\to \mathscr D(C_{Z/X})^{\mon}, \quad F\mapsto \psi_\pi(j_*\pr_1^*F)\]
    where $\psi_\pi$ denotes the nearby cycle functor associated to $\pi$ \cite{SGA7-nearby}. 
    By \cite[(SP1)]{verdier}, the complex $\sp_{Z/X}F$ is monodromic.
    See also \cite[Appendix~A.1]{bezrukavnikov2026affinespringerfibersmall}.
    
    When an algebraic group $G$ acts on $X$ and $Z$ is $G$-stable, we have a $G$-equivariant analogue
    \[\sp_{Z/X}:\mathscr D_G(X) \lra \mathscr D_{G}(C_{Z/X})^{\mon}.\]
    We often omit $X$ in the notation and write it simply as $\sp_Z$.
\end{definition}

Let $i_Z:Z\hookrightarrow X$ and $0_Z:Z\hookrightarrow C_{Z/X}$ be
the inclusion and the zero section. By \cite[(SP5)]{verdier},
\beq\label{eq:sp.zero}
0_Z^*\sp_{Z/X}F\simeq i_Z^*F.
\eeq

The specialization functor satisfies several basic properties, including smooth base change.
By Verdier \cite{verdier-riemann-roch} (see also \cite[Chapter~5]{Fu98}), $M_{Z/X}$ admits a description by relative spectrum:
\[M_{Z/X}\simeq \Spec_{X\times \A^1} \bigoplus_{n\in \ZZ}\mathcal I_{Z/X}^nt^{-n}\]
where $t$ is the coordinate of $\A^1$ and $\mathcal I_{Z/X}$ denotes the ideal sheaf defining $Z\subset X$. Here $\mathcal I^n_{Z/X}=\mathcal O_X$ for $n\leq 0$. 
In particular, when we have a diagram 
\beq\label{eq:XYWZ}\begin{tikzcd}
		Z\arrow[r,hook]\arrow[d,"f_Z"']& X\arrow[d,"f"]\\ W\arrow[r,hook]& Y
	\end{tikzcd}\eeq
such that $Z$ and $W$ are closed subschemes of $X$ and $Y$,
respectively, the induced morphisms fit into the following commutative diagram:
\[\begin{tikzcd}
	C_{Z/X}\arrow[r,hook]\arrow[d,"C(f)"']&M_{Z/X}\arrow[d,"\tilde f"] &X\times \Gm \arrow[l,hook',"j_X"']\arrow[d,"f\times \id_{\Gm}"] \\ C_{W/Y}\arrow[r,hook]\arrow[d]&M_{W/Y}\arrow[d,"\pi"]& Y\times \Gm \arrow[l,hook',"j_Y"']\arrow[d]\\ \{0\}\arrow[r,hook] & \A^1& \Gm \arrow[l,hook']
\end{tikzcd}\]

\begin{proposition}[Smooth pullback for specialization]\label{prop:sp.smooth}
	Suppose we have a commutative diagram 
	\eqref{eq:XYWZ}
	which is not necessarily Cartesian. 
	Suppose that both $f$ and $f_Z$ are smooth. 
	Then $\tilde f$ is smooth, and 
	there is a natural isomorphism $C(f)^* \sp_{W/Y}\simeq\sp_{Z/X}f^*$.
\end{proposition}
\begin{proof}
	\'Etale locally on the source and the target, we may assume that $X= Y\times \A^d$ with $f$ being the first projection and $Z= W\times \A^{d-c}$ where $\A^{d-c}\subset \A^{d}$ is the affine subspace where the last $c$ coordinates vanish. Then one can check that $\tilde f$ is smooth with fibers $\A^d$. Thus,
	\[C(f)^*\sp_{W/Y}= C(f)^*\psi_\pi j_{Y*}\pr_1^*\simeq \psi_{\pi\circ \tilde f}j_{X*}(f\times \id_{\Gm})^*\pr_1^* \simeq \psi_{\pi\circ \tilde f}j_{X*}\pr_1^*f^*= \sp_{Z/X}f^*,\]
	where the displayed isomorphisms follow by combining the smooth pullback property for the nearby cycle functor (cf.~\cite[Lemma~4.1.7(2)]{Achar}) and smooth base change \cite[Theorem~2.2.2]{Achar}.
\end{proof}
When \eqref{eq:XYWZ} is Cartesian, this specializes to \cite[(SP3)]{verdier}. 
See also \cite[Proposition~4.2.5]{KS94}.

\subsection{Fourier transform and microlocalization}\label{ss:micro}
Let $S$ be a scheme of finite type over $\Bbbk$, and let $V\to S$ be a vector bundle of rank $r$ with dual $V^\vee\to S$.
Denote by
\[\begin{tikzcd}
	&V\times_SV^\vee\arrow[rr,"m"]\arrow[ld,"\pr_1"']\arrow[rd,"\pr_2"]&&\A^1\\V&&V^\vee
\end{tikzcd}\]	
the first and second projection maps and the natural pairing.
We write $\mathscr D(V)^{\mon}\subset\mathscr D(V)$ for the full subcategory of monodromic complexes with respect to the fiberwise scaling $\Gm$-action on $V$, as defined in Subsection~\ref{ss:specialization}.

\begin{definition}[Fourier transform]\label{def:Fourier}
	For $\Bbbk=\CC$ or $\overline \FF_q$, the \emph{Fourier transform} is defined as
	\[\Four_V:\mathscr D(V)^{\mon}\xrightarrow{~\sim~}\mathscr D(V^\vee)^{\mon},\quad F\mapsto \pr_{2!}(\pr_1^* F\otimes  m^* \mathcal J)[r],\]
	where $\mathcal J$ is a sheaf on $\A^1$ defined as follows. 
	
	When $\Bbbk=\C$, $\mathcal J$ is the constant sheaf on the (non-algebraic) subset $\{z\in \C:\Re(z)\leq 0\}\subset \C$, extended by zero to $\C$. In this case, $\Four_V$ is called the \emph{Fourier--Sato transform} \cite[Section~3.7]{KS94}.
	Although $\mathcal J$ is not algebraic, the Fourier--Sato transform preserves algebraically constructible complexes since
	\beq\label{eq:Fourier.van}\Four_V\simeq \phi_{\pr_2'}(\pr_2\times m)_*\pr_1^*[r-1]\eeq
	by \cite[(10.3.31)]{KS94}, where $\phi_{\pr_2'}$ is the vanishing cycle functor \cite{SGA7-nearby} associated to $\pr_2':V^\vee\times \C\to \C$.
	
	When $\Bbbk=\overline \FF_q$, $\mathcal J$ is the Artin--Schreier sheaf associated to a nontrivial additive character $\chi:\FF_q\to \overline\QQ_\ell^*$. In this case, $\Four_V$ is called the \emph{Fourier--Deligne transform} \cite{Laumon}. 
	
	Over an arbitrary algebraically closed field, we use Wang's Fourier transform on monodromic complexes \cite{Wang}; 
	see also \cite[Appendix~A]{bezrukavnikov2026affinespringerfibersmall}.
	In positive characteristic, it is isomorphic to the Fourier--Deligne transform on this subcategory \cite[Corollary~4.5]{Wang}.

	When an algebraic group $G$ acts on $S$ and $V$ is a $G$-equivariant vector bundle, we have 
	\[\Four_{V}: \mathscr D_{G}(V)^{\mon}\xrightarrow{~\sim~} \mathscr D_{G}(V^\vee)^{\mon}.\]
\end{definition}
	
\begin{lemma}\label{lem:Four.subbdl}
	Let $U\subset V$ be a subbundle, and let $U^\perp\subset V^\vee$ be the annihilator of $U$. Then 
	\[\Four_V(\underline \EE_{U}[\dim U])\simeq \underline \EE_{U^\perp}[\dim U^\perp](-\rank U).\]
Here the constant sheaves are extended by zero along the closed
embeddings $U\subset V$ and $U^\perp\subset V^\vee$.
	The same holds when $S$ is a $G$-scheme and $U\subset V$ is a $G$-equivariant subbundle.
\end{lemma}
\begin{proof}
	See \cite[Corollary~6.9.14]{Achar}.
\end{proof}

We will also need the compatibility of Fourier transform with proper direct image. 
\begin{lemma}\label{lem:Four.pullback}
	If $f:S'\to S$ is proper and $V$ is a vector bundle over $S$, then
\[\Four_V\circ(f_V)_*\simeq(f_{V^\vee})_*\circ\Four_{f^*V},\]
where $f_V:V\times_SS'\to V$ and $f_{V^\vee}:V^\vee\times_SS'\to V^\vee$ are the induced morphisms.
\end{lemma}
\begin{proof}
	See \cite[Proposition~6.9.15]{Achar}.
\end{proof}

\begin{example}\label{ex:Four.pullback}
	When $f$ is $\cB\to \Spec\Bbbk$ and $V=\fg$, the following diagram  commutes.
\[
\begin{tikzcd}
	\mathscr D_G(\cB\times \fg)^{\mon}\arrow[r,"\pr_{2*}"]\arrow[d,"\Four_{\cB\times \fg/\cB}"']&\mathscr D_G(\fg)^{\mon}\arrow[d,"\Four_\fg"]\\
	\mathscr D_G(\cB\times \fg^\vee)^{\mon}\arrow[r,"\pr_{2*}"]&\mathscr D_G(\fg^\vee)^{\mon}
\end{tikzcd}\]

\end{example}

\begin{definition}[Microlocalization]
	Let $X$ be a scheme of finite type over $\Bbbk$, and let $Z\subset X$ be a regularly embedded closed subscheme. 
	The \emph{microlocalization} functor along $Z$ is defined as the composition of the specialization functor and the Fourier transform
	\[\mu_{Z/X}:\mathscr D(X)\xrightarrow{~\sp_{Z/X}~}\mathscr D(N_{Z/X})^{\mon}\xrightarrow{\Four_{N_{Z/X}}}\mathscr D(N_{Z/X}^\vee)^{\mon}.\]
	
	When a group $G$ acts on $X$ and $Z$ is $G$-stable, we have a $G$-equivariant analogue
	\[\mu_{Z/X}:\mathscr D_G(X)\lra \mathscr D_{G}(N_{Z/X}^\vee)^{\mon}.\]
	When the ambient space $X$ is clear from the context, we simply write $\mu_Z$.
\end{definition}

\begin{remark}\label{rem:mu.texact}
	Our normalization of the Fourier 
	transform, and hence of the microlocalization functor, differs from that of \cite{KS94} by a shift. 
	Indeed, we include the shift $[r]$ in Definition~\ref{def:Fourier} in order to 
	make $\Four_V$ $t$-exact for the perverse $t$-structure.
	Consequently, since $\sp_{Z/X}$ is $t$-exact for the perverse $t$-structure, so is $\mu_{Z/X}$ (\cite[Proposition~10.3.19]{KS94}).
\end{remark}

\section{Hessenberg character}\label{s:hess}
We first recall the Hecke category and its relation to the Iwahori--Hecke algebra.
We then define the Hessenberg functor and the Hessenberg character, and prove Theorems~\ref{thm:smooth.intro} and~\ref{thm:main}.

\subsection{Hecke category and Hessenberg character}
\label{ss Hecke}
Let $S\subset W$ be the set of simple reflections determined by the Borel subgroup $B$, and write $\ell$ and $\leq$ for the length function and Bruhat order, respectively.
Set $\cA=\mathbb Z[u^{\pm1/2}]$.
The Iwahori--Hecke algebra $\hecke$ of $(W,S)$ is the free $\cA$-module with basis $\{T_w\}_{w\in W}$ and multiplication determined by
\begin{align*}
  &T_w T_{w'} = T_{ww'} \qquad\qquad  \text{if } \ell(ww')=\ell(w)+\ell(w'), \\
  &(T_s+1)(T_s-u)=0 \quad \text{for } s\in S.
\end{align*}

There is an involution $\overline{(-)}:\hecke\to\hecke$ defined by
$\overline{u^{1/2}}=u^{-1/2}$ and $\overline{T_w}=T_{w^{-1}}^{-1}$.
The Kazhdan--Lusztig basis $\{C'_w\}_{w\in W}$ is the unique $\cA$-basis of $\hecke$ such that
\[\overline{C'_w}=C'_w \and C'_w = u^{-\ell(w)/2}\left(T_w + \sum_{x< w} P_{x,w}(u)\,T_x\right)\]
for some $P_{x,w}\in\mathbb{Z}[u]$ satisfying $P_{w,w}=1$ and $\deg P_{x,w}\le \frac{1}{2}(\ell(w)-\ell(x)-1)$ for $x<w$.
The polynomials $P_{x,w}$ are uniquely determined, and are called the Kazhdan--Lusztig polynomials.
We set
\[\widehat C_w := u^{\ell(w)/2} C'_w.\]

Relative Schubert varieties $\overline{\mathfrak{X}}_w$ are the closures of the $G$-orbits $G\cdot (B,wB)$ in $\cB\times\cB$, and the simple objects in $\Perv_G(\cB \times \cB)$ are given by 
$\ICperv_{\overline{\fX}_w}=\IC_{\overline \fX_w}[\dim \cB+\ell(w)]$.

The convolution product endows $\mathscr D_G(\cB \times \cB)$ with a monoidal structure, and is given by
\[F \star G:= \mathrm{pr}_{13*} \big(\mathrm{pr}_{12}^* F \otimes \mathrm{pr}_{23}^*G\big)[-\dim \mathcal B],\]
for $F, G \in \mathscr D_G(\cB \times \cB)$ (\cite[Section~7.2, Exercise~7.7.1]{Achar}). 
\beq\label{eq:prij}
\begin{tikzcd}
	&\cB\times\cB\times \cB\arrow[ld,"\pr_{12}"']\arrow[d,"\pr_{13}"]\arrow[rd,"\pr_{23}"]&\\ 
	\cB\times \cB &\cB\times \cB &\cB\times \cB
\end{tikzcd}
\eeq
Here $\pr_{ij}$ denotes the projection to the $i$-th and $j$-th factors.

Moreover, the \emph{geometric Hecke category}, the full subcategory of semisimple complexes
\[\mathsf H_G:=\Semis_G(\cB\times \cB)\subset \mathscr D_G(\cB \times \cB)\]
is closed under $\star$ (\cite[Proposition~7.2.6]{Achar}). In particular, the convolution product induces a ring structure on 
$\mathscr K_{\oplus}(\mathsf H_G)$,
which is naturally an $\cA$-module via $u^{1/2}\cdot [F] := [F[-1]]$.

\begin{proposition}[{\cite[Theorem~7.3.8]{Achar}}]  \label{t hecke schubert}  
The assignment
\[\mathscr K_{\oplus}(\mathsf H_G) \xrightarrow{~\sim~} \hecke, \quad [\ICperv_{\overline{\fX}_w}]\mapsto C_w'
\]
is an isomorphism of $\cA$-algebras. 
It sends the class of $\mathcal F_w=\IC_{\overline \fX_w}[\dim \cB]$ to $\widehat C_w$.
\end{proposition}

We now use this to define the Hessenberg functor.
Recall the moment map
\[
\eta:T_\cB^\vee\lra\fg^\vee,\quad(B',\xi)\mapsto\xi,
\]
whose image is the coadjoint nilpotent cone $\mathcal N^\vee$ defined in the introduction.
When there exists a $G$-equivariant isomorphism $\fg^\vee\to \fg$, we may identify $(\fg/\fb)^\vee$ with the nilpotent radical $\mathfrak u\subset \fg$, and the moment map $\eta$ with the Springer map 
\[G\times^B\mathfrak u \lra \mathcal N, \quad [g,u]\mapsto \Ad(g) u.\] 

\begin{definition}\label{def:Hess}
The \emph{Hessenberg functor} is defined to be the composition 
\beq\label{eq:Hess.funct}
	\mathsf{Hess}:~\mathsf H_G\hookrightarrow \mathscr D_G(\cB \times \cB) \xrightarrow{~\mu_{\Delta}~} \mathscr D_{G} (T^\vee_{\cB})^{\mon}  \xrightarrow{~\eta_{*}~} \mathscr D_{G}(\mathcal N^\vee).\eeq

\end{definition}
Taking graded perverse cohomology groups and using the identification of Proposition~\ref{t hecke schubert}, we obtain an $\cA$-linear map
\[
\hbar:\hecke\lra \mathscr K_\oplus(\Perv_G(\mathcal N^\vee))\otimes\cA,
\quad [F]\mapsto\sum_i[{}^p\cH^i(\mathsf{Hess}(F))]u^{i/2}.
\]
To obtain a map to $\Rep W\otimes \cA$, consider the coadjoint Springer sheaf
\[\Spr^\vee:=\eta_*\underline{\EE}_{T^\vee_{\cB}}[2\dim\cB],\]
endowed with the $W$-action induced by Fourier transform. 
To be precise, write
\[
\mathrm{Groth}:=\pi_{\fb*}\underline{\EE}_{G\times^B\fb}[\dim\fg].
\]
By \cite[Theorem~13.5]{jantzen-nilpotent}, this is the intermediate extension of
$\pi_{\fb*}^{\mathrm{rs}}\underline{\EE}_{X_\fb^{\mathrm{rs}}}[\dim\fg]$, so the $W$-action induced by the $W$-torsor extends uniquely to $\mathrm{Groth}$.
Lemmas~\ref{lem:Four.subbdl} and~\ref{lem:Four.pullback} give
\[\Spr^\vee\simeq\Four_\fg(\mathrm{Groth})(b)\]
where $b:=\dim\fb$.
We equip $\Spr^\vee$ with the transported $W$-action.

For $A\in\Perv_G(\mathcal N^\vee)$, set
\[\mathsf S(A):=\Hom_{\Perv_G(\mathcal N^\vee)}(\Spr^\vee,A),\]
with $W$-action given by $w\cdot f=f\circ w^{-1}$.
We extend $\mathsf S$ $\cA$-linearly to $\mathscr K_\oplus(\Perv_G(\mathcal N^\vee))\otimes\cA$.

\begin{definition}
	We define the \emph{Hessenberg character} by
\[\widehat h:=\mathsf S\circ\hbar.\]
\end{definition}

\begin{theorem}\label{thm:smooth} 
For smooth $w\in W$, we have $\widehat h(\widehat C_w)=P_{H_w}(u)$.
\end{theorem}

\begin{proof}
Let $H:=H_w$ and $I:=H^\perp\subset\fg^\vee$.
Write
\[K_H:
=\pi_{H*}\underline{\EE}_{X_H}[\dim \fg].\]
By definition $\widehat h(\widehat C_w)$ is the image of $\Hess(\mathcal F_w)$ in $\Rep W\otimes \cA$.
We first claim that
\[\Hess(\mathcal F_w)\simeq \Four_{\fg}\left(K_H\right)(b).\]
Indeed, since $w$ is smooth, the normal bundle of $\Delta$ in $\overline{\fX}_w$ is $G\times^B(H/\fb)$. Thus by Lemma~\ref{lem:Four.subbdl},
\[\begin{split}
	\mu_\Delta(\mathcal F_w)&\simeq \underline{\EE}_{G\times^B I}[2\dim \cB-2\ell(w)](-\ell(w))\simeq \Four_{\cB\times \fg/\cB}\left(\underline\EE_{G\times^BH}\right)[\dim \fg](b),
\end{split}
\]
where we regard complexes on $T_\cB^\vee$ as complexes on $\cB\times \fg^\vee$.
Pushing forward to $\fg^\vee$ and using the commutative diagram in Example~\ref{ex:Four.pullback}, we obtain the claim.

The decomposition theorem implies that $K_H$ is semisimple.
Since $\pi_H^{\mathrm{rs}}$ is smooth and proper, $R^i\pi_{H*}\underline{\EE}_{X_H}|_{\grs}\simeq R^i\pi_{H*}^{\mathrm{rs}}\underline{\EE}_{X_H^\rs}$ are local systems on $\grs$, and hence
\[{}^p\cH^i(K_H)|_{\grs}\simeq R^i\pi_{H*}^{\mathrm{rs}}\underline{\EE}_{X_H^\rs}[\dim\fg].\]
By Theorem~\ref{thm:Hess.monodromy}, the local system on the right corresponds to $V_i:=\coh^i(X_H(s))$.

Since $\Groth$ is the intermediate extension of the local system corresponding to the regular representation of $W$ and ${}^p\cH^i(K_H)$ is semisimple, restriction to $\grs$ gives
\[
\Hom_{\Perv_G(\fg)}
\big(\Groth,{}^p\cH^i(K_H)\big)
\simeq \Hom_{\EE[W]}(\EE[W],V_i)
\simeq V_i.
\]
The last isomorphism is evaluation at $1$ and is $W$-equivariant.
Since Fourier transform is a perverse $t$-exact equivalence, this is also isomorphic to $\mathsf S({}^p\cH^i(\Hess(\mathcal F_w)))$.
Therefore
\[
\widehat h(\widehat C_w)=\sum_i\big[\mathsf S({}^p\cH^i(\Hess(\mathcal F_w)))\big]u^{i/2}
=\sum_i[V_i]u^{i/2}=P_H(u).\qedhere
\]
\end{proof}

\subsection{Comparison with Lusztig's character functor}

We now compare the Hessenberg functor with Lusztig's character functor.
Recall that $\overline{Y}_w \subset G \times \cB$ is the subvariety consisting of pairs $(x,gB)$ such that $g^{-1}xg \in \overline{BwB}$.
Define
\[\sigma_G:G\times \cB \lra \cB\times \cB,\quad(x,gB)\mapsto (gB,xgB).\]
The map $\sigma_G$ is smooth and
$\overline{Y}_w=\sigma_G^{-1}(\overline{\fX}_w)$, so
$\sigma_G^*\IC_{\overline{\fX}_w}\simeq \IC_{\overline{Y}_w}$.
Lusztig's character functor is defined by the correspondence
\beq\label{eq:Lus.funct}
	\Ch:~\mathsf H_G\hookrightarrow \mathscr D_G(\cB\times \cB)\xrightarrow{\ \sigma_G^*\ }\mathscr D_G(G\times \cB)\xrightarrow{\ \pr_{1*}\ }\mathscr D_G(G),\eeq
where $\pr_1:G\times \cB\to G$ is the first projection.

There is a natural linear analogue of $\sigma_G$, which is the quotient vector bundle morphism
\[\sigma_{\fg}:\fg\times \cB\lra T_{\cB}\simeq G\times^B(\fg/\fb), \quad (v,gB)\mapsto [g,\Ad(g^{-1})v \bmod \fb].\]
Specialization relates these maps as follows.

\begin{proposition}\label{prop:lus.hess}
Consider the diagram
\[
\begin{tikzcd}
\mathscr D_G(\cB\times \cB)\arrow[r,"\sigma_G^*"]\arrow[d,"\sp_\Delta"']&
\mathscr D_G(G\times \cB)\arrow[r,"\pr_{1*}"]\arrow[d,"\sp_{e\times \cB}"']&
\mathscr D_G(G)\arrow[d,"\sp_e"]\\
\mathscr D_{G}(T_{\cB})^{\mon}\arrow[r,"\sigma_{\fg}^*"]&
\mathscr D_{G}(\fg\times \cB)^{\mon}\arrow[r,"\pr_{1*}"]&
\mathscr D_{G}(\fg)^{\mon}.
\end{tikzcd}
\]
Then both squares commute.
\end{proposition}

\begin{proof}
The map $\sigma_G$ is smooth and restricts to an isomorphism $\{e\}\times \cB \to \Delta$.
The commutativity of the left square follows from Proposition~\ref{prop:sp.smooth}.
The commutativity of the right square follows from the compatibility of specialization with proper direct image \cite[(SP2)]{verdier}.
\end{proof}
We next apply the Fourier transforms to the categories that appear in the lower row of the diagram in Proposition~\ref{prop:lus.hess}.
Using the identifications 
\[T_{\cB}\simeq G\times^B(\fg/\fb) \and \fg\times \cB\simeq G\times^B \fg,\] the map $\sigma_{\fg}:\fg\times \cB\to T_\cB$ may be viewed as the vector bundle morphism 
$G\times^B\fg\to G\times^B(\fg/\fb)$
induced by the canonical quotient $\fg\twoheadrightarrow \fg/\fb$.
The dual of the map $\sigma_{\fg}$ is the vector bundle map
${}^t\sigma_\fg:G\times^B(\fg/\fb)^\vee \to G\times^B\fg^\vee$, and it is naturally identified with the morphism
\[{}^t\sigma_\fg:T_\cB^\vee\lra  \fg^\vee\times \cB, \quad [g,v^*]\mapsto (g\cdot v^*,gB).\]

Let $\underline \eta: T_\cB^\vee\to \fg^\vee\times \cB$ denote the morphism induced by $\eta$; it sends $u=[g,v^*]$ to $(\eta(u),gB)$.
\begin{lemma}\label{lem:dual}
We have ${}^t\sigma_\fg=\underline \eta$.
\end{lemma}
\begin{proof}
	This follows immediately from the definition of the moment map $\eta$. 
\end{proof}

The compatibility of Fourier transform with pullback along vector bundle morphisms and with proper pushforward along projections gives the following.

\begin{proposition}\label{prop:fourier.comparison}
The diagram
\[
\begin{tikzcd}
\mathscr D_{G}(T_{\cB})^{\mon}
\arrow[r,"\sigma_{\fg}^!"]
\arrow[d,"\Four_{T_{\cB}}"']
&
\mathscr D_{G}(\fg\times \cB)^{\mon}
\arrow[r,"\pr_{1*}"]
\arrow[d,"\Four_{\fg\times\cB/\cB}"]
&
\mathscr D_{G}(\fg)^{\mon}
\arrow[d,"\Four_{\fg}"]
\\
\mathscr D_{G}(T_{\cB}^\vee)^{\mon}
\arrow[r,"\underline{\eta}_*\lbrack b\rbrack"]
&
\mathscr D_{G}(\fg^\vee\times \cB)^{\mon}
\arrow[r,"\pr_{1*}"]
&
\mathscr D_{G}(\fg^\vee)^{\mon}
\end{tikzcd}
\]
commutes, where $b:=\dim \fb$ is the relative dimension of $\sigma_\fg$. In particular,  the composition $\mathscr D_{G}(T_\cB)^{\mon}\to \mathscr D_{G}(\fg^\vee)^{\mon}$ factors through the full subcategory $\mathscr D_{G}(\mathcal N^\vee)^{\mon}\subset \mathscr D_{G}(\fg^\vee)^{\mon}$.
\end{proposition}

\begin{proof}
The left square commutes by \cite[Lemma~A.3.5]{bezrukavnikov2026affinespringerfibersmall} 
and Lemma~\ref{lem:dual}. 
For the right square, see Example~\ref{ex:Four.pullback}.
The last assertion holds since the image of $\eta:T_{\cB}^\vee\to \fg^\vee$ is $\mathcal N^\vee$.
\end{proof}

Combining Propositions~\ref{prop:lus.hess} and \ref{prop:fourier.comparison}, we obtain the following relation between Lusztig's character functor and the Hessenberg functor.

\begin{theorem}\label{thm:lus.hess}
The Hessenberg functor $\mathsf{Hess}$ in \eqref{eq:Hess.funct} is isomorphic, up to shift and Tate twist, to the composition of Lusztig's character functor $\Ch$ in \eqref{eq:Lus.funct} followed by the microlocalization $\mu_e$ at the identity.
More precisely, letting $b=\dim \fb$, we have
\[\mu_e\circ\Ch\simeq \mathsf{Hess}[-b](-b):\mathsf H_G\lra \mathscr D_G(\mathcal N^\vee).\]
\end{theorem}

\begin{proof}
This follows from Propositions~\ref{prop:lus.hess} and \ref{prop:fourier.comparison}, since $\sigma_\fg^!=\sigma_\fg^*[2b](b)$.
\end{proof}

\begin{corollary}[Trace property]\label{cor:trace}
For $F,G\in\mathsf H_G$, there is a natural isomorphism
\[\mathsf{Hess}(F\star G)\simeq \mathsf{Hess}(G\star F).\]
In particular, $\widehat h(ab)=\widehat h(ba)$ for all $a,b\in\hecke$.
\end{corollary}

\begin{proof}
Lusztig's character functor $\Ch$ satisfies the trace property \cite[Lemma~2.7.1]{trinh2021heckecategoryunipotentlocus}. Apply microlocalization at the identity and use Theorem~\ref{thm:lus.hess}.
\end{proof}

\subsection{Semisimplicity and the generalized Springer correspondence}
\label{ss:Hess.blocks}

In this subsection, we work over $\C$ with characteristic-zero coefficients and suppress Tate twists.
We choose a $G$-invariant nondegenerate bilinear form on $\fg$ and use it to identify $\mathcal N^\vee$ with $\mathcal N$.

We call a simple $G$-equivariant perverse sheaf on $\mathcal N$ \emph{principal} if it is a direct summand of $\Spr^\vee$, and \emph{nonprincipal} otherwise.
Let $\mathscr D_G(\mathcal N)_{\mathrm{pr}}$ and $\mathscr D_G(\mathcal N)_{\mathrm{np}}$ be the full triangulated subcategories generated by the principal and nonprincipal simple perverse sheaves, respectively.
Lusztig's generalized Springer theory gives the orthogonal decomposition
\[
\mathscr D_G(\mathcal N)\simeq
\mathscr D_G(\mathcal N)_{\mathrm{pr}}
\oplus\mathscr D_G(\mathcal N)_{\mathrm{np}}
\]
(see \cite[Theorem~3.5]{rider-russell}).
Composing $\Hess$ with the projections onto these subcategories, we obtain its principal and nonprincipal parts, $\Hess_{\mathrm{pr}}$ and $\Hess_{\mathrm{np}}$.

\begin{proposition}\label{prop:Hess.semisimple}\label{cor:Hess.blocks}
The Hessenberg functor takes values in $\Semis_G(\mathcal N)$. Moreover, 
\begin{equation}\label{eq:Hess.blocks}
\Hess\simeq\Hess_{\mathrm{pr}}\oplus\Hess_{\mathrm{np}}.
\end{equation}
For every $F\in\mathsf H_G$ and every $i$, we have
$\mathsf S({}^p\cH^i(\Hess_{\mathrm{np}}(F)))=0$.
\end{proposition}

\begin{proof}
Every object $F$ of $\mathsf H_G$ is a direct sum of shifts of $\IC_{\overline \fX_w}$ and their pullbacks by $\sigma_G$ are $\IC_{\overline Y_w}$.
The decomposition theorem therefore gives $\Ch(F)\simeq\bigoplus_i A_i[n_i]$ for simple perverse sheaves $A_i$. 
By Remark~\ref{rem:mu.texact}, each $\mu_eA_i$ is perverse and its support is contained in $\mathcal N$ by Theorem~\ref{thm:lus.hess}.
In particular, $\mu_eA_i\in \Perv_G(\mathcal N)$. 
Since 
$\Perv_G(\mathcal N)$ is semisimple (\cite[Section~9.3, Theorem~61]{CEE}),
applying $\mu_e$ to the decomposition of $\Ch(F)$ proves the first assertion.

The two projections preserve the perverse $t$-structure and send simple perverse sheaves to themselves or zero.
Hence they give \eqref{eq:Hess.blocks}, and orthogonality gives the claimed vanishing.
\end{proof}

\begin{remark}
	The Hessenberg character $\widehat h$ determines the principal part of $\Hess(\mathcal F_w)$. 
	Moreover, $\Hess_{\mathrm{np}}(\mathcal F_w)=0$ for smooth $w$ by the result of Precup and Sommers \cite[Theorem~1.1]{PS}, and for every $w$ when $G=\mathrm{GL}_n$, since every simple $G$-equivariant perverse sheaf on $\mathcal N$ occurs in the Springer correspondence.
	However, it can be nonzero; see Proposition~\ref{prop:G2.nonprincipal}. 
\end{remark}

\section{Lusztig's character formula}
\label{s lusztig character}
We recall unipotent character sheaves and Lusztig's
character $\varepsilon$.
We then use Theorem~\ref{thm:lus.hess} and Lusztig's multiplicity formula to compute the Hessenberg character.
To obtain the formula for regular semisimple Lusztig varieties, we identify the stalks of $\overline K_w$ over $G^{\mathrm{rs}}$ with the intersection cohomology of the fibers and show that the monodromy action factors through $W$.
These arguments prove Corollary~\ref{cor:lus.formula}.

For background on character sheaves, see \cite{Lus85,Lus85ii,Lus86iv,Lus86v,MarSpr89,Sho88}.
We keep the notation of Subsection~\ref{ss notations}.

\subsection{Character sheaves}
Recall that the $G$-orbits on $\mathcal B \times \mathcal B$ are in bijection with the Weyl group $W$.
For $w\in W$, let
$\fX_w := G \cdot (B, wB)$
denote the corresponding $G$-orbit. 
Define
\[Y_w := \{(g, x) \in G \times \mathcal B : (x, gx) \in \fX_w\},\]
and let $\pi_w : Y_w \to G$ be the first projection.
Set
\[K_w := \pi_{w!} \underline{\EE} \;\in \mathscr D_G(G).\]

\begin{definition}
	An irreducible perverse sheaf $A$ on $G$ is called a \emph{unipotent character sheaf} if $A$ is a constituent of ${}^p\cH^i(K_w)$ for some $w \in W$ and some $i\in\ZZ$.
	We denote by $\CS_1(G)$ the set of isomorphism classes of unipotent character sheaves. 
\end{definition}

Similarly, let $\overline{\fX}_w$ denote the closure of $\fX_w$ in $\cB\times \cB$, and  set
\[\overline{Y}_w := \{(g, x) \in G \times \mathcal B : (x, gx) \in \overline{\fX}_w\}.\]
Let $\overline{\pi}_w : \overline{Y}_w \to G$ be the first projection, and define
\[\overline{K}_w := \overline{\pi}_{w*} \mathrm{IC}_{\overline{Y}_w} \;\in \mathscr D_G(G).\]

The following results are due to Lusztig.

\begin{proposition}[{\cite[Propositions~12.7 and~12.8]{Lus85ii}}]\label{prop:Kbarw}
	For $w\in W$, the complex $\overline{K}_w$ is semisimple. Moreover, an irreducible perverse sheaf $A$ on $G$ is a unipotent character sheaf if and only if $A$ occurs as a direct summand of ${}^p\cH^i(\overline{K}_w)$ for some $w \in W$ and some $i\in \ZZ$.
\end{proposition}

\subsection{Lusztig's character and multiplicity formula}
\label{s families}\label{ss:lusztig.character}

Lusztig constructs an $\cA$-linear map
\[
\varepsilon:\hecke\lra\mathscr K_\oplus(\Perv_G(G))\otimes\cA
\]
which we call \emph{Lusztig's character}. It satisfies
\beq\label{eq:epsilon.KL}
\varepsilon(\widehat C_w) = u^{-\dim G/2}\sum_i(-1)^i[{}^p\cH^i(\overline K_w)]u^{i/2}
\eeq
for every $w\in W$
\cite[Proposition~12.10]{Lus85ii}.
It also satisfies the trace property
\[
\varepsilon(ab)=\varepsilon(ba) \quad\text{for }a,b\in\hecke
\]
\cite[Proposition~6.2]{Lus85}.
These statements hold over the algebraically closed
base field $\Bbbk$.

For the reader's convenience, we briefly recall Lusztig's original construction of $\varepsilon$ in the unipotent case and how he proves that it is well defined and satisfies the trace property.
Let $\underline w=(s_1,\ldots,s_r)$ be a finite sequence of simple reflections, not necessarily reduced. Set
\[
\overline Y_{\underline w}
=\left\{(g,x_0,\ldots,x_r)\in G\times\cB^{r+1}:
(x_{j-1},x_j)\in\overline{\fX}_{s_j}\;\text{ for }1\leq j\leq r \;\text{ and }\;x_r=gx_0
\right\}.\]
For the empty sequence, $\overline Y_{\underline w}=Y_e$.
The variety $\overline Y_{\underline w}$ is smooth, and the projection
$\overline\pi_{\underline w}:\overline Y_{\underline w}\to G$
is proper \cite[Sections~2.6--2.8]{Lus85}.
Let
\[
\overline K_{\underline w} := \overline\pi_{\underline w*}\underline\EE_{\overline Y_{\underline w}}.
\]
The complex $\overline K_{\underline w}$ is semisimple, and the simple constituents of its perverse cohomology are unipotent character sheaves
\cite[Propositions~2.9 and~2.17]{Lus85}.

In his original construction, Lusztig defines $\varepsilon$ by prescribing
\beq\label{eq:epsilon.BS}
\varepsilon\big((1+T_{s_1})\cdots(1+T_{s_r})\big) = u^{-\dim G/2}\sum_i(-1)^i[{}^p\cH^i(\overline K_{\underline w})]u^{i/2}
\eeq
for every such sequence.
Then he proves that these values define a unique $\cA$-linear map and that this map satisfies the trace property
\cite[Section~6.1 and Proposition~6.2]{Lus85}.

To prove well-definedness and the trace property, Lusztig first reduces to the case of an algebraic closure of a finite field, using \cite[Section~6]{BBD}; see the proof of \cite[Proposition~6.2]{Lus85}.
In this case, we equip $K_w$ and $\overline K_{\underline w}$ with their natural mixed structures over a sufficiently large finite field.
Write $\Gr^{\W}_j{}^p\cH^i(K_w)$ for the $j$-th graded piece of the weight filtration on ${}^p\cH^i(K_w)$.

Lusztig defines an $\cA$-linear map $\varepsilon': \hecke\to\mathscr K(\Perv_G(G))\otimes\cA$ on the standard basis by
\beq\label{eq:epsilon.prime}
\varepsilon'(T_w) = u^{-\dim G/2}\sum_{i,j}(-1)^i\{\Gr^{\W}_j{}^p\cH^i(K_w)\}u^{j/2}.
\eeq
Here the braces denote classes in the usual Grothendieck group after forgetting the Frobenius structures.
Since the elements $T_w$ form an $\cA$-basis of $\hecke$, this defines $\varepsilon'$ uniquely.

He then proves that, for every sequence $\underline w=(s_1,\ldots,s_r)$,
\[
\varepsilon'\big((1+T_{s_1})\cdots(1+T_{s_r})\big)
=u^{-\dim G/2}\sum_i(-1)^i
\{{}^p\cH^i(\overline K_{\underline w})\}u^{i/2}
\]
\cite[(6.3.1)]{Lus85}.
Hence, after composing $\varepsilon$ as defined in \eqref{eq:epsilon.BS} with the natural map from the split Grothendieck group to the usual Grothendieck group, we obtain $\varepsilon'$, which is well defined and satisfies the trace property \cite[Section~6.4]{Lus85}.

This natural map is injective on the subgroup spanned by unipotent character sheaves, since they are simple perverse sheaves.
Moreover, each ${}^p\cH^i(\overline K_{\underline w})$ is a direct sum of unipotent character sheaves, and the products in \eqref{eq:epsilon.BS} span $\hecke$ over $\cA$, since $T_s=(1+T_s)-1$. 
It follows that \eqref{eq:epsilon.BS} defines a unique well-defined
$\cA$-linear map $\varepsilon$ satisfying the trace property.

\begin{remark}
Each $\Gr^{\W}_j{}^p\cH^i(K_w)$ is pure and is semisimple after forgetting its Frobenius structure \cite[Theorem~5.3.8]{BBD}.
Its simple constituents are unipotent character sheaves, since it is a subquotient of ${}^p\cH^i(K_w)$.
The injectivity above therefore shows that \eqref{eq:epsilon.prime}, with $[\Gr^{\W}_j{}^p\cH^i(K_w)]$ in place of $\{\Gr^{\W}_j{}^p\cH^i(K_w)\}$, gives $\varepsilon(T_w)$ in the split Grothendieck group.
\end{remark}

We identify a unipotent character sheaf with its class in the split Grothendieck group.
After extending scalars to $\EE(u^{1/2})$, the trace property gives 
\[
\varepsilon(a)
=\sum_{A\in\CS_1(G)}\sum_{\chi\in\Irr W}
c_{A,\chi}\,\widetilde\chi(a)A,
\quad\text{for some }~ c_{A,\chi}\in\EE(u^{1/2}).
\]
The coefficients are uniquely determined.
We now recall their description by Lusztig's Fourier pairing and deduce the formula for the graded multiplicities without alternating signs.

The set $\Irr W$ is partitioned into \emph{families}.
To each family $\mathscr F$, Lusztig associates a finite group $\mathscr G_{\mathscr F}$ and an embedding $\mathscr F\hookrightarrow\mathscr M(\mathscr G_{\mathscr F})$  written as $\chi\mapsto m_{\chi}$; see \cite[Chapter~4]{Lus84}.
For a finite group $\mathscr G$, the set $\mathscr M(\mathscr G)$ consists of pairs $(x,\tau)$, where $x\in\mathscr G$ and $\tau$ is an irreducible representation of the centralizer $Z(x)$ of $x$ in $\mathscr G$, up to simultaneous conjugation.

Thus we have an embedding
\begin{equation}\label{eq:family.IrrW}
\Irr W\hookrightarrow
\coprod_{\mathscr F}\mathscr M(\mathscr G_{\mathscr F}),
\quad \chi\mapsto m_\chi.
\end{equation}
The Fourier pairing on $\mathscr M(\mathscr G)$ is
\[
\{(x,\tau),(y,\sigma)\}
=\frac{1}{|Z(x)||Z(y)|}
\sum_{\substack{g\in\mathscr G\\
x(gyg^{-1})=(gyg^{-1})x}}
\Tr(g^{-1}x^{-1}g,\sigma)\Tr(gyg^{-1},\tau).
\]
Here $\Tr(g,\tau)$ denotes the trace of $g$ in the representation $\tau$.
We extend the pairing by zero between distinct families.
For $\mathscr G\in\{S_2,S_3,S_4,S_5\}$, the matrices of this pairing are given in \cite[Section~4.15]{Lus84}.

Lusztig's classification gives a bijection
\beq\label{eq:family}
\CS_1(G)\longleftrightarrow\coprod_{\mathscr F}\mathscr M(\mathscr G_{\mathscr F}), \quad A\leftrightarrow\overline m_A
\eeq
\cite[Theorem~23.1]{Lus86v}.
We use Lusztig's parametrization and write
\[
\{A,\chi\}:=\{\overline m_A,m_\chi\}
\and
\{\psi,\chi\}:=\{m_{\psi},m_\chi\}.
\]
The restriction on the characteristic in \cite[(23.0.1)]{Lus86v}, used there to ensure cleanness, is removed by \cite[Theorem~1.2 and Section~3.10]{Lus12}.
The parametrization and the formula below therefore apply over an algebraically closed field of arbitrary characteristic.

\begin{theorem}[Lusztig]\label{p coefficient}
For every $w\in W$,
\begin{equation}\label{e 8}
\sum_i[{}^p\cH^i(\overline K_w)]u^{i/2}
=u^{\dim G/2}
\sum_{A\in\CS_1(G)}\sum_{\chi\in\Irr W}
\{A,\chi\}\widetilde\chi(\widehat C_w)A.
\end{equation}
\end{theorem}
\begin{proof}
Lusztig's multiplicity formula \cite[Corollary~14.11]{Lus85ii} and \cite[Theorem~23.1]{Lus86v} identify
\[
c_{A,\chi} = (-1)^{\dim G}\widehat\varepsilon_A\{A,\chi\},
\]
where $\widehat\varepsilon_A=(-1)^{\codim(\supp A)}$.
Combining this with \eqref{eq:epsilon.KL} gives
\begin{equation}\label{e 7}
\sum_i(-1)^i[{}^p\cH^i(\overline K_w)]u^{i/2}
=(-1)^{\dim G}u^{\dim G/2}
\sum_{A\in\CS_1(G)}\sum_{\chi\in\Irr W}
\widehat\varepsilon_A\{A,\chi\}
\widetilde\chi(\widehat C_w)A.
\end{equation}

By Lusztig's parity property \cite[Sections~15.11--15.13]{Lus85ii}, \cite[Theorem~23.1]{Lus86v}, \cite[Section~3.10]{Lus12},
if $A$ occurs in ${}^p\cH^i(\overline K_w)$, then $(-1)^{i+\dim G}=\widehat\varepsilon_A$.
Comparing the coefficient of each $A$ in \eqref{e 7} gives \eqref{e 8}.
\end{proof}

\subsection{Unipotent character sheaves with full support and representations of $W$}\label{ss:fullsupp}
We record the following standard consequence of Lusztig's classification of character sheaves. 

Write $Y_e^{\mathrm{rs}}:=\pi_e^{-1}(\Grs)$ and let $\pi_e^{\mathrm{rs}}:Y_e^{\mathrm{rs}}\to\Grs$ be the restriction of $\pi_e$.
This is the group version of the Grothendieck--Springer $W$-torsor, and it induces a surjection
\beq\label{eq:quotient.W}
	\pi_1(\Grs,\mathbf{s})\twoheadrightarrow W,
\eeq
well defined up to an inner automorphism of $W$.
\begin{lemma}\label{lem:fullsupp.CS}
Let $j:\Grs\hookrightarrow G$. For $\chi\in\Irr W$, let $\rho_\chi: W\to\mathrm{GL}(V_\chi)$ be an irreducible representation with character $\chi$, and consider the associated local system
$L_\chi:=Y_e^{\mathrm{rs}}\times^W V_\chi$
on $\Grs$.
Then the unipotent character sheaves on $G$ having support equal to $G$ are precisely
\[\mathbb{IC}(G,L_\chi)= j_{!*}\left(L_\chi[\dim G]\right)\]
for $\chi\in\Irr W$.
Moreover, the monodromy representation of $L_\chi$ at a closed point $\mathbf{s}\in \Grs$ is the composite
$\pi_1(\Grs,\mathbf{s})\twoheadrightarrow W\xrightarrow{\rho_\chi}\mathrm{GL}(V_\chi)$,
where the first map is \eqref{eq:quotient.W}. 
\end{lemma}
\begin{proof}
The Grothendieck--Springer morphism $\pi_e:Y_e\to G$ is small. Hence
\[K_e[\dim G]\simeq j_{!*}\left(\pi_{e*}^{\mathrm{rs}}\underline \EE_{Y_e^{\mathrm{rs}}}[\dim G]\right).\]
The $W$-action on the $W$-torsor $\pi_e^{\mathrm{rs}}$ gives an isotypic decomposition
\[K_e[\dim G]\simeq\bigoplus_{\chi\in\Irr W}V_\chi^\vee\otimes\mathbb{IC}(G,L_\chi).\]
Thus $\ICperv(G,L_\chi)$ are precisely the simple summands of $K_e[\dim G]$.

Conversely, let $A\in\CS_1(G)$ have support equal to $G$.
By Lusztig's parabolic induction theorem for character sheaves
\cite[Theorem~4.4(a)]{Lus85}, there exist a parabolic subgroup
$P\subset G$ with Levi subgroup $L$ and a cuspidal character sheaf $A_1$ on $L$ such that $A$ is a direct summand of $\ind_P^G(A_1)$, where the parabolic induction $\ind_P^G(-)$ is defined in \cite[Section~4.1]{Lus85}.

By Lusztig's support formula~\cite[(4.3.1)]{Lus85}, $A$ can have full support only when $L=T$.
Since $A$ is unipotent, \cite[Proposition~4.8(b)]{Lus85} and \cite[Proposition~11.2(c)]{Lus85ii} imply that
\[A_1=\underline \EE_T[\dim T].\]

When $L=T$ and $P=B$, 
by the definition of parabolic induction,
\[
	\ind_B^GA_1 
	\simeq \pi_{e*}\underline \EE_{Y_e}[\dim G]
	=K_e[\dim G].
\]
Thus $A$ is a simple summand of $K_e[\dim G]$, and hence $A=\mathbb{IC}(G,L_\chi)$ for some $\chi\in \Irr W$.

The final assertion follows directly from the definition of the local system $L_\chi$.
\end{proof}

By Lusztig's parametrizations~\eqref{eq:family.IrrW} and~\eqref{eq:family}, $\overline m_{\ICperv(G,L_\chi)}=m_\chi$;
see \cite[Theorem~23.1(c)]{Lus86v}.

\begin{lemma}\label{lem:deform.rs}
Let $M:=M_{\{e\}/G}$ be the total space of the deformation of $G$ to the normal cone $C_{\{e\}/G}\simeq\fg$.
The subset
\[
M^{\mathrm{rs}}:=(\Grs\times\Gm)\sqcup\grs
\]
is open in $M$.
Moreover, the Grothendieck--Springer $W$-torsors over $\Grs$ and $\grs$ fit into a single $W$-torsor over $M^{\mathrm{rs}}$.
Consequently, for every $\chi\in\Irr W$, 
$\big(\sp_e\ICperv(G,L_\chi)\big)|_{\grs} \simeq \mathcal L_\chi[\dim\fg]$.

\end{lemma}

\begin{proof}
Conjugation on $G$ extends to an action on $M$.
The inclusion $B\subset G$ induces a closed immersion
$M_{\{e\}/B}\hookrightarrow M$, and hence a morphism
\[
\pi_M:\widetilde{M}:=G\times^B M_{\{e\}/B}\lra M,\quad[g,z]\mapsto g\cdot z,
\]
whose fibers over $\Gm$ and $0$ are $\pi_e:G\times^B B\to G$ and $\pi_\fb:G\times^B\fb\to\fg$, respectively.
The map $\widetilde M\to M\times\cB$, $[g,z]\mapsto(g\cdot z,gB)$, is a closed immersion.
Since $\cB$ is proper, $\pi_M$ is proper.

Let $\widetilde{M}_{\text{\'et}}\subset\widetilde{M}$ be the \'etale locus of $\pi_M$.
Both $\widetilde{M}$ and $M$ are smooth over $\A^1$, so $\pi_M$ is \'etale at a point if and only if its restriction to the corresponding fiber is \'etale there.
At $[e,b]\in G\times^B B$ and $[e,x]\in G\times^B\fb$, the latter condition is the invertibility of
\[
\Ad(b)-\id:\fg/\fb\lra\fg/\fb
\and
\ad(x):\fg/\fb\lra\fg/\fb,
\]
respectively.
These conditions are equivalent to $b\in\Grs$ and $x\in\grs$, respectively.
Hence,
\[\widetilde{M}_{\text{\'et}}=\pi_M^{-1}(M^\rs) \and \pi_M(\widetilde{M}_{\text{\'et}})=M^\rs.\]
Since \'etale morphisms are open, $M^{\mathrm{rs}}=\pi_M(\widetilde{M}_{\text{\'et}})$ is open in $M$.

Moreover, the restriction $\pi_M|_{\widetilde{M}_{\text{\'et}}}$ is proper and \'etale, hence finite \'etale.
To identify its $W$-action, let $\mathcal T^{\mathrm{rs}}\subset M_{\{e\}/T}$ be the inverse image of $M^{\mathrm{rs}}$.
The inclusion $T\subset B$ induces 
\[
G\times^T\mathcal T^{\mathrm{rs}}\lra {\widetilde{M}_{\text{\'et}}}.
\]
Let $T^{\mathrm{rs}}=T\cap G^{\mathrm{rs}}$.
On every geometric fiber over $\A^1$, this is the standard isomorphism $G\times^T T^{\mathrm{rs}}\simeq Y_e^{\mathrm{rs}}$ or $G\times^T\ft^{\mathrm{rs}}\simeq\tgrs$.
Since both the source and target are smooth over $\A^1$ and the morphism is an isomorphism on every geometric fiber, it is 
an isomorphism.

For $w\in W$, with a representative $\dot w\in N_G(T)$,
$[g,z]\cdot w=[g\dot w,\dot w^{-1}\cdot z]$
defines a $W$-action on ${\widetilde{M}_{\text{\'et}}}$ over $M^{\mathrm{rs}}$.
It restricts to the $W$-actions on the fibers, so $\pi_M|_{\widetilde{M}_{\text{\'et}}}$ is the required $W$-torsor.

The local system on $M^{\mathrm{rs}}$ associated with $\chi$ restricts to $L_\chi\boxtimes\underline\EE_{\Gm}$ over $\Gm$ and to $\mathcal L_\chi$ over zero.
Since $M^{\mathrm{rs}}$ is smooth over $\A^1$, the last assertion follows from \cite[2.1.5]{SGA7-nearby}.
\end{proof}

\begin{proof}[Proof of Corollary~\ref{cor:hess.formula}]

By Proposition~\ref{prop:lus.hess},
$\sp_e(K_e[\dim G])\simeq\mathrm{Groth}$.
Thus $\sp_e\ICperv(G,L_\chi)$ is a direct summand of $\mathrm{Groth}$.
Since $\mathrm{Groth}$ is the intermediate extension of the local system corresponding to the regular representation, by Lemma~\ref{lem:deform.rs}, restriction gives
\[\Hom_{\Perv_G(\fg)}\big(\mathrm{Groth},\sp_e\ICperv(G,L_\chi)\big)\simeq V_\chi\]
as $W$-representations.
If $A\in\CS_1(G)$ has proper support, then $\sp_eA=0$ when $e\notin\mathrm{Supp}A$. Otherwise, its support is contained in the tangent cone $C_e(\mathrm{Supp}A)$, which has dimension $\dim\mathrm{Supp}A<\dim G$. Thus $\sp_eA$ also has proper support, and
$\Hom_{\Perv_G(\fg)}(\mathrm{Groth},\sp_eA)=0$.
Thus, using Fourier transform and Theorem~\ref{thm:lus.hess}, we obtain
\[
\widehat h(\widehat C_w)
=\sum_{\chi\in\Irr W}\sum_i
\dim\Hom_{\Perv_G(G)}\bigl(\ICperv(G,L_\chi),
{}^p\cH^{i+\dim G}(\overline K_w)\bigr)u^{i/2}\chi.
\]

Applying \eqref{e 8} and using $\overline m_{\ICperv(G,L_\chi)}=m_\chi$ gives
\[
\widehat h(\widehat C_w)=\sum_{\chi,\psi\in\Irr W}\{\chi,\psi\}\widetilde{\psi}(\widehat C_w)\chi.\qedhere
\]
\end{proof}

\subsection{Restriction to a regular semisimple fiber}
\label{ss:rs.fiber}

We show that, over any algebraically closed field $\Bbbk$, $\cH^i(\overline K_w)|_{\Grs}$ is a local system for every $i\in \ZZ$, and identify its stalk at $\mathbf{s}\in \Grs$ with the intersection cohomology $\IH^i(\overline Y_w(\mathbf{s}))$ of the corresponding fiber $\overline Y_w(\mathbf{s})$. 

Over $\C$, Abreu--Nigro obtain these conclusions from the topological fibrations in \cite[Propositions~8.2 and~8.3]{AN}. 
Their proof uses Whitney stratifications, Thom's first isotopy lemma, and Ehresmann's fibration theorem.
We provide an algebraic proof valid over $\Bbbk$.

We then prove that the resulting monodromy representations factor through the quotient
\[
	\pi_1(\Grs,\mathbf{s})\twoheadrightarrow W
\]
induced by the Grothendieck--Springer $W$-torsor $\pi_e^{\mathrm{rs}}:Y_e^{\mathrm{rs}}\to \Grs$. 
Over $\C$, analogous results for regular semisimple Hessenberg varieties \cite{BrCh18,balibanu-crooks} and Lusztig varieties \cite{AN} are proved using torus-equivariant cohomology and localization. Here we instead use Lemma~\ref{lem:fullsupp.CS}.

\begin{proposition}\label{prop:fiber.IC}
	Let $w\in W$. For every $i\in \ZZ$, $\cH^i(\overline K_w)|_{\Grs}$ is a local system, and it vanishes when $i$ is odd. For every  $\mathbf{s}\in \Grs$, there is a canonical isomorphism 
	$\cH^i(\overline K_w)_{\mathbf{s}}\simeq \IH^i(\overline Y_w(\mathbf{s}))$.
\end{proposition}
\begin{proof}
	Choose a reduced expression $\underline w$ of $w$.
	Let $\rho_{\underline w}:\overline Y_{\underline w}\to \overline Y_w$ denote the Bott--Samelson resolution. It is proper and birational, and it is fiberwise a resolution of singularities over $\Grs$. In particular, $\overline \pi_{\underline w}:\overline Y_{\underline w}\to G$ is proper, and smooth over $\Grs$.
	
	By the decomposition theorem, $\IC_{\overline Y_w}$ is a direct summand of $\rho_{\underline w *}\underline \EE_{\overline Y_{\underline w}}$. Thus $\mathcal H^i(\overline K_w)=\mathcal H^i(\overline\pi_{w*}\IC_{\overline Y_w})$ is a direct summand of $\mathcal H^i(\overline \pi_{\underline w*}\underline \EE)$. 
	Since $\overline \pi_{\underline w}$ is smooth and proper over $\Grs$, the cohomology sheaves of $(\overline \pi_{\underline w*}\underline \EE_{\overline Y_{\underline w}})|_{\Grs}$ are local systems. 
	Since a direct summand of a local system is a local system, $\mathcal H^i(\overline K_w)|_{\Grs}$ is also a local system.
	
	For the vanishing in odd degrees, fix $\mathbf s\in\Grs$ and let $T_{\mathbf s}:=Z_G(\mathbf s)^\circ$.
	The smooth projective variety $\overline Y_{\underline w}(\mathbf s)$ has finitely many $T_{\mathbf s}$-fixed points, since it is a $T_{\mathbf s}$-stable subvariety of a product of copies of $\cB$.
	Thus the Bia{\l}ynicki--Birula theorem gives an affine paving, so its odd cohomology vanishes.
	By the above observation and proper base change, $\cH^i(\overline K_w)|_{\Grs}=0$ for odd $i$.
	
	To prove the last statement, we consider the morphism
	\[\varphi_{\mathbf{s}}:G\times\cB\lra \cB\times \cB, \quad (g,x)\mapsto (gx,g\mathbf{s}x).\]
	The restriction of $\varphi_{\mathbf{s}}$ to the identity section $x\mapsto (e,x)$ is the graph $\cB\to \cB\times \cB$, $x\mapsto (x,\mathbf{s}x)$.
	By \cite[Lemma~1.1 and the following remark]{Lus}, this graph is transverse to every diagonal $G$-orbit $\fX_w$ in $\cB\times \cB$. 
	It follows that $d\varphi_{\mathbf{s}}$ is surjective along the identity section. 
	Since $\varphi_{\mathbf{s}}$ is $G$-equivariant, its differential is surjective everywhere. Hence $\varphi_{\mathbf{s}}$ is smooth.
	
	The inverse image of $\overline{\fX}_w$ under $\varphi_{\mathbf s}$ is $G\times\overline Y_w(\mathbf s)$.
	Let $e_{\mathbf s}:\overline Y_w(\mathbf s)\to G\times\overline Y_w(\mathbf s)$ be the identity section, and let $i_{\mathbf s}:\overline Y_w(\mathbf s)\hookrightarrow\overline Y_w$ be the fiber inclusion. 
	Restricting the maps to $\overline Y_w(\mathbf s)$, we have $\sigma_G\circ i_{\mathbf s}=\varphi_{\mathbf s}\circ e_{\mathbf s}$.
	Smooth pullback and the product formula therefore give
	\[
	i_{\mathbf s}^*\IC_{\overline Y_w}
	\simeq i_{\mathbf s}^*\sigma_G^*\IC_{\overline{\fX}_w}
	\simeq e_{\mathbf s}^*\varphi_{\mathbf s}^*\IC_{\overline{\fX}_w}
	\simeq e_{\mathbf s}^*\IC_{G\times\overline Y_w(\mathbf s)}
	\simeq\IC_{\overline Y_w(\mathbf s)}.
	\]
	By proper base change, we finally have
	\[
	\cH^i(\overline K_w)_{\mathbf s}
	\simeq\mathbb H^i\bigl(\overline Y_w(\mathbf s),
	i_{\mathbf s}^*\IC_{\overline Y_w}\bigr)
	\simeq\IH^i(\overline Y_w(\mathbf s)).\qedhere
	\]
\end{proof}

Via the identification in Proposition~\ref{prop:fiber.IC}, the monodromy of the local system $\cH^i(\overline K_w)|_{\Grs}$ induces an action of $\pi_1(\Grs,\mathbf{s})$ on $\IH^i(\overline Y_w(\mathbf{s}))$.

\smallskip

\begin{proposition}\label{prop:mono.fact}
	The monodromy action of $\pi_1(\Grs,\mathbf{s})$ on $\IH^*(\overline Y_w(\mathbf{s}))$ factors through \eqref{eq:quotient.W}.
\end{proposition}

\begin{proof}
By Proposition~\ref{prop:Kbarw}, $\overline K_w$ is a direct sum of shifts of unipotent character sheaves. 
Since $\cH^i(\overline K_w)|_{\Grs}$ is a local system, every summand of $\overline K_w$ with proper support restricts trivially to $\Grs$.
Hence, by Lemma~\ref{lem:fullsupp.CS}, each $\cH^i(\overline K_w)|_{\Grs}$ is a finite direct sum of local systems of the form $L_\chi$, with $\chi\in\Irr W$.
By the last assertion of Lemma~\ref{lem:fullsupp.CS}, the monodromy of each $L_\chi$ factors through \eqref{eq:quotient.W}, thus the monodromy on $\cH^i(\overline K_w)_{\mathbf{s}}$ factors through \eqref{eq:quotient.W}.
The claim now follows from the identification $\cH^i(\overline K_w)_{\mathbf{s}}\simeq \IH^i(\overline Y_w(\mathbf{s}))$ in Proposition~\ref{prop:fiber.IC}.
\end{proof}

\begin{proof}[Proof of Corollary~\ref{cor:lus.formula}]
Since $\Grs$ is smooth and $\cH^i(\overline K_w)|_{\Grs}$ is a local system, we have
\[{}^p\cH^{i+\dim G}(\overline K_w)_{\mathbf s}\simeq\cH^{i}(\overline K_w)_{\mathbf s}[\dim G]\simeq\IH^i(\overline{Y}_w(\mathbf s))[\dim G],\]
where the second isomorphism follows from Proposition~\ref{prop:fiber.IC}.
By Proposition~\ref{prop:mono.fact}, the monodromy on $\IH^i(\overline Y_w(\mathbf{s}))$ factors through $W$.
Taking stalks at $\mathbf{s}$ in \eqref{e 8}, only the character sheaves with full support contribute. By
Lemma~\ref{lem:fullsupp.CS}, the stalk $\mathbb{IC}(G,L_{\psi})_{\mathbf{s}}$ is $V_{\psi}[\dim G]$. 
We therefore obtain an equality in $\Rep W\otimes\cA$
\[\sum_i 
\IH^{i}(\overline Y_w(\mathbf s)) u^{i/2} =  \sum_{\psi \in \Irr W}\sum_{\chi \in \Irr W}   \{\psi,\chi\}\widetilde\chi(  \widehat C _w  ) \psi.\]
The vanishing in odd degrees in Proposition~\ref{prop:fiber.IC} gives the stated formula.
\end{proof}

\section{Hessenberg schemes defined by cones}
\label{s:nonlinear}
We study Hessenberg schemes defined by nonempty pure-dimensional $B$-stable closed cones in $\fg/\fb$.
For the morphisms $\pi_C:X_C\to\fg$ constructed below, we show that $\pi_{C*}\IC_{X_C}$ is a direct sum of shifts of nilpotent character sheaves on $\fg$.
This is analogous to the decomposition of $\overline K_w$ into shifts of unipotent character sheaves on $G$.
From this, we deduce that the cohomology sheaves of $\pi_{C*}\IC_{X_C}$ are locally constant over $\grs$, and their monodromy factors through $W$, giving a $W$-action on the intersection cohomology of each regular semisimple fiber.

We then apply this construction to the tangent cones of Schubert varieties at $eB$.
When specialization at $eB$ takes the intersection cohomology complex of $\overline X_w$ to that of its tangent cone, the resulting direct image is the specialization of $\overline K_w$ at $e$, and the graded character of each regular semisimple fiber is $\widehat h(\widehat C_w)$.
We give a sufficient local condition and examples.

\subsection{Character sheaves and monodromy}
\label{ss:nonlinear.support}

Let $\overline C\subset\fg/\fb$ be a nonempty pure-dimensional $B$-stable closed cone, and let $C\subset\fg$ be its inverse image under $\fg\to\fg/\fb$.
Set
\[
X_{\overline C}=G\times^B\overline C\subset T_{\cB} \and X_{C}=G\times^BC.
\]
Under the isomorphism $G\times^B\fg\simeq\fg\times\cB$, we have $X_{C}=\sigma_\fg^{-1}(X_{\overline C})$.
The projection to $\fg$ restricts to a proper morphism $\pi_{C}:X_{C}\to\fg$.
Its fibers
\[X_{C}(s):=\pi_{C}^{-1}(s)\]
are Hessenberg schemes in the sense of \cite[Definition~5.1]{BHL25}.

We call a simple $G$-equivariant perverse sheaf $A$ on $\fg$, monodromic for scalar multiplication, a \emph{nilpotent character sheaf} if the support of its Fourier transform is contained in $\mathcal N^\vee$:
\[
\supp\Four_\fg(A)\subseteq\mathcal N^\vee.
\]
Equivalently, $\Four_\fg(A)$ is the intersection cohomology extension of an irreducible $G$-equivariant local system on a nilpotent coadjoint orbit.
These are Lie algebra analogues of unipotent character sheaves on $G$; see \cite[Section~2]{AHJR-overview}.

Over $\C$, the singular support of a character sheaf on $G$ lies in $G\times\mathcal N^\vee\subset T_G^\vee$ \cite[Theorem~4.4]{MV88}.
Since $A$ is monodromic for scalar multiplication, the support of $\Four_\fg(A)$ is the projection of its singular support to $\fg^\vee$ \cite[Theorem~5.5.5]{KS94}.
Thus the defining condition above is equivalent to the analogous requirement that the singular support of $A$ lie in $\fg\times\mathcal N^\vee\subset T_\fg^\vee$.

Basic examples are $\ICperv(\fg,\mathcal L_\chi)$ for $\chi\in\Irr W$, and these are precisely those with full support; see Lemma~\ref{lem:nilpotent.character.rs} below.
In general, nilpotent character sheaves can also have proper support.

\begin{theorem}\label{thm:cone.character}
Let $\overline C\subset\fg/\fb$ be a nonempty pure-dimensional $B$-stable closed cone, and let $C\subset\fg$ be its inverse image.
Then the following hold.
\begin{enumerate}
\item The complex $\pi_{C*}\IC_{X_C}$ is a direct sum of shifts of nilpotent character sheaves on $\fg$.
\item For every $s\in\grs$, the fiber $X_C(s)$ is pure-dimensional of dimension $\dim\overline C$, and
\[
i_s^*\IC_{X_C}\simeq\IC_{X_C(s)}
\]
where $i_s:X_C(s)\hookrightarrow X_C$ is the inclusion.
\item For every $i$, $\mathcal H^i(\pi_{C*}\IC_{X_C})|_{\grs}$ is a local system with stalk $\IH^i(X_C(s))$ at $s$.
Its pullback to $\tgrs$ is constant.
In particular, its monodromy factors through $\pi_1(\grs,s)\twoheadrightarrow W$.
\end{enumerate}
\end{theorem}

\begin{proof}
(1) The morphism $\pi_C$ is equivariant for both $G$ and scalar multiplication.
By the decomposition theorem, $\pi_{C*}\IC_{X_C}$ is a direct sum of shifts of simple $G$-equivariant monodromic perverse sheaves.
Proposition~\ref{prop:fourier.comparison} identifies $\Four_\fg(\pi_{C*}\IC_{X_C})$, up to shift and Tate twist, with $\eta_*\Four_{T_{\cB}}(\IC_{X_{\overline C}})$.
Its support is therefore contained in $\eta(T^\vee_{\cB})=\mathcal N^\vee$.
The same holds for the Fourier transform of each simple perverse constituent, proving (1).

(2) For a flag $F$, write $\fb_F$ for its Borel subalgebra.
Consider
\[
\varphi_s:G\times\cB\longrightarrow T_{\cB},
\quad
(g,F)\mapsto g\cdot(s\bmod\fb_F).
\]
Its differential at $(e,F)$ is surjective if and only if $[\fg,s]+\fb_F=\fg$.
This equality was verified in the proof of Theorem~\ref{thm:DMPS}.
By equivariance, $\varphi_s$ is smooth of relative dimension $b$.

The inverse image of $X_{\overline C}$ is
$G\times X_C(s)$.
Since $X_{\overline C}$ is pure-dimensional of dimension $\dim\cB+\dim\overline C$, it follows that $X_C(s)$ is pure-dimensional of dimension $\dim\overline C$.
Smooth pullback along $\varphi_s$, together with $\sigma_\fg^*\IC_{X_{\overline C}}\simeq\IC_{X_C}$ and $G$-equivariance, gives
\[
\IC_{G\times X_C(s)} \simeq \pr_2^*i_s^*\IC_{X_C}.
\]
Restricting to $\{e\}\times X_C(s)$ proves the asserted
isomorphism.

(3) By (1) and Lemma~\ref{lem:nilpotent.character.rs}, proved below, the restriction of $\pi_{C*}\IC_{X_C}$ to $\grs$ is a direct sum of shifts of the local systems $\mathcal L_\chi$.
Hence each cohomology sheaf is a local system whose pullback to $\tgrs$ is constant.
Proper base change and (2) identify its stalk at $s$ with
\[
\mathbb H^i(X_C(s),i_s^*\IC_{X_C}) \simeq \IH^i(X_C(s)). \qedhere
\]
\end{proof}

For the monodromy assertion, we use the following Lie algebra analogue of Lemma~\ref{lem:fullsupp.CS}. 
Over $\C$, see \cite[Lemma~6.1]{Mautner} and the related result for linear Hessenberg spaces in \cite[Proposition~2]{VX}.
For completeness, we give a proof over $\Bbbk$.


\begin{lemma}\label{lem:nilpotent.character.rs}
The nilpotent character sheaves on $\fg$ with full support are precisely $\ICperv(\fg,\mathcal L_\chi)$ for $\chi\in\Irr W$.
Every other nilpotent character sheaf restricts to zero on $\grs$.
\end{lemma}

\begin{proof}
	Let $A$ be a nilpotent character sheaf.
Write $U$ for the unipotent radical of $B$, with Lie
algebra $\fu$, and consider the natural inclusion $i:\fb\hookrightarrow\fg$ and quotients
\[
q:\fb\twoheadrightarrow\ft \and q':\fu^\perp\twoheadrightarrow\ft^\vee.
\]
Compatibility of Fourier transform with linear maps (cf.~\cite[Lemma~4.2]{Mir04}), together with base change, gives, up to Tate twist,
\[
\Four_\ft(q_!i^*A) \simeq q'_!\big((\Four_\fg A)|_{\fu^\perp}\big).
\]
Since $\mathcal N^\vee\cap\fu^\perp\subseteq\fb^\perp$ and $q'|_{\fb^\perp}=0$,
this shows that $\Four_\ft(q_!i^*A)$ is supported at zero (cf.~\cite[Lemma~4.5]{Mir04}).
Hence, $q_!i^*A$ is a constant complex on $\ft$.

On the other hand, the adjoint action map $G\times \fg\to \fg$ restricts to an isomorphism
\[
U\times\ft^{\mathrm{rs}}
\xrightarrow{~\sim~}q^{-1}(\ft^{\mathrm{rs}}),
\quad (u,h)\mapsto\Ad(u)h,
\]
under which, by the $G$-equivariance of $A$,  $(i^*A)|_{q^{-1}(\ft^{\rs})}$ corresponds to $\underline\EE_U\boxtimes A|_{\ft^{\mathrm{rs}}}$.
Since the isomorphism identifies $q$ with the projection to $\ft^{\mathrm{rs}}$ and $R\Gamma_c(U,\EE)\simeq\EE[-2\dim U](-\dim U)$, 
\[
(q_!i^*A)|_{\ft^{\mathrm{rs}}}
\simeq
A|_{\ft^{\mathrm{rs}}}[-2\dim U](-\dim U)
\]
(cf.~\cite[Lemma~4.3]{Mir04}).
This shows that  $A|_{\ft^{\mathrm{rs}}}$ is constant.
By $G$-equivariance, the pullback of $A$ along $G\times\ft^{\mathrm{rs}}\to\grs$ is also constant.
Since $G\times\ft^{\mathrm{rs}}\to G\times^T\ft^{\mathrm{rs}}$ is a torsor under the connected group $T$, the pullback of $A$ to $\tgrs\simeq G\times^T\ft^{\mathrm{rs}}$ has constant cohomology sheaves.

If $A|_{\grs}\ne0$, it is an irreducible local system shifted by $\dim\fg$, since $A$ is simple and perverse.
Its monodromy factors through $W$, so $A|_{\grs}\simeq\mathcal L_\chi[\dim\fg]$ for some $\chi\in\Irr W$.
By simplicity, $A\simeq\ICperv(\fg,\mathcal L_\chi)$.
Conversely, these intersection cohomology complexes are simple summands of the Grothendieck sheaf, whose Fourier transform is supported on $\mathcal N^\vee$.
\end{proof}

\subsection{Tangent cones of Schubert varieties and specialization}
\label{ss:nonlinear.specialization}

We now apply this construction to the tangent cone $C_{eB}\overline X_w$ of a Schubert variety at $eB$. Let $D_w\subset \fg $ denote its inverse image. If $w$ is smooth, then $D_w=H_w$ and $X_{D_w}(s)=X_{H_w}(s)$.

To extend \eqref{eq:intro.smooth.comparison} and~Theorem~\ref{thm:smooth} to this setting, we consider the following condition:
\begin{equation}\tag{SP}\label{eq:IC.specialization.condition}
\sp_{eB/\cB}\IC_{\overline X_w}\simeq\IC_{C_{eB}\overline X_w}
\quad\text{in }\mathscr D_B(\fg/\fb).
\end{equation}

\begin{proposition}\label{prop:nonlinear.specialization}
Let $w\in W$, and suppose that \eqref{eq:IC.specialization.condition} holds. Then
\beq\label{eq:nonlinear.direct.image}
\sp_e\overline K_w\simeq \pi_{D_w*}\IC_{X_{D_w}}\quad\text{in }\mathscr D_G(\fg),
\eeq
and
\beq\label{eq:nonlinear.Hess}
\Hess\big(\mathcal F_w\big)\simeq\Four_\fg\big(\pi_{D_w*}\IC_{X_{D_w}}[\dim \fg]\big)(b).
\eeq
For every $s\in\grs$ and $\mathbf s\in\Grs$, there is
an isomorphism of graded $W$-representations
\beq\label{eq:nonlinear.isom}
\IH^*(X_{D_w}(s))\simeq\IH^*(\overline Y_w(\mathbf s)).
\eeq
In particular,
$\IH^{\mathrm{odd}}(X_{D_w}(s))=0$ and
$\sum_{j\geq0}\IH^{2j}(X_{D_w}(s))u^j
=\widehat h(\widehat C_w)$.
\end{proposition}

\begin{proof}
By \eqref{eq:IC.specialization.condition}, $\sp_\Delta\IC_{\overline{\fX}_w}\simeq\IC_{G\times^BC_{eB}\overline X_w}$.
This pulls back to 
$\sigma_\fg^*\sp_\Delta\IC_{\overline{\fX}_w} \simeq\IC_{X_{D_w}}$ since $\sigma_\fg$ is smooth.
Then \eqref{eq:nonlinear.direct.image} follows from Proposition~\ref{prop:lus.hess}, and \eqref{eq:nonlinear.Hess} follows from Theorem~\ref{thm:lus.hess}.

By Theorem~\ref{thm:cone.character}(3), Lemma~\ref{lem:deform.rs}, and \eqref{eq:nonlinear.direct.image}, the cohomology sheaves of $\overline K_w|_{\Grs}$ and $(\pi_{D_w*}\IC_{X_{D_w}})|_{\grs}$ are direct sums of the local systems $L_\chi$ and $\mathcal L_\chi$, respectively, with the same multiplicities for each $\chi$ in each degree.
Together with Proposition~\ref{prop:fiber.IC}, this proves the asserted isomorphism \eqref{eq:nonlinear.isom}.
The remaining assertions follow from Corollary~\ref{cor:lus.formula}.
\end{proof}

\subsection{A sufficient local condition}\label{ss:nonlinear.conical}

To apply Proposition~\ref{prop:nonlinear.specialization}, we must verify \eqref{eq:IC.specialization.condition}. 
We give a sufficient condition in terms of the singularity of $\overline X_w$ at $eB$, and apply it to several classes of Schubert varieties.

By an \emph{affine cone}, we mean a closed subscheme of a vector space defined by a homogeneous ideal, where every coordinate has degree one.
Pointed varieties $(X,x)$ and $(Y,y)$ are \emph{smoothly equivalent} if some pointed variety $(Z,z)$ admits morphisms to both that are smooth at $z$. 

\begin{definition}\label{def:local.conical}
We consider the following local condition on $w$:
\begin{equation}\tag{LC}\label{eq:local.conical.condition}
(\overline X_w,eB)\text{ is smoothly equivalent to }(D,0)
\text{ for some affine cone }D.
\end{equation}
The smooth equivalence is not required to be equivariant.
\end{definition}
Such a cone is automatically normal and irreducible: normality of $\overline X_w$ and smooth equivalence imply normality at the vertex, and the contracting scalar action gives normality everywhere.
Every irreducible component contains the vertex, so the cone is also irreducible.

\begin{theorem}\label{thm:local.conical}
If $w$ satisfies \eqref{eq:local.conical.condition}, then $C_{eB}\overline X_w$ is normal and irreducible and \eqref{eq:IC.specialization.condition} holds.
\end{theorem}

\begin{proof}
For a cone $D$, the deformation to its tangent cone at the vertex $0$ is $M_{0/D}\cong D\times\A^1$. 
Scalar equivariance of $\IC_D$ identifies the complex on the family over $\Gm$ with the exterior product $\IC_D\boxtimes\underline\EE_{\Gm}$. 
Thus $\sp_{0/D}\IC_D\simeq\IC_D$.
This isomorphism is preserved under smooth equivalence. 
To see this, let $f:(Z,z)\to(X,x)$ be a pointed smooth morphism of relative dimension $r$. 
Then, the induced morphism of tangent cones is the projection
\[
C(f):C_zZ\simeq C_xX\times\A^r\longrightarrow C_xX.
\]
By Proposition~\ref{prop:sp.smooth} and smooth pullback for 
intersection cohomology complexes,
\[
C(f)^*\sp_{x/X}\IC_X\simeq\sp_{z/Z}\IC_Z
\and C(f)^*\IC_{C_xX}\simeq\IC_{C_zZ}.
\]
Thus an isomorphism $\sp_{x/X}\IC_X\simeq\IC_{C_xX}$ pulls back to $\sp_{z/Z}\IC_Z\simeq \IC_{C_zZ}$. 
Conversely, an isomorphism for $(Z,z)$ gives one for $(X,x)$ by restriction to a section of $C(f)$. 
Moreover, since $C_zZ\simeq C_xX\times \A^r$, $C_zZ$ is normal and irreducible if and only if $C_xX$ is. 

Applying these observations to the smooth morphisms that give the equivalence of $(\overline X_w,eB)$ and $(D,0)$, we obtain normality and irreducibility of $C_{eB}\overline X_w$, together with
\[\sp_{eB/\overline X_w}\IC_{\overline X_w}\simeq\IC_{C_{eB}\overline X_w}.\]
Compatibility of specialization with proper direct image, applied to $\overline X_w\hookrightarrow\cB$, gives the isomorphism in \eqref{eq:IC.specialization.condition} in $\mathscr D(\fg/\fb)$.
After shifting by $\ell(w)$, both sides are perverse.
Since $B$ is connected, the forgetful functor $\Perv_B(\fg/\fb)\to\Perv(\fg/\fb)$ is fully faithful (\cite[Proposition~6.2.15]{Achar}), so the isomorphism is $B$-equivariant.
\end{proof}

\begin{remark}\label{rem:local.conical.inverse}
Condition~\eqref{eq:local.conical.condition} holds for $w$ if and only if it holds for $w^{-1}$.
Indeed, identifying the opposite big cell of $\cB$ with $U^-$, the map $u\mapsto u^{-1}$ fixes the identity and identifies $U^-\cap\overline{BwB}$ with $U^-\cap\overline{Bw^{-1}B}$.
Therefore, it induces an isomorphism $(\overline X_w,eB)\simeq(\overline X_{w^{-1}},eB)$ of germs.
\end{remark}

We next give several classes of elements that satisfy the hypothesis of Theorem~\ref{thm:local.conical}.
For a parabolic subgroup $P\supset B$, let $W_P$ be the Weyl group of its Levi subgroup containing $T$, write $W^P$ for the minimal representatives of $W/W_P$, and let $w_P$ be the longest element of $W_P$.
For $v\in W^P$, write $\overline X_v^P$ for the corresponding Schubert variety in $G/P$.
We use the terms minuscule and cominuscule as in
\cite[Definitions~3.1 and~3.3]{Perrin07}.

\begin{corollary}\label{cor:local.conical.classes}
The condition~\eqref{eq:local.conical.condition} holds in each of the
following cases:
\begin{enumerate}
\item $w$ is smooth;
\item $w=vw_P$, where $P\supset B$ is minuscule or cominuscule and $v\in W^P$;
\item $G=\mathrm{GL}_n$ and $w\in S_n$ is covexillary, that is, $3412$-avoiding.
\end{enumerate}
\end{corollary}

\begin{proof}
The assertion is immediate when $\overline X_w$ is smooth.

Suppose that $G/P$ is cominuscule.
For use also in (3), we show that the local condition holds at every point (not only $eP$) of every Schubert variety in $G/P$.

The opposite big cell at $eP$ is an affine space on which a cocharacter of the center of the Levi subgroup containing $T$ acts with a single positive weight~\cite[Proposition~2.9]{GK15}.
Given a $T$-fixed point, choose $n\in N_G(T)$ taking $eP$ to that point.
Translating the cell by $n$ and conjugating the cocharacter by $n$ gives an affine chart centered at that point.
The conjugated cocharacter lies in $T$ and acts on this chart with the same single positive weight.

The intersection of a Schubert variety with such a chart centered at one of its $T$-fixed points is affine.
As a nonempty open subset of the Schubert variety, it is normal and irreducible \cite{RR85}.
Since the Schubert variety is $T$-stable, this intersection is preserved by the cocharacter action, 
hence it is an affine cone.
Every point of a Schubert variety lies in the $B$-orbit of a $T$-fixed point, so the local condition holds at every point.
Since the projection $\cB\to G/P$ is smooth and $\overline X_{vw_P}$ is the inverse image of $\overline X_v^P$, this proves (2) when $P$ is cominuscule.

If $G/P$ is minuscule but not cominuscule, it is isomorphic to a cominuscule flag variety for another group:
$C_n/P_1\simeq\PP^{2n-1}$
and
$B_n/P_n\simeq D_{n+1}/P_{n+1}$.
These isomorphisms preserve Schubert varieties \cite[Remark~3.4]{Perrin07}. The preceding argument and the smooth projection $\cB\to G/P$ therefore prove (2) in the remaining minuscule cases.

For covexillary $w$, the germ $(\overline X_w,eB)$ is smoothly equivalent to a germ of a Schubert variety in a Grassmannian \cite[Section~5 and Theorem~7.1]{AIJK21}.
Since Grassmannians are cominuscule and the preceding argument applies at every point, this proves (3).
\end{proof}

For $G=\mathrm{GL}_n$, the elements in (2) are the permutations with at most one ascent.
These permutations are covexillary, so (2) is a special case of (3).

The following example shows that \eqref{eq:local.conical.condition} can also hold for noncovexillary permutations.
\begin{example}[A noncovexillary example]
\label{ex:nonlinear.3412}
Let $G=\mathrm{GL}_4$. We identify the opposite big cell at $eB$ with the lower triangular matrices whose diagonal entries are all $1$, and write $x_{ij}$ for their entries below the diagonal.
For $w=3412$, write
\[
(a,b,c,d,e,f)=(x_{21},x_{31},x_{32},x_{41},x_{42},x_{43}).
\]
The intersection of $\overline X_w$ with this chart has defining ideal $\mathcal I=(d,acf-ae-bf)$.
Then the coordinate change $e'=e-cf$ identifies this intersection with
\[
\{ae'+bf=0\}\times\A^1_c.
\]
This is a normal irreducible affine cone, so \eqref{eq:local.conical.condition} holds.
Note that the original ideal $\mathcal I$ is not homogeneous for the standard grading.
\end{example}

Condition~\eqref{eq:IC.specialization.condition} does not hold for every $w$.
Indeed, \eqref{eq:sp.zero} and normality of $\overline X_w$ imply that \eqref{eq:IC.specialization.condition} forces $C_{eB}\overline X_w$ to be irreducible.
There exist examples with $C_{eB}\overline X_w$ reducible; see \cite
{fuchs.et.al}.

\subsection{Supports of the direct images}
\label{ss:nonlinear.domination}

For $G=\mathrm{GL}_n$, we record an extension of the full-support theorem of Precup and Sommers \cite{PS} to families of Hessenberg schemes defined by pure-dimensional $B$-stable closed cones in $\fg/\fb$, not necessarily of the form $C_{eB}\overline X_w$.

\begin{proposition}\label{prop:cone.typeA.support}
Let $G=\mathrm{GL}_n$, and let $\overline C\subset\fg/\fb$ be a nonempty pure-dimensional $B$-stable closed cone with inverse image $C$ in $\fg$.
Every simple perverse constituent of $\pi_{C*}\IC_{X_{C}}$ is of the form $\ICperv(\fg,\mathcal L_\chi)$ for some $\chi\in\Irr S_n$, and hence has support equal to $\fg$.
\end{proposition}

\begin{proof}
By Theorem~\ref{thm:cone.character}(1), it suffices to show that every nilpotent character sheaf on $\fg$ has full support.
For $G=\mathrm{GL}_n$, every simple $G$-equivariant perverse sheaf on $\mathcal N^\vee$ is a direct summand of $\Spr^\vee$.
Its inverse Fourier transform is therefore of the form $\ICperv(\fg,\mathcal L_\chi)$.
\end{proof}

\begin{remark}
	The full-support theorem can fail for other groups even when $\overline C=C_{eB}\overline X_w$ and \eqref{eq:IC.specialization.condition} holds.
See Corollary~\ref{cor:G2.full.support} for examples in type $G_2$.
\end{remark}

\section{Codimension-one Hessenberg varieties}
\label{s:applications}\label{s:codim1}

In this section, we compute $P_H(u)$ for codimension-one Hessenberg spaces in arbitrary type and determine when every graded piece is a permutation representation.
We work over $\C$ and fix a $G$-invariant nondegenerate bilinear form on $\fg$, identifying $\fg^\vee$ with $\fg$ and $\mathcal N^\vee$ with $\mathcal N$.
By Remark~\ref{rem:Hess.dot.action}, the formulas below hold over any algebraically closed field for Hessenberg spaces defined by the same sets of roots.

For an irreducible root system, $H^\perp$ is the root space corresponding to the highest root, so the calculation reduces to the Springer correspondence for the minimal nilpotent orbit.
We also determine the stalk at the origin of the intersection cohomology complex of its closure.

For $m\geq1$, write 
\[[m]=[m]_u=1+u+\cdots+u^{m-1}.\]
For a variety $Z$ with vanishing odd cohomology, write
\[p_Z(u):=\sum_i\dim \coh^{2i}(Z)u^i.\]

Let $\Phi=\Phi(G,T)$ be the root system of $G$, and let $\Phi^+$ be the set of positive roots determined by $B$. 
Let $\Phi^-=-\Phi^+$.
For a root $\alpha$, write $\mathfrak g_\alpha$ for its root space, $s_\alpha$ for its reflection, and $\alpha^\vee$ for its coroot.

Assume that $\Phi$ is irreducible.
Let $\theta$ be the highest root, $I=\mathfrak g_\theta$, $H=I^\perp$, and $P$ the stabilizer of $I$.
Then $G/P=\PP(\mathcal O_{\min})$, where $\mathcal O_{\min}$ is the minimal nonzero nilpotent orbit.
Write $\chi_{\min}$ for the Springer character of $(\mathcal O_{\min},\underline\C)$. 
Our convention assigns the trivial character to the zero orbit and differs from that of \cite[Section~13.3]{Carter} by a sign twist. 

Let $r=\operatorname{rank}\Phi$, and let $m_1\leq\cdots\leq m_r$ be the exponents of $W$ \cite[Section~3.16]{humph}.

\begin{proposition}\label{prop:codim1}
Set $k=(\dim G/P-1)/2$. Then
\beq\label{eq:codim1}
P_H(u)=[m_r]\prod_{j=1}^{r-1}[m_j+1]\,\mathbf1
       +u^kp_{P/B}(u)\chi_{\min}.
\eeq
\end{proposition}
\begin{proof}
The map $\pi_{I,P}:G\times^P I\to\overline{\mathcal O}_{\min}$, $[g,x]\mapsto\Ad(g)x$, is an isomorphism over $\mathcal O_{\min}$ and contracts the zero section $G/P$ to zero. The decomposition theorem therefore gives
\beq\label{eq:codim1.decom}
(\pi_{I,P})_*\underline\C[2k+2]
\simeq\ICperv_{\overline{\mathcal O}_{\min}}
      \oplus~\underline\C_0\otimes F,
\eeq
where $F$ is a complex of vector spaces. The first summand occurs once, without shift.

Using the notation of Remark~\ref{rem:partial.Hess}, we apply Fourier transform to \eqref{eq:codim1.decom}, restrict to $\grs$, and remove the shift $[\dim\fg+2k]$ to obtain
\[
\left.\pi_{H,P*}\underline\C\right|_{\grs}
\simeq
\mathcal L_{\chi_{\min}}[-2k]
\oplus
\underline\C_{\grs}\otimes F[-2k],
\]
by Lemmas~\ref{lem:Four.subbdl} and~\ref{lem:Four.pullback} and the Springer correspondence.
Hence
$P_{H,P}(u)=c(u)\mathbf1+u^k\chi_{\min}$
for a polynomial $c(u)$.
By \eqref{eq:flagbundle}, this gives
\[
P_H(u)=p_{P/B}(u)c(u)\mathbf1 + u^kp_{P/B}(u)\chi_{\min}.
\]

It remains to determine the coefficient of $\mathbf1$ in $P_H(u)$.
By \cite{ST,AHMMS}, this coefficient is the product of $[e+1]$ over the ideal exponents $e$ of $\Phi^+\setminus\{\theta\}$.
Since $\theta$ is the unique root of maximal height, removing it lowers the largest ideal exponent by one.
Thus these exponents are $m_1,\ldots,m_{r-1},m_r-1$, which gives the first term of \eqref{eq:codim1}.
\end{proof}

For simply laced $\Phi$, $\chi_{\min}$ is the reflection character.
In types $B_r,C_r,F_4,G_2$, it is inflated from the reflection character of the quotient by the subgroup generated by short-root reflections; these quotients are $S_r,S_2,S_3,S_2$, respectively \cite[Section~13.3]{Carter}, \cite[Section~6.2]{Juteau}.

\begin{corollary}\label{cor:codim1:permutation}
Let $\Phi$ be irreducible and $H\subset\fg$ be a Hessenberg space of codimension one.
Every graded piece of $\coh^*(X_H(s))$ is a permutation representation if and only if $\Phi$ has type $A_r$, $B_r$, $C_r$, $F_4$, or $G_2$. 
In types $D_r$, $E_6$, $E_7$, and $E_8$, neither $\coh^{2k}(X_H(s))$ nor the total cohomology $\coh^*(X_H(s))$ is a permutation representation.
\end{corollary}
\begin{proof}
Write \eqref{eq:codim1} as $P_H=a\mathbf1+b\chi_{\min}$, where $b=u^kp_{P/B}$.
In the five types listed, $\mathbf1+\chi_{\min}$ is a permutation character.
Indeed, it is the character of the natural action on $r+1$ letters in type $A_r$.
In the other types, it is inflated from the permutation character of the corresponding symmetric-group quotient, listed above, of $W$, given by the sum of its trivial and reflection characters.
Thus
\[
P_H=(a-b)\mathbf1+b(\mathbf1+\chi_{\min}).
\]
Since $p_{P/B}$ has nonnegative coefficients, it suffices to show that $(a-b)/p_{P/B}$ is a polynomial with nonnegative coefficients.
These quotients are given in the following table.
\[
\begin{array}{c|l}
\Phi&(a-b)/p_{P/B}\\\hline
A_r&[r]^2-u^{r-1}\\
B_r&[r-1]_{u^2}[2r-1]-u^{2r-3}\\
C_r&[2r-1]-u^{r-1}\\
F_4&(1+u^4)[11]-u^7\\
G_2&[5]-u^2
\end{array}
\]
Each polynomial in the table has nonnegative coefficients.

For types $D_r$ ($r\geq4$), $E_6,E_7,E_8$, choose a standard parabolic subgroup $W_4$ of type $D_4$, and let $\chi_4$ be its reflection character. 
Since $W_4$ acts trivially on the orthogonal complement of the span of its four simple roots, we have
\[
\chi_{\min}|_{W_4}=\chi_4+(r-4)\mathbf1.
\]
We claim that $c\mathbf1+d\chi_{\min}$ cannot be a permutation character when $c\geq0$ and $d>0$.
Otherwise, its restriction to $W_4$ would be a permutation character with no irreducible constituents other than
$\mathbf1$ and $\chi_4$, and each transitive orbit would have character $\mathbf1+m\chi_4$ for some $m\geq0$.

Since $d>0$, some orbit would have $m>0$.
Let $z\in W_4$ be the element acting as $-1$ on its reflection representation. Then $\chi_4(z)=-4$, so this orbit would have character value $1-4m<0$ at $z$, a contradiction.
Since \eqref{eq:codim1} gives
\[
[\coh^{2k}(X_H(s))]=a_k\mathbf1+\chi_{\min}
\and
P_H(1)=a(1)\mathbf1+p_{P/B}(1)\chi_{\min},
\]
where $a_k$ is the coefficient of $u^k$ in $a(u)$, and since $p_{P/B}(1)>0$, the claim applies to both.
\end{proof}

The decomposition \eqref{eq:codim1.decom} also determines the stalk of $\IC_{\overline{\mathcal O}_{\min}}$ at the origin.
We continue to write $P$ for the stabilizer of $I=\fg_\theta$, and write $b_j=\dim\coh^{2j}(G/P)$ with $b_{-1}=0$.

\begin{corollary}\label{cor:min.orbit.stalk}
The stalk of $\IC_{\overline{\mathcal O}_{\min}}$ at the origin has cohomology only in even degrees, and
\begin{equation}\label{eq:min.orbit}
\sum_j
\dim\cH^{2j}
\bigl(\IC_{\overline{\mathcal O}_{\min}}\bigr)_0\,u^j
=\sum_{j=0}^{k}(b_j-b_{j-1})u^j.
\end{equation}
Moreover, $\dim\chi_{\min}=b_k$. 
In particular, $\overline{\mathcal O}_{\min}$ is rationally smooth if and only if $\dim\chi_{\min}=1$.
\end{corollary}

\begin{proof}
By proper base change, the stalk of the left-hand side of \eqref{eq:codim1.decom} at the origin is $R\Gamma(G/P,\underline\C)[2k+2]$.
Since $\coh^{\mathrm{odd}}(G/P)=0$, both $F$ and
$(\ICperv_{\overline{\mathcal O}_{\min}})_0$
have cohomology only in even degrees.
Let $f(u)$ and $g(u)$ be their respective Poincar\'e polynomials, with degree $2j$ recorded by $u^j$. They are Laurent polynomials with nonnegative integer coefficients, and
\[
u^{-k-1}p_{G/P}(u)=f(u)+g(u).
\]
Since $G\times^P I$ is smooth of dimension $2k+2$ and $\pi_{I,P}$ is proper, $(\pi_{I,P})_*\underline\C[2k+2]$ is Verdier self-dual \cite[Corollary~2.8.8]{Achar}.
The same holds for $\ICperv_{\overline{\mathcal O}_{\min}}$ \cite[Lemma~3.3.13]{Achar}.
Thus \eqref{eq:codim1.decom} implies that $F$ is self-dual, and hence $f(u)=f(u^{-1})$.
Since the stalk of $\ICperv_{\overline{\mathcal O}_{\min}}$ at the origin has no cohomology in nonnegative degrees \cite[Lemma~3.3.11]{Achar}, 
$g(u)$ has only negative powers.
Consequently, the nonnegative powers in $u^{-k-1}p_{G/P}(u)$ determine $f(u)$.

For $0\leq j\leq k$, the coefficient of $u^{j-k-1}$ in $g(u)$ is therefore $b_j-b_{2k+2-j}=b_j-b_{j-1}$, by Poincar\'e duality on $G/P$, where $b_{2k+2}=0$.
Since $\ICperv_{\overline{\mathcal O}_{\min}}=\IC_{\overline{\mathcal O}_{\min}}[2k+2]$, this proves
\eqref{eq:min.orbit}.

As in the proof of Proposition~\ref{prop:codim1}, Fourier transform gives
$P_{H,P}(u)=u^k(f(u)\mathbf1+\chi_{\min})$.
Since $X_{H,P}(s)$ and $G/P$ have the same Euler characteristic,
evaluating at $u=1$ 
gives
\[
f(1)+\dim\chi_{\min}=p_{G/P}(1)=f(1)+g(1).
\]
Hence $\dim\chi_{\min}=g(1)=\sum_{j=0}^{k}(b_j-b_{j-1})=b_k$.
Since \eqref{eq:min.orbit} has nonnegative coefficients, constant term one, and value $b_k=\dim\chi_{\min}$ at $u=1$, it equals $1$ if and only if $\dim\chi_{\min}=1$.
Since $\overline{\mathcal O}_{\min}$ is smooth away from the origin, this proves the last assertion.
\end{proof}

\section{Explicit computations}\label{s:characters}

We first compute $P_H(u)$ for all Hessenberg spaces $H$ in types $G_2$ and $C_3$.
The character formula determines $P_{H_w}(u)$ for every smooth element $w\in W$.
We obtain the remaining characters from the modular law \eqref{eq:modular.law}, together with Proposition~\ref{prop:codim1} and a calculation involving
nilpotent orbits in type $C_3$.
In type $G_2$, we also determine $\widehat h(\widehat C_w)$ and the nonprincipal part of $\Hess(\mathcal F_w)$ for every $w\in W$.

We then explain in Subsection~\ref{ss:higher.rank} how to combine the character formula and modular law with Precup--Sommers' decompositions of Hessenberg ideal fibers to compute $P_H(u)$ in higher rank.
We have verified computationally that these methods determine $P_H(u)$ for every Hessenberg space in types $B_n$, $C_n$, and $D_n$ for $n\leq10$, and in all exceptional types.

We work over $\C$ and write $[m]=1+u+\cdots+u^{m-1}$.
As in Section~\ref{s:codim1}, our convention for the Springer correspondence assigns the trivial character to the zero orbit.

\subsection{Character formula and the modular law}
Recall from Corollaries~\ref{cor:hess.formula} and~\ref{cor:lus.formula} that
\beq\label{eq:char.formula}
\widehat h(\widehat C_w)
=\sum_{\chi,\psi\in\Irr W}\{\chi,\psi\} \widetilde\psi(\widehat C_w)\chi
=\sum_i\IH^{2i}(\overline Y_w(\mathbf s))u^i,
\eeq
where $\{ , \}$ is Lusztig's Fourier pairing.
If $w$ is smooth, this equals $P_{H_w}(u)$ by Theorem~\ref{thm:smooth}.

Since $\widehat h$ is a trace, it factors through the cocenter $\hecke/[\hecke,\hecke]$.
For each conjugacy class $C\subset W$, choose an element $w_C\in C$ of minimal length.
The classes of $T_{w_C}$ form an $\cA$-basis of the cocenter \cite[Theorem~8.2.3]{GP}, which we use to evaluate the Hecke characters in \eqref{eq:char.formula}.

We also use the modular law of Precup--Sommers \cite{PS}.
For a simple root $\alpha$, let $P_\alpha\supset B$ be the corresponding minimal parabolic subgroup.
A triple $(I_0,I_1,I_2)$ of $B$-stable ideals in $\mathfrak u$ is called a \emph{modular triple} if, for some simple root $\alpha$ and positive root $\beta$,
\[
I_0=I_1\oplus\mathfrak g_\beta,\quad
I_1=I_2\oplus\mathfrak g_{\alpha+\beta},\quad
\langle\beta,\alpha^\vee\rangle=-1,
\]
and $I_0$ is $P_\alpha$-stable.
For such a triple, \cite[Proposition~6.5]{PS} gives
\beq\label{eq:modular.law}
(1+u)P_{I_1^\perp}(u)=P_{I_2^\perp}(u)+uP_{I_0^\perp}(u).
\eeq

\subsection{Type \texorpdfstring{$G_2$}{G2}}\label{ss:G2.characters}
Let $G$ be the simple group of type $G_2$.
Let $\alpha$ and $\beta$ be the short and long simple roots, with reflections $s$ and $t$, respectively, and write $w_0=(st)^3$.
We write
\[
\Irr W=\{\mathbf1,\epsilon,\chi_1,\chi_2,\epsilon_1,\epsilon_2\},
\]
where $\epsilon$ is the sign character, $\chi_1$ is the reflection character, and
\[
(\epsilon_1(s),\epsilon_1(t))=(1,-1),\quad
(\epsilon_2(s),\epsilon_2(t))=(-1,1),\quad
\chi_2=\chi_1\otimes\epsilon_1.
\]

There are two singleton families, $\{\mathbf1\}$ and $\{\epsilon\}$, and one family $\mathscr F=\{\chi_1,\chi_2,\epsilon_1,\epsilon_2\}$.
Lusztig's embedding $\mathscr F\hookrightarrow\mathscr M(S_3)$ sends these four characters to $(1,1),(g_2,1),(1,r),(g_3,1)$, respectively, where $g_i$ is an $i$-cycle and $r$ is the two-dimensional irreducible character of $S_3$ \cite[p.~95]{Lus84}.
In this ordering, the restriction of the Fourier pairing to $\mathscr F\times\mathscr F$ has matrix \cite[p.~110]{Lus84}
\[
\frac16
\begin{pmatrix}
1&3&2&2\\
3&3&0&0\\
2&0&4&-2\\
2&0&-2&4
\end{pmatrix}.
\]
On each singleton family, the Fourier pairing has matrix $(1)$.

Every Schubert variety in type $G_2$ is rationally smooth, so $\widehat C_w=\sum_{v\leq w}T_v$.
We reduce these sums to the cocenter basis represented by $T_e,T_s,T_t,T_{st},T_{(st)^2},T_{w_0}$ (where $w_0$ denotes the longest element $(st)^3=(ts)^3$) and use the Hecke character table in \cite[p.~270]{GP}.
For example, $\widehat C_{st}=(1+T_s)(1+T_t)$ has trivial and sign character values $[2]^2$ and $0$, while
\[
(\widetilde\chi_1,\widetilde\chi_2,
 \widetilde\epsilon_1,\widetilde\epsilon_2)(\widehat C_{st})
=(3u,u,0,0).
\]
Substituting these values into \eqref{eq:char.formula}, we obtain
\[
\widehat h(\widehat C_{st})
=[2]^2\mathbf1+u(\chi_1+2\chi_2+\epsilon_1+\epsilon_2).
\]
The same calculation gives the remaining values in Table~\ref{tab:G2}.

\begin{table}[htbp]
\centering
\renewcommand{\arraystretch}{1.18}
\begin{tabular}{@{}clc@{}}
\toprule
$w$ & $\widehat h(\widehat C_w)$ & $\codim H_w$\\
\midrule
$e$ & $\mathbf1+\epsilon+\epsilon_1+\epsilon_2+2\chi_1+2\chi_2$ & $6$\\
$s$ & $[2](\mathbf1+\epsilon_1+\chi_1+\chi_2)$ & $5$\\
$t$ & $[2](\mathbf1+\epsilon_2+\chi_1+\chi_2)$ & $5$\\
$st,ts$ & $[2]^2\mathbf1+u(\chi_1+2\chi_2+\epsilon_1+\epsilon_2)$ & $4$\\
$sts$ & $[2][3]\mathbf1+u[2](\chi_2+\epsilon_2)$ & ---\\
$tst$ & $[2][3]\mathbf1+u[2](\chi_2+\epsilon_1)$ & $3$\\
$stst,tsts$ & $[2][4]\mathbf1+u^2(\chi_2+\epsilon_1+\epsilon_2)$ & ---\\
$ststs$ & $[2][5]\mathbf1+u^2[2]\epsilon_1$ & ---\\
$tstst$ & $[2][5]\mathbf1+u^2[2]\epsilon_2$ & ---\\
$w_0$ & $[2][6]\mathbf1$ & $0$\\
\bottomrule
\end{tabular}
\par
\vspace{8pt}
\caption{$\widehat h(\widehat C_w)$ in type $G_2$.
A dash in the last column indicates that $w$ is not smooth.}
\label{tab:G2}
\end{table}

The smooth elements in Table~\ref{tab:G2} determine $P_H(u)$ for six of the eight Hessenberg spaces.
It remains to consider the spaces of codimensions one and two.
Let $H_i$ be the Hessenberg space of codimension $i$, for $i=0,1,2$, and let $I_i=H_i^\perp$. 
Then $I_0=0$, $I_1=\fg_{3\alpha+2\beta}$, and $I_2=\fg_{3\alpha+\beta}\oplus\mathfrak g_{3\alpha+2\beta}$.
The stabilizer $P$ of $I_1$ satisfies $P/B\simeq\PP^1$.

The trivial local system on the minimal nonzero nilpotent orbit corresponds to $\epsilon_1$ under the Springer correspondence \cite{Carter}.
Hence Proposition~\ref{prop:codim1} gives
\beq\label{eq:G2.codim1}
P_{H_1}(u)=[2][5]\mathbf1+u^2[2]\epsilon_1.
\eeq
Since $I_2$ is $P_\beta$-stable and $s_\beta(3\alpha+\beta)=3\alpha+2\beta$, the triple $(I_2,I_1,I_0)$ forms a modular triple. 
Hence
$[2]P_{H_1}(u)=P_{H_0}(u)+uP_{H_2}(u)$.
Using $P_{H_0}(u)=[2][6]\mathbf1$, we obtain
\[
P_{H_2}(u)=[2][4]\mathbf1+u[2]^2\epsilon_1.
\]
Together with Table~\ref{tab:G2}, these formulas recover all eight characters in \cite[Table~5]{xue}.

\medskip

We also determine the nonprincipal part of $\Hess(\mathcal F_w)$ for all $w\in W$  in type $G_2$.
The result in particular shows that the full-support theorem of Precup and Sommers does not extend to families of Hessenberg schemes defined by $B$-stable closed cones in $\fg/\fb$ in general.

\begin{proposition}[Nonprincipal summands in type $G_2$]
\label{prop:G2.nonprincipal}
Let $A_{\mathrm{cusp}}$ be the unipotent character sheaf labelled by $(1,\mathrm{sgn}_{S_3})$ in Lusztig's parametrization, and set $\mathcal Q=\mu_e A_{\mathrm{cusp}}$.
Then $\mathcal Q$ is a nonzero simple perverse sheaf on $\mathcal N$ with $\mathsf S(\mathcal Q)=0$, and there are noncanonical isomorphisms
\[
\Hess_{\mathrm{np}}(\mathcal F_w)\simeq
\begin{cases}
\mathcal Q[-4],
    & w=stst,\ tsts,\\
\mathcal Q[-4]\oplus\mathcal Q[-6],
    & w=ststs,\ tstst,\\
0,  & \text{otherwise}.
\end{cases}
\]
\end{proposition}

\begin{proof}
By \cite[Appendix~C]{DLM}, $A_{\mathrm{cusp}}$ is the perverse intersection cohomology complex of the subregular unipotent orbit with the local system corresponding to the sign representation of $S_3$.
Let $\mathcal O_{\mathrm{sub}}$ be the subregular nilpotent orbit, with the corresponding local system $\mathcal E_{\mathrm{sgn}}$. 

We claim that there is an isomorphism
\beq\label{eq:G2.cuspidal.specialization}
\sp_e A_{\mathrm{cusp}}
\simeq \ICperv(\overline{\mathcal O}_{\mathrm{sub}},\mathcal E_{\mathrm{sgn}})=:K.
\eeq
This implies that $\mathcal Q\simeq\Four_\fg(K)$ is nonzero, simple, and perverse.
Since $A_{\mathrm{cusp}}$ is a unipotent character sheaf and has  proper support, $\mathcal Q=\mu_eA_{\mathrm{cusp}}$
is supported on $\mathcal N$ and $\mathsf S(\mathcal Q)=0$.

\smallskip

Since nilpotent orbits are stable under scalar multiplication, $\overline{\mathcal O}_{\mathrm{sub}}$ is an affine cone.
Its deformation to the tangent cone at zero is therefore isomorphic to $\overline{\mathcal O}_{\mathrm{sub}}\times\A^1$.
To compute specialization using this product structure, we will also establish scalar equivariance of $K$.
This is an assertion about $\mathcal E_{\mathrm{sgn}}$, which does not follow from scalar stability of its support.

Fix $x\in\mathcal O_{\mathrm{sub}}$.
There exists a cocharacter $\lambda:\Gm\to G$ such that $\Ad(\lambda(t))x=tx$.
(The subregular orbit is even \cite[Example~8.2.13]{CMcG} and the derived subgroup of $G$ is adjoint, so a Dynkin cocharacter for $x$ is divisible by two.)
Let $G\times\Gm$ act on $\mathcal O_{\mathrm{sub}}$ by $(g,t)\cdot y=t\Ad(g)y$.
The identity above shows that the map
\[
t\longmapsto(\lambda(t)^{-1},t)
\]
is a section of the projection from the stabilizer of $x$ to $\Gm$.
Thus the stabilizer is isomorphic to $Z_G(x)\rtimes\Gm$.
Since $\Gm$ is connected, the inclusion of $Z_G(x)$
induces an isomorphism on component groups.
By the classification of equivariant local systems on homogeneous
spaces \cite[Proposition~6.2.13]{Achar}, the $G$-equivariant
structure on $\mathcal E_{\mathrm{sgn}}$ therefore extends to
a $G\times\Gm$-equivariant structure.
In particular, $K$ is equivariant for scalar multiplication.

Under the identification $M_{0/\overline{\mathcal O}_{\mathrm{sub}}}\simeq \overline{\mathcal O}_{\mathrm{sub}}\times\A^1$ above,
the canonical open immersion
\[
j:\overline{\mathcal O}_{\mathrm{sub}}\times\Gm
\hookrightarrow M_{0/\overline{\mathcal O}_{\mathrm{sub}}}
\]
is given by $j(x,t)=(t^{-1}x,t)$.
Scalar equivariance of $K$ therefore identifies the complex on the punctured deformation with $K\boxtimes\underline\EE_{\Gm}$.
Taking nearby cycles gives
$\sp_{0/\overline{\mathcal O}_{\mathrm{sub}}}K\simeq K$.

The exponential map restricts to a $G$-equivariant algebraic isomorphism from $\overline{\mathcal O}_{\mathrm{sub}}$ onto the closure of the subregular unipotent class, identifying $K$ with $A_{\mathrm{cusp}}$.
Since $\overline{\mathcal O}_{\mathrm{sub}}$ is a cone and the differential of the exponential map at zero is the identity, the tangent cone of the latter closure at $e$ is identified with $\overline{\mathcal O}_{\mathrm{sub}}\subset\fg$.
By compatibility of specialization with this isomorphism and with closed immersions, the isomorphism $\sp_{0/\overline{\mathcal O}_{\mathrm{sub}}}K\simeq K$ gives \eqref{eq:G2.cuspidal.specialization}.

\smallskip

We now determine the multiplicities of shifts of $\mathcal Q$ in $\Hess_{\mathrm{np}}(\mathcal F_w)$.
Lusztig's Fourier pairings of $A_{\mathrm{cusp}}$ with $\chi_1,\chi_2,\epsilon_1,\epsilon_2$ are $1/6,-1/2,1/3,1/3$, respectively
\cite[pp.~95, 110]{Lus84}.
Set
\[
\tau=\frac16\widetilde\chi_1-\frac12\widetilde\chi_2
     +\frac13\widetilde\epsilon_1
     +\frac13\widetilde\epsilon_2.
\]
By \eqref{e 8}, $\tau(\widehat C_w)$ records the graded multiplicity of $A_{\mathrm{cusp}}$ in $\overline K_w[\dim G]$, with perverse degree $i$ recorded by $u^{i/2}$.
Using $\widehat C_w=\sum_{v\leq w}T_v$, the Hecke character table \cite[p.~270]{GP} gives
\[
\tau(\widehat C_w)=
\begin{cases}
u^2,&w=stst,\ tsts,\\
u^2+u^3,&w=ststs,\ tstst,\\
0,&\text{otherwise}.
\end{cases}
\]
By \eqref{e 8} and Theorem~\ref{thm:lus.hess}, the values $u^2$ and $u^2+u^3$ give the summands $\mathcal Q[-4]$ and $\mathcal Q[-4]\oplus\mathcal Q[-6]$, respectively.

Finally, every unipotent character sheaf with proper support other than $A_{\mathrm{cusp}}$ is supported on the closure of a conjugacy class with nonidentity semisimple part \cite[Appendix~C]{DLM}.
Such a closure does not contain the identity, so the corresponding sheaf specializes to zero at the identity.
Since the full-support unipotent character sheaves microlocalize to simple summands of $\Spr^\vee$, as in the proof of Corollary~\ref{cor:hess.formula}, all nonprincipal summands come from $A_{\mathrm{cusp}}$.
\end{proof}

\begin{corollary}[Failure of full support]\label{cor:G2.full.support}
For $w\in\{stst,tsts,tstst\}$, condition~\eqref{eq:IC.specialization.condition} holds, but $\pi_{D_w*}\IC_{X_{D_w}}$ has a nonzero direct summand with proper support.
\end{corollary}
\begin{proof}
Let $P_t\supset B$ be the parabolic subgroup with Weyl group $\langle t\rangle$.
The projection
\[p:G/B\lra G/P_t\simeq Q^5\]
is a $\PP^1$-bundle \cite[Proposition~4.1]{anderson11}, where $Q^5\subset\PP(V)$ is the quadric of isotropic lines in the seven-dimensional representation $V$.
For the fixed $B$-stable isotropic flag $F_\bullet$,
\[
\overline X_{stst}
=p^{-1}\big(Q^5\cap\PP(F_2^\perp)\big)
\and
\overline X_{tstst}
=p^{-1}\big(Q^5\cap\PP(F_1^\perp)\big)
\]
\cite[Section~4.2]{anderson11}.
The quadratic form on $F_2^\perp$ has rank three and radical $F_2$, while $\PP(F_1^\perp)$ is the tangent hyperplane to $Q^5$ at $[F_1]$.
Since $p$ is a $\PP^1$-bundle, we obtain isomorphisms of germs
\[
\begin{aligned}
(\overline X_{stst},eB)
&\simeq
\bigl(\A^2\times\{xy-z^2=0\},0\bigr),\\
(\overline X_{tstst},eB)
&\simeq
\bigl(\A^1\times
\{x_1x_2+x_3x_4+x_5^2=0\},0\bigr).
\end{aligned}
\]
Both are normal irreducible affine cones.
The same is true for $w=tsts$ by Remark~\ref{rem:local.conical.inverse}.
Thus \eqref{eq:local.conical.condition} holds in all three cases, and Theorem~\ref{thm:local.conical} gives \eqref{eq:IC.specialization.condition}.

For each $w\in\{stst,tsts,tstst\}$, $\mathcal Q[-4]$ is a direct summand of $\Hess(\mathcal F_w)$ by Proposition~\ref{prop:G2.nonprincipal}.
Hence, by \eqref{eq:nonlinear.Hess} and \eqref{eq:G2.cuspidal.specialization}, $\ICperv(\overline{\mathcal O}_{\mathrm{sub}},
        \mathcal E_{\mathrm{sgn}})[-\dim\fg-4]$ is a direct summand of $\pi_{D_w*}\IC_{X_{D_w}}$.
\end{proof}

\subsection{Type \texorpdfstring{$C_3$}{C3}}\label{ss:C3}

Let $G=\operatorname{Sp}_6(\C)$, and let $e_1,e_2,e_3$ be the standard basis of $\mathbb R^3$. 
We choose the simple roots $\alpha_1=e_1-e_2$, $\alpha_2=e_2-e_3$, and $\alpha_3=2e_3$.
Let $s_i=s_{\alpha_i}$ for $i=1,2,3$, and set
\[
t=t_0:=s_3,\quad t_1:=s_2ts_2,\and t_2:=s_1t_1s_1.
\]
The generators denoted by $t,s_1,s_2$ in \cite[Example~3.4.8 and Table~10.3]{GP} correspond to $t,s_2,s_1$, respectively, in our notation.

We label the irreducible characters of $W$ by bipartitions $(\lambda,\mu)$ of $3$ and set
\[
\mathbf1=\chi_{(3,\emptyset)},\quad
\rho=\chi_{(2,1)},\quad \sigma=\chi_{(21,\emptyset)},\and
\epsilon_t=\chi_{(\emptyset,3)}.
\]
Here $\rho$ is the three-dimensional reflection character, $\sigma$ has dimension two, and $\epsilon_t$ is the one-dimensional character with $\epsilon_t(s_1)=\epsilon_t(s_2)=1$ and $\epsilon_t(t)=-1$.
The sign character is $\chi_{(\emptyset,111)}$.

There are two nonsingleton families,
\[
(\rho,\sigma,\epsilon_t) \and
(\chi_{(1,11)},\chi_{(\emptyset,21)},\chi_{(111,\emptyset)}),
\]
and four singleton families. With the orderings above, the restriction
of the Fourier pairing to either nonsingleton family has matrix
\cite[Section~4.5]{Lus84}
\[
\frac12
\begin{pmatrix}1&1&1\\1&1&-1\\1&-1&1\end{pmatrix}.
\]
\begin{example}\label{ex:C3.character}
We compute $P_{H_w}(u)$ for the smooth element $w=t_0t_1s_1$ using \eqref{eq:char.formula}.
The Bruhat interval below $w$ gives the factorization $\widehat C_w =(1+T_t)(1+T_{s_2}+T_{s_2t}+T_{t_1})(1+T_{s_1})$.
The quadratic and braid relations and cyclically permuting factors give, modulo $[\hecke,\hecke]$,
\[
\begin{aligned}
\widehat C_w\equiv{}&
T_e+(u+1)T_{t_0}+T_{t_0t_1}+(u+2)T_{s_2}
 +3uT_{t_0s_2}\\
&\quad+T_{t_0s_1}+T_{t_0t_1s_1}
 +(u+1)T_{s_2s_1}+2uT_{t_0s_2s_1}.
\end{aligned}
\]
The elements indexing the right-hand side are minimal-length representatives of distinct conjugacy classes
\cite[Example~3.4.8]{GP}.
The Hecke character table \cite[Table~10.3]{GP} now gives
\[
\widetilde{\mathbf1}(\widehat C_w)=[2]^2[4],\quad
\widetilde\rho(\widehat C_w)=u[2]^3,\and 
\widetilde\sigma(\widehat C_w)=u[2](1+u^2),
\]
while the other irreducible character values vanish.
Applying \eqref{eq:char.formula}, we obtain
\beq\label{eq:C3.ex}
P_{H_w}(u)=\widehat h(\widehat C_w)
=[2]^2[4]\mathbf1+u[2][3](\rho+\sigma)+u^2[2]\epsilon_t.
\eeq

Formula~\eqref{eq:C3.ex} can also be obtained geometrically.
See Remark~\ref{rem:C3.blowup} below. 
\end{example}

\smallskip
We now compute $P_H(u)$ for all twenty Hessenberg spaces of type $C_3$.
We number the corresponding ideals as in \cite[Section~8]{PS}, and set $H_i=I_i^\perp$ and $P_i=P_{H_i}$.
Table~\ref{tab:C3.ideals} lists the minimal elements of $\Phi_{I_i}:=\{\alpha\in\Phi^+:\mathfrak g_\alpha\subset I_i\}$.

\begin{table}[H]
\centering
\small
\renewcommand{\arraystretch}{1.2}
\begin{tabular}{@{}clc@{\qquad}clc@{}}
\toprule
$i$ & $\min\Phi_{I_i}$ & $\dim H_i/\mathfrak b$ &
$i$ & $\min\Phi_{I_i}$ & $\dim H_i/\mathfrak b$\\
\midrule
$20$ & $\emptyset$ & $9$ & $11$ & $001$ & $3$\\
$19$ & $221$ & $8$ & $9$ & $100,021$ & $3$\\
$18$ & $121$ & $7$ & $8$ & $110,011$ & $3$\\
$17$ & $021$ & $6$ & $7$ & $010$ & $2$\\
$16$ & $111$ & $6$ & $6$ & $100,011$ & $2$\\
$15$ & $110$ & $5$ & $5$ & $110,001$ & $2$\\
$13$ & $111,021$ & $5$ & $4$ & $100,001$ & $1$\\
$14$ & $100$ & $4$ & $3$ & $010,001$ & $1$\\
$12$ & $011$ & $4$ & $2$ & $100,010$ & $1$\\
$10$ & $110,021$ & $4$ & $1$ & $100,010,001$ & $0$\\
\bottomrule
\end{tabular}
\par
\vspace{8pt}
\caption{Hessenberg ideals in type $C_3$.
Here $abc$ denotes $a\alpha_1+b\alpha_2+c\alpha_3$.}
\label{tab:C3.ideals}
\end{table}

For $i\neq13,17$, there is a smooth element $w$ with $H_i=H_w$ \cite[Example~5.6]{BHL25}, so $P_i(u)$ is given by \eqref{eq:char.formula}.
To compute $P_{13}$ and $P_{17}$, consider the ideals
\[
\begin{split}
I_{17}&=\mathfrak g_{2\alpha_2+\alpha_3}
 \oplus\mathfrak g_{\alpha_1+2\alpha_2+\alpha_3}
 \oplus\mathfrak g_{2\alpha_1+2\alpha_2+\alpha_3},\\
I_{13}&=I_{17}\oplus\mathfrak g_{\alpha_1+\alpha_2+\alpha_3},\;\and
I_{10}=I_{13}\oplus\mathfrak g_{\alpha_1+\alpha_2}.
\end{split}
\]
The ideal $I_{10}$ is $P_{\alpha_3}$-stable and
$s_{\alpha_3}(\alpha_1+\alpha_2)=\alpha_1+\alpha_2+\alpha_3$.
Thus the triple $(I_{10},I_{13},I_{17})$ is a modular triple. The modular law gives
\beq\label{eq:C3.modular.law}
[2]P_{13}(u)=P_{17}(u)+uP_{10}(u).
\eeq
Since $P_{10}$ is known, it remains to compute $P_{17}$.

\begin{proposition}\label{prop:C3.codim3}
We have 
\[P_{17}(u)=[2][3][4]\mathbf1+u[2]^2[3]\epsilon_t+u^2[2]^2\rho.\]
\end{proposition}

\begin{proof}
Let $V=\C^6$ with its standard symplectic form.
Under the identification $\mathfrak{sp}(V)\simeq\Sym^2V$, we have $I_{17}=\Sym^2L_0$, where $L_0\subset V$ is the $B$-stable isotropic two-plane. Let $Q$ be the stabilizer of $L_0$.
Then $I_{17}$ is $Q$-stable. Write $P_{H_{17},Q}(u)$ for the graded character of
\[X_{H_{17},Q}(s)=\{gQ:\Ad(g^{-1})s\in H_{17}\}.\]
Since $Q/B\simeq\PP^1\times\PP^1$, the same Leray--Hirsch argument as in \eqref{eq:flagbundle} gives
\[
P_{17}(u)=[2]^2P_{H_{17},Q}(u) \and \dim X_{H_{17},Q}(s)=4.
\]
Consider the proper map
\[
\phi=\pi_{I_{17},Q}:G\times^Q I_{17}\lra \overline{\mathcal O}_{[2^2,1^2]},
\quad[g,x]\mapsto\Ad(g)x,
\]
where the subscript records the sizes of the Jordan blocks in the nilpotent orbit.
Its fiber over $x$ consists of isotropic two-planes $L$ satisfying $\image x\subset L$ and $x(L^\perp)=0$.
If $x$ has rank two, necessarily $L=\image x$.
The inverse $x\mapsto\image x$ is regular on this orbit, so $\phi$ is an isomorphism over $\mathcal O_{[2^2,1^2]}$.
It follows that $\ICperv_{\overline{\mathcal O}_{[2^2,1^2]}}$ occurs in $\phi_*\underline\C[10]$ with multiplicity one and no shift.
No summand corresponds to a nontrivial local system on this orbit.

The boundary consists of $\mathcal O_{[2,1^4]}$ and $\{0\}$.
Over $x\in\mathcal O_{[2,1^4]}$, the fiber is $\PP((\image x)^\perp/\image x) \simeq\PP^3$. 
As $x$ varies over this orbit, the lines $\image x$ form a subbundle of the constant bundle with fiber $V$.
The quotients $(\image x)^\perp/\image x$ therefore form a rank-four vector bundle, and the restriction of $\phi$ is its projectivization.
Powers of the relative hyperplane class trivialize the cohomology local systems. The decomposition theorem then excludes summands associated with nontrivial local systems on this orbit.
Under the Springer correspondence, the trivial local systems on $\mathcal O_{[2^2,1^2]}$, $\mathcal O_{[2,1^4]}$, and $\{0\}$ correspond to $\rho$, $\epsilon_t$, and $\mathbf1$, respectively \cite[pp.~424--425]{Carter}.
Applying Fourier transform to $\phi_*\underline\C[10]$, as in the proof of Proposition~\ref{prop:codim1}, gives
\beq\label{eq:C3.supp.red}
P_{H_{17},Q}(u)=a(u)\mathbf1+b(u)\epsilon_t+u^2\rho.
\eeq
Here the multiplicity of $\rho$ is $u^2$: since $\dim(G\times^QH_{17})=25=\dim\mathfrak g+4$, the perverse summand $\ICperv_{\overline{\mathcal O}_{[2^2,1^2]}}$ corresponds to one copy of $\rho$ in $\coh^4(X_{H_{17},Q}(s))$.

The multiplicity of the trivial representation in $P_{17}(u)$ is given by the product of $[e+1]$ for the ideal exponents $e$ \cite{ST,AHMMS}.
The positive roots $\alpha$ for which $\mathfrak g_{-\alpha}\subset H_{17}$ form the set
\[
\Phi^+_{H_{17}}
=\Phi^+\setminus\{2e_1,e_1+e_2,2e_2\}.
\]
It has three roots of height one, two of height two, and one of height three.
Its ideal exponents are $1,2,3$, so the trivial multiplicity is $[2][3][4]$.
It follows that
$a(u)=[2][3][4]/[2]^2=[3](1+u^2)$.
The dimensions of the attracting cells give the Poincar\'e polynomial of $X_{H_{17}}(s)$ as
\[
\sum_{v\in W}u^{|v(\Phi^+_{H_{17}})\cap\Phi^-|}
=[2]^2(1+2u+6u^2+2u^3+u^4).
\]
Dividing by $[2]^2$ and using $\dim\rho=3$ and $\dim\epsilon_t=1$ in
\eqref{eq:C3.supp.red} yields $b(u)=u[3]$.
Multiplication by $[2]^2$ proves the formula.
\end{proof}

The character formula gives
\[
P_{10}(u)=[2]^2[3]\mathbf1
 +u[2]^2(\rho+\sigma+\epsilon_t+\chi_{(1,2)}).
\]
Substitution into \eqref{eq:C3.modular.law} now yields
\[
P_{13}(u)=[2][3]^2\mathbf1+u[2]([3]+u)\epsilon_t+u^2[2](\sigma+2\rho+\chi_{(1,2)}).
\]

The characters $P_i(u)$ are listed in Table~\ref{fig C3}.

\begin{table}
\centering
\footnotesize
\setlength{\tabcolsep}{2.6pt}
\renewcommand{\arraystretch}{1.2}
\begin{tabular}{@{}c*{10}{c}@{}}
\toprule
$i$&$(3,\emptyset)$&$(21,\emptyset)$&$(\emptyset,3)$&$(2,1)$&$(1,2)$&
$(111,\emptyset)$&$(11,1)$&$(\emptyset,21)$&$(1,11)$&$(\emptyset,111)$\\
\midrule
$20$ & $[2][4][6]$ &  &  &  &  &  &  &  &  & \\
$19$ & $[2][4][5]$ &  & $u^2[2][4]$ &  &  &  &  &  &  & \\
$18$ & $[2][4]^2$ & $u^3[2]$ & $u^2[2][3]$ & $u^3[2]$ &  &  &  &  &  & \\
$17$ & $[2][3][4]$ &  & $u[2]^2[3]$ & $u^2[2]^2$ &  &  &  &  &  & \\
$16$ & $[2][3][4]$ & $u^2[2]^2$ & $u^2[2]^2$ & $u^2[2]^2$ &  &  &  &  &  & \\
$15$ & $[2]^2[4]$ & $u[2][3]$ & $u^2[2]$ & $u[2][3]$ &  &  &  &  &  & \\
$13$ & $[2][3]^2$ & $u^2[2]$ & $u[2]([3]+u)$ & $2u^2[2]$ & $u^2[2]$ &  &  &  &  & \\
$14$ & $[2][4]$ & $[2][4]$ &  & $[2][4]$ &  &  &  &  &  & \\
$12$ & $[2]^2[3]$ & $u[2]^2$ & $u[2]^2$ & $u[2]^2$ & $u[2]^2$ &  &  &  &  & \\
$10$ & $[2]^2[3]$ & $u[2]^2$ & $u[2]^2$ & $u[2]^2$ & $u[2]^2$ &  &  &  &  & \\
$11$ & $[2][3]$ &  & $[2][3]$ & $[2][3]$ & $[2][3]$ &  &  &  &  & \\
$9$ & $[2][3]$ & $[2][3]$ & $u[2]$ & $[2][3]$ & $u[2]$ &  &  & $u[2]$ &  & \\
$8$ & $[2]^3$ & $2u[2]$ & $u[2]$ & $2u[2]$ & $2u[2]$ &  & $u[2]$ &  &  & \\
$7$ & $[2]^2$ & $[2]^2$ &  & $[2]^2$ & $[2]^2$ &  & $[2]^2$ &  &  & \\
$6$ & $[2]^2$ & $[3]+2u$ & $u$ & $[3]+2u$ & $2u$ & $u$ & $2u$ & $u$ & $u$ & \\
$5$ & $[2]^2$ & $u$ & $[2]^2$ & $[3]+2u$ & $[3]+2u$ &  & $u$ & $u$ & $u$ & \\
$4$ & $[2]$ & $[2]$ & $[2]$ & $2[2]$ & $2[2]$ &  & $[2]$ & $[2]$ & $[2]$ & \\
$3$ & $[2]$ & $[2]$ & $[2]$ & $2[2]$ & $2[2]$ &  & $[2]$ & $[2]$ & $[2]$ & \\
$2$ & $[2]$ & $2[2]$ &  & $2[2]$ & $[2]$ & $[2]$ & $2[2]$ &  & $[2]$ & \\
$1$ & $1$ & $2$ & $1$ & $3$ & $3$ & $1$ & $3$ & $2$ & $3$ & $1$\\
\bottomrule
\end{tabular}
\par
\vspace{8pt}
\caption{Graded multiplicities of $\chi_{(\lambda,\mu)}$ in $P_i(u)$ in type $C_3$.
The ideals $I_i$ are listed in Table~\ref{tab:C3.ideals}.
Blank entries are zero.}
\label{fig C3}
\end{table}

\medskip

The table also determines when $\coh^{2k}(X_{H_i}(s))$ is a permutation representation.

\begin{corollary}\label{cor:C3.permutation}
Every graded piece of $\coh^*(X_{H_i}(s))$ is a permutation representation if and only if $i\ne17$.
For $i=17$, this fails exactly in degrees $4$, $6$, and $8$, and the total cohomology is not a permutation representation either.
\end{corollary}

\begin{proof}
For $i\in\{1,2,3,4,7,11,14,20\}$, the space $H_i$ is a parabolic subalgebra. If $H_i=\operatorname{Lie}P$, the flag-bundle argument in \eqref{eq:flagbundle} gives $P_i(u)=p_{P/B}(u)\Ind_{W_P}^{W}\mathbf1$, so every graded piece is a permutation representation.
For the remaining cases, set $\Pi_0=\mathbf1$ and $\Pi_j=\Ind_{K_j}^{W}\mathbf1$, where
\[
\begin{array}{c|l@{\qquad}c|l}
j&K_j&j&K_j\\\hline
1&\langle s_1,s_2,ts_2t\rangle
 &5&\langle s_1,t_2s_1t_2\rangle\\
2&\langle s_1,t,t_1,t_2\rangle
 &6&\langle s_1,t\rangle\\
3&\langle s_2,t\rangle
 &7&\langle t_1,t_2\rangle\\
4&\langle s_1,s_2\rangle
 &8&\langle s_1\rangle.
\end{array}
\]
Comparing these induced characters with Table~\ref{fig C3} gives Table~\ref{tab:C3.permutation}.
Every coefficient polynomial has nonnegative integer coefficients, proving the assertion for $i\ne17$.

For $i=17$, let $z=t_0t_1t_2\in W$.
This element changes the signs of all three coordinates, so $\epsilon_t(z)=-1$ and $\rho(z)=-3$.
Hence
\beq\label{eq:P17}
P_{17}(u)(z)=1+2u-u^2-4u^3-u^4+2u^5+u^6.
\eeq
Since permutation characters are nonnegative on every group element, the cohomology in degrees $4$, $6$, and $8$ cannot be a permutation representation. 
The remaining nonzero graded pieces have characters $\mathbf1$ in degrees $0$ and $12$, and $3\mathbf1+\epsilon_t=2\Pi_0+\Pi_1$ in degrees $2$ and $10$.

Finally, the total cohomology has character 
\[P_{17}(1)=24\mathbf1+12\epsilon_t+4\rho, \and P_{17}(1)(z)=0\]
by \eqref{eq:P17}.
If it were a permutation representation, its underlying $W$-set would have $48$ elements and $24$ orbits, since each orbit contributes one trivial constituent.
The vanishing at $z$ excludes singleton orbits, so all $24$ orbits would have exactly two points.
This contradicts the occurrence of the three-dimensional representation $\rho$.
Thus the total cohomology is not a permutation representation.
\end{proof}

\begin{table}[htbp]
\centering
\renewcommand{\arraystretch}{1.2}
\begin{tabular}{@{}cl@{}}
\toprule
$i$&$P_i(u)$\\
\midrule
$19$&$[2]^2[4](1+u^3)\Pi_0+u^2[2][4]\Pi_1$\\
$18$&$[2]^3(1+u^2+u^4)\Pi_0
       +u^2[2][3]\Pi_1+u^3[2]\Pi_3$\\
$16$&$[2]^3(1+u^3)\Pi_0+u^2[2]^2(\Pi_1+\Pi_3)$\\
$15$&$[2]^2(1+u^3)\Pi_0+u^2[2]\Pi_1+u[2][3]\Pi_3$\\
$13$&$[2]^2(1+u^3)\Pi_0
       +u[2][3]\Pi_1+u^2[2](\Pi_3+\Pi_4)$\\
$10,12$&$[2](1+u^3)\Pi_0+u[2]^2(\Pi_2+\Pi_4)$\\
$9$&$[4]\Pi_3+u[2]\Pi_5$\\
$8$&$(1+u^3)\Pi_0+u[2](\Pi_2+\Pi_4+\Pi_6)$\\
$6$&$(1+u^2)\Pi_3+u(\Pi_7+\Pi_8)$\\
$5$&$[3]\Pi_4+u\Pi_8$\\
\bottomrule
\end{tabular}
\par
\vspace{8pt}
\caption{Permutation-character expressions for nonparabolic
cases with $i\neq 17$.}
\label{tab:C3.permutation}
\end{table}

\begin{remark}\label{rem:C3.blowup}
Formula~\eqref{eq:C3.ex} can be obtained geometrically. 
Let $w=t_0t_1s_1$ and let $V=\C^6$ with symplectic form $\omega$.
The variety $X_{H_w}(s)$ consists of isotropic flags $F_1\subset F_2\subset F_3$ satisfying $sF_1\subset F_2$.
Forgetting $F_3$ gives a $\PP^1$-bundle over the blowup of the smooth quadric in $\PP(V)$
\[
Q_s=\{[v]\in\PP(V):\omega(sv,v)=0\} 
\]
at the six points corresponding to the eigenlines of $s$ (cf.~\cite{kiem-lee}).
These constructions extend over $\grs$, and the blowup and projective bundle formulas are compatible with monodromy.

Powers of the hyperplane class contribute $[5]\mathbf1$ to the graded character of $\coh^*(Q_s)$.
The remaining primitive middle cohomology is spanned by the difference of the classes of planes in the two families \cite[Theorem~1.13]{Reid72}.
Over the Grothendieck--Springer cover $\tgrs$, the parametrized $s$-stable Lagrangian subspace determines a plane on $Q_s$ in each fiber.
Its class and the square of the hyperplane class trivialize the middle cohomology local system on this cover, so its monodromy factors through $W$.
The generators $s_1,s_2$ preserve each family of planes, whereas $t$ replaces one eigenline by its symplectic partner.
The two corresponding planes meet in a line and hence belong to different families \cite[Theorem~1.2(b)]{Reid72}.
Thus the primitive middle cohomology has character $\epsilon_t$.

The permutation character on the six eigenlines is $\mathbf1+\rho+\sigma$.
The blowup formula and the $\PP^1$-bundle therefore give
$P_{H_w}(u)=[2]\big([5]\mathbf1+u^2\epsilon_t+u[3](\mathbf1+\rho+\sigma)\big)$,
which agrees with \eqref{eq:C3.ex}.
\end{remark}

\subsection{Computations in higher rank}
\label{ss:higher.rank}

To extend these computations to higher rank, we also use Precup--Sommers' decomposition of Hessenberg ideal fibers.

For a nilpotent element $z$, let $A(z)=Z_G(z)/Z_G(z)^\circ$ be its component group, and set
\[
G_{I,z}(u):=\sum_{j\geq0}[H^{2j}(\pi_I^{-1}(z))]u^j~\in(\Rep A(z))[u].
\]
For fixed $z$, as $I$ varies, the fibers $\pi_I^{-1}(z)$ admit stratifications into vector bundles over a common finite collection of smooth projective varieties.
Consequently, the characters $G_{I,z}(u)$ are linear combinations of the graded $A(z)$-characters of these varieties, with polynomial coefficients computable from root data \cite[Proposition~3.3 and Corollary~3.8]{PS}.
We treat these graded $A(z)$-characters as unknowns.

The cohomology of the Springer fiber at $z$ carries commuting actions of $W$ and $A(z)$.
Let $K$ be the matrix whose entry $K_{\psi,(z,\eta)}(u)$ records the graded multiplicity of $\psi\otimes\eta$ in this cohomology, where $\psi\in\Irr W$ and $(z,\eta)$ runs over the pairs occurring in the Springer correspondence.
With a suitable ordering, $K$ is triangular with powers of $u$ on the diagonal.
For $K$, we use the Springer convention of \cite[Section~4]{PS}, which assigns the trivial character to the regular nilpotent orbit.
This differs by a sign twist from our convention, which assigns the trivial character to the zero orbit.
The matrix $K$ is computable by the Lusztig--Shoji algorithm
\cite[Section~24]{Lus86v}.

For a Hessenberg space $H$, we have
\beq\label{eq:ideal.fiber.K}
u^{\dim(H/\fb)}[\eta:G_{H^\perp,z}(u)]=\sum_{\psi\in\Irr W}[\psi\otimes\operatorname{sgn}:P_H(u)]
K_{\psi,(z,\eta)}(u)
\eeq
for every such pair $(z,\eta)$ \cite[(4.7) and Proposition~5.2]{PS}.
Here $[\chi:-]$ denotes the multiplicity polynomial of $\chi$, and $\operatorname{sgn}$ is the sign character of $W$.
The multiplicities associated with pairs outside the Springer correspondence vanish \cite[Section~4]{PS}.

With $K$ taken as input, \eqref{eq:ideal.fiber.K} determines $G_{H_w^\perp,z}(u)$ from the values of $P_{H_w}(u)$ computed by our character formula for smooth $w$.
Substituting these values into Precup--Sommers' decompositions gives linear equations for the unknown graded $A(z)$-characters.  
These equations need not determine these graded $A(z)$-characters 
individually, but they may determine the linear combinations defining the characters $G_{I,z}(u)$ for additional ideals $I$.
Whenever $G_{I,z}(u)$ is determined for a further ideal $I$ and all nilpotent $z$, \eqref{eq:ideal.fiber.K} determines $P_{I^\perp}(u)$.
The modular law \eqref{eq:modular.law} can then be used to obtain further characters from known ones.

We have verified computationally that these methods determine $P_H(u)$ for every Hessenberg space $H$ in types $B_n$, $C_n$, and $D_n$ for $n\leq10$, and in all exceptional types.
We expect the same conclusion to hold in every Lie type and in arbitrary rank.

\begin{example}[Alternative complete computation in type $C_3$]
\label{ex:C3.ideal}
In type $C_n$, the codimension-one Hessenberg space is of the form $H_w$ for a smooth element $w$. 
In particular, in type $C_3$, the formula for $P_{19}$ follows from our character formula.
Below, we recover $P_{17}$ from the values of $P_7$ and $P_{14}$.
We keep the notation of Subsection~\ref{ss:C3} and write $G_{i,z}=G_{I_i,z}$.

Let $x$ have Jordan type $(2^2,1^2)$.
Then $A(x)\simeq\mathbb Z/2\mathbb Z$.
Denote its nontrivial character by $\varepsilon$.
The Precup--Sommers stratifications of $\pi_{I_i}^{-1}(x)$, for $i=7,14,17$, use two smooth projective varieties as bases.
Write $a(u)$ and $b(u)$ for the graded $A(x)$-characters of their cohomology.
For $I_{14}$, the strata over the first base have bundle ranks $0,1,2,3$.
For $I_7$, the corresponding ranks are $1,2,2,3$, together with one rank-zero stratum over the second base; for $I_{17}$, only the latter stratum occurs.
Thus \cite[(3.5) and Corollary~3.8]{PS} gives
\[
G_{14,x}=[4]a,\quad
G_{7,x}=u[2]^2a+b,\and
G_{17,x}=b.
\]
Our character formula determines $P_{14}$ and $P_7$.
Applying \eqref{eq:ideal.fiber.K} to these values yields
\[
G_{14,x}=[2][4](\mathbf1+\varepsilon)
\and
G_{7,x}=[2]^2[3]\mathbf1+u[2]^3\varepsilon.
\]
Comparing these, we obtain $a=[2](\mathbf1+\varepsilon)$ and $b=[2]^2\mathbf1$, and hence $G_{17,x}=[2]^2\mathbf1$.

Similarly, for $y$ of Jordan type $(2,1^4)$, we obtain
\[
G_{17,y}=[2]G_{14,y}=[2]^2[4]\mathbf1 \and G_{17,0}=p_{G/B}(u)=[2][4][6],
\]
and the fibers over all other nilpotent orbits are empty.
These characters for Hessenberg ideal fibers, together with the triangularity of $K$ \cite[Section~4]{PS}, imply that
\[
P_{17}(u)=c(u)\mathbf1+d(u)\epsilon_t+e(u)\rho.
\]
Using \eqref{eq:ideal.fiber.K} with $\dim(H_{17}/\fb)=6$,
the relevant entries of $K$ give
\[
\begin{aligned}
u^6[2]^2&=u^4e,\\
u^6[2]^2[4]&=u^6d+u^4e,\\
u^6[2][4][6]&=u^9c+u^6d+(u^4+u^6+u^8)e.
\end{aligned}
\]
Thus $c=[2][3][4]$, $d=u[2]^2[3]$, $e=u^2[2]^2$,  recovering the formula for $P_{17}$ in Proposition~\ref{prop:C3.codim3}.
\end{example}

\section{Positivity and unimodality conjectures}
\label{s:positivity}
Lesnevich asks whether the graded cohomology of regular semisimple Hessenberg varieties in types $B$ and $C$ is double-$h$ positive \cite[Question~A.4]{Lesnevich}.
He proves this positivity in degree two \cite[Theorem~8.6 and Table~2]{Lesnevich}.
Our computations support an affirmative answer and suggest that the coefficient polynomials are also unimodal.
We formulate positivity and unimodality conjectures for Hessenberg varieties and the intersection cohomology of regular semisimple Lusztig varieties in types $B$, $C$, and $D$, together with a positivity conjecture for their partial flag versions.
We first describe the character bases used in these statements.

\subsection{Types $B$ and $C$}\label{ss:positivity.BC}

Let $W_n=\{\pm1\}^n\rtimes S_n$ be the Weyl group of types $B_n$ and $C_n$.

For a partition $\lambda=(\lambda_1,\ldots,\lambda_\ell)$, write $|\lambda|=\sum_i\lambda_i$, and let $W_\lambda:=\prod_i W_{\lambda_i}\subset W_{\lvert\lambda\rvert}$ act on disjoint coordinate blocks of sizes $\lambda_i$.
Let $\delta_\lambda$ be the restriction to $W_\lambda$ of the character
\[W_{\lvert\lambda\rvert}\lra \{\pm 1\}, \quad ((\varepsilon_1,\ldots,\varepsilon_{|\lambda|}),\sigma)
\mapsto \prod_{j=1}^{|\lambda|}\varepsilon_j.\]

For partitions $\lambda$ and $\mu$ with $|\lambda|+|\mu|=n$, allowing either partition to be empty, set
\[
M_{\lambda,\mu}
:=
\Ind_{W_\lambda\times W_\mu}^{W_n}
(\mathbf1\otimes\delta_\mu).
\]
These characters form an integral basis of $\Rep W_n$. 
Under the Frobenius characteristic map for types $B$ and $C$, $M_{\lambda,\mu}$ corresponds to $h_\lambda(X)h_\mu(Y)$; see \cite[Appendix~A and Table~2]{Lesnevich}.
Here $X$ and $Y$ are two independent sets of variables, $h_\lambda=\prod_i h_{\lambda_i}$, and $h_m$ is the complete homogeneous symmetric function of degree $m$.

\subsection{Type $D$}\label{ss:positivity.D}
Let $n\geq4$, and identify $W_n^D:=W(D_n)$ with the kernel of the character $\delta_{(n)}$ of $W_n$.
We use a basis obtained by restricting the $W_n$-characters $M_{\lambda,\mu}$ to $W_n^D$ and splitting the restriction into two characters when $\lambda=\mu$.

For each unordered pair $\{\lambda,\mu\}$ of distinct partitions with $|\lambda|+|\mu|=n$, set
\[
M^D_{\{\lambda,\mu\}}
:=
\Res_{W_n^D}^{W_n}M_{\lambda,\mu}.
\]
This is independent of the order of $\lambda$ and $\mu$, since tensoring $M_{\lambda,\mu}$ with $\delta_{(n)}$ gives $M_{\mu,\lambda}$.

Now suppose  $\lambda=\mu$, so $n=2m$ with $\lambda$ a partition of $m$.
We have two characters $M_\lambda^\pm$ with
\beq\label{eq:D.Young.split}
\Res_{W_n^D}^{W_n}M_{\lambda,\lambda}
=
M_\lambda^++M_\lambda^-.
\eeq
To define them, set
\[
Q_\lambda
=(W_\lambda\times W_\lambda)\cap W_n^D
=\{(x,y)\in W_\lambda\times W_\lambda:
\delta_\lambda(x)\delta_\lambda(y)=1\}.
\]
Define a character $\theta_\lambda$ of $Q_\lambda$ by
\[
\theta_\lambda(x,y)=\delta_\lambda(x)=\delta_\lambda(y).
\]
Let $\tau=(1,m+1)\cdots(m,2m)\in W_n^D$ exchange the two blocks of $m$ coordinates.
Since conjugation by $\tau$ exchanges $x$ and $y$, $\theta_\lambda$ is invariant under this conjugation.
It therefore has two extensions $\theta_\lambda^\pm$ to $Q_\lambda\rtimes\langle\tau\rangle$, specified by $\theta_\lambda^\pm(\tau)=\pm1$.
Define
\[
M_\lambda^\pm
:=
\Ind_{Q_\lambda\rtimes\langle\tau\rangle}^{W_n^D}
\theta_\lambda^\pm.
\]
The decomposition \eqref{eq:D.Young.split} follows from Mackey's restriction formula \cite[Section~7.3, Proposition~22]{SerreReps} and transitivity of induction.

The characters $M^D_{\{\lambda,\mu\}}$ and $M_\lambda^\pm$ form an integral basis of $\Rep W_n^D$. 
Indeed, Young's rule \cite[Theorem~2.11.2]{Sagan}, together with the standard construction of the irreducible representations of $W_n^D$ \cite[Sections~5.5--5.6]{GP}, gives a unitriangular transition matrix to the irreducible-character basis.

\begin{remark}
This basis can be viewed as a type $D$ analogue of the double-$h$ basis.
To explain this, let $S^\lambda$ be the irreducible representation of $S_a$ indexed by a partition $\lambda$ of $a$, and 
\[
M^\lambda:=\Ind_{S_\lambda}^{S_a}\mathbf1
\]
be the permutation module associated with $S_\lambda=\prod_i S_{\lambda_i}$.
Their Frobenius characteristics are the Schur function $s_\lambda$ and the complete homogeneous symmetric function $h_\lambda$, respectively.

Regard these as representations of $W_a$ through the quotient $W_a\to S_a$.
For $a=|\lambda|$ and $b=|\mu|$, the irreducible $W_n$-representation indexed by the bipartition $(\lambda,\mu)$ is
\[
\Ind_{W_a\times W_b}^{W_n}
\bigl(S^\lambda\otimes
(\delta_{(b)}\otimes S^\mu)\bigr).
\]
See \cite[Chapter~I, Appendix~B, Section~9]{Macdonald}.
Replacing $S^\lambda$ and $S^\mu$ by $M^\lambda$ and $M^\mu$ gives a representation with character $M_{\lambda,\mu}$.

The same replacement gives the type-$D$ basis above.
For $\lambda\ne\mu$, restricting the two representations above to $W_n^D$ gives the irreducible representation indexed by $\{\lambda,\mu\}$ and a representation with character $M^D_{\{\lambda,\mu\}}$, respectively.
For $\lambda=\mu$, the two irreducible representations are constructed using the block-interchange actions $v\otimes w\mapsto\pm w\otimes v$ and induction (cf.~\cite[Chapter~I, Appendix~B, Section~8]{Macdonald} and \cite[Section~8.2, Proposition~25]{SerreReps}).
Replacing $S^\lambda$ by $M^\lambda$ in this construction gives representations with characters $M_\lambda^\pm$.
\end{remark}

\begin{remark}
In Shoji's symmetric-function description of characters of complex reflection groups, the bases above correspond to the specializations at $t=0$ of his $q$-function bases for $W_n=G(2,1,n)$ and $W_n^D=G(2,2,n)$, respectively; see \cite[Sections~2.15 and~3.1--3.3]{Shoji02}.
\end{remark}

\subsection{The conjectures}\label{ss:positivity.conjectures}

We call an element of $(\Rep W)[u]$ \emph{$M$-positive} if its coefficients in the corresponding character basis defined above belong to $\mathbb Z_{\geq0}[u]$.
We call it \emph{$M$-unimodal} if each coefficient polynomial has a unimodal coefficient sequence.

\begin{conjecture}\label{conj:classical.positivity}
Suppose that the root system of $G$ has type $B_n$, $C_n$, or $D_n$.
\begin{enumerate}
	\item For every Hessenberg space $H$, $P_H(u)$ is $M$-positive and $M$-unimodal.
	\item For every $w\in W$, $\widehat h(\widehat C_w)$ is $M$-positive and $M$-unimodal.
\end{enumerate}
\end{conjecture}

In types $B_n$ and $C_n$, the positivity assertion in part~(1) is the question posed by Lesnevich \cite[Question~A.4]{Lesnevich}.

\begin{example}\label{ex:M.logconcavity}
The unimodality assertion cannot be strengthened to log-concavity in either part of Conjecture~\ref{conj:classical.positivity}.
For a common counterexample, consider type $D_9$,
with $\alpha_i=e_i-e_{i+1}$ for $1\leq i\leq8$
and $\alpha_9=e_8+e_9$, and set
\[
H=\fb
\oplus\bigoplus_{i=1}^{3}\fg_{-\alpha_i}
\oplus\fg_{-(\alpha_1+\alpha_2)}
\oplus\fg_{-(\alpha_2+\alpha_3)}.
\]
The element $w=s_1s_2s_3s_1s_2$ is smooth and satisfies $H_w=H$.
Thus $\widehat h(\widehat C_w)=P_H(u)$ by
Theorem~\ref{thm:smooth}.
The coefficient of $M^D_{\{(1),(3,1^5)\}}$ in this
character is
\[
1+3u+10u^2+10u^3+3u^4+u^5.
\]
This polynomial is unimodal but not log-concave,
since $3^2<1\cdot10$.

For the Hessenberg space defined by the same five roots
in types $B_9$ and $C_9$, the same polynomial occurs
as the coefficient of $M_{(1),(3,1^5)}$.
\end{example}

We also propose a partial flag extension of the positivity assertion.
We use $P_{H,P}(u)$ as defined in Remark~\ref{rem:partial.Hess}.
For the analogous construction for Lusztig varieties, let $P\supset B$ be a parabolic subgroup and suppose that $w$ is the longest element of $W_PwW_P$.
For $\mathbf s\in\Grs$, set
\[
\overline Y_{w,P}(\mathbf s)
=
\{gP\in G/P:
g^{-1}\mathbf s g\in\overline{PwP}\}.
\]
Since $\overline{BwB}=\overline{PwP}$, the natural projection 
$\overline Y_w(\mathbf s)\to \overline Y_{w,P}(\mathbf s)$
is the base change of $G/B\to G/P$.
Smooth pullback for intersection complexes and the Leray--Hirsch argument give
\[
\IH^*(\overline Y_w(\mathbf s))
\simeq
\IH^*(\overline Y_{w,P}(\mathbf s))
\otimes\coh^*(P/B).
\]
This isomorphism is compatible with monodromy, which acts trivially on the second factor.
Thus Proposition~\ref{prop:mono.fact} implies that monodromy on $\IH^*(\overline Y_{w,P}(\mathbf s))$ factors through $W$.
By Corollary~\ref{cor:lus.formula}, its odd intersection cohomology vanishes, and
\[
\sum_{j\geq0}
\IH^{2j}(\overline Y_{w,P}(\mathbf s))u^j
=
\frac{\widehat h(\widehat C_w)}{p_{P/B}(u)}.
\]

\begin{conjecture}\label{conj:classical.partial.positivity}
Suppose that the root system of $G$ has type $B_n$, $C_n$, or $D_n$.
\begin{enumerate}
	\item For every Hessenberg space $H$ and every parabolic subgroup $P\supset B$ such that $H$ is $P$-stable, the graded character $P_{H,P}(u)$ is $M$-positive.

	\item For every parabolic subgroup $P\supset B$ and every $w\in W$ that is the longest element of $W_PwW_P$, the graded character $\widehat h(\widehat C_w)/p_{P/B}(u)$ is $M$-positive.
\end{enumerate}
\end{conjecture}

This extension does not include unimodality.
For example, for the codimension-one Hessenberg space in type $B_2$ and the stabilizer $P$ of $H^\perp$, Proposition~\ref{prop:codim1} gives
\[P_{H,P}(u)=(1+u^2)M_{(2),\emptyset}+uM_{(1,1),\emptyset},\]
which is not $M$-unimodal.
If $\alpha_1$ is long and $\alpha_2$ is short, then $w=s_1s_2s_1$ is smooth, satisfies $H_w=H$, and is the longest element of $W_PwW_P$.
Thus the same counterexample occurs for Lusztig varieties, since $\widehat h(\widehat C_w)/p_{P/B}(u)=P_{H,P}(u)$.

\begin{remark}
In type $A$, $M$-positivity and $M$-unimodality of $P_H(u)$ are the representation-theoretic form of the $e$-positivity and $e$-unimodality conjecture of Shareshian and Wachs \cite[Conjecture~1.3]{SW}, which refines the Stanley--Stembridge conjecture \cite{stanley-stembridge}.
The latter was proved by Hikita \cite{Hikita}, and independently by Griffin, Mellit, Romero, Weigl, and Wen \cite{GMRWW}.
A proof of the graded $e$-positivity conjecture has recently been announced in \cite{OpenAI}.
The corresponding question for Hessenberg varieties in $G/P$ was posed by Kiem and Lee \cite[Question~4.9]{kiem-lee}.

For Lusztig varieties in $\operatorname{GL}_n(\C)/B$, the corresponding positivity and unimodality assertions are equivalent to those in Haiman's conjecture \cite[Conjecture~2.1]{haiman}, with the geometric interpretation given by Abreu and Nigro \cite{AN}.
A different extension to other Weyl groups, using Springer representations, is proposed in \cite{trinh2026haimansconjecturespringersrepresentations}.
\end{remark}

For $G$ over $\C$, we have computationally verified part~(1) of both conjectures for $n\leq10$ and part~(2) of both conjectures for $n\leq7$.


\bibliographystyle{plain}
\bibliography{character}
\end{document}